\documentclass[12pt]{article}
\usepackage{amsmath}
\usepackage{graphicx}
\usepackage{amsfonts}
\usepackage{amssymb}
\newtheorem{theorem}{Theorem}[section]

\newtheorem{corollary}[theorem]{Corollary}

\newtheorem{definition}[theorem]{Definition}
\newtheorem{example}[theorem]{Example}

\newtheorem{lemma}[theorem]{Lemma}
\newtheorem{notation}[theorem]{Notation}

\newtheorem{proposition}[theorem]{Proposition}
\newtheorem{remark}[theorem]{Remark}

\newenvironment{proof}[1][Proof]{\textbf{#1.} }{\ \rule{0.5em}{0.5em}}

\begin{document}

\title{Hardy Algebras, $W^{\ast}$-Correspondences and Interpolation Theory}
\author{Paul S. Muhly\thanks{Supported in part by grants from the National Science
Foundation and from the U.S.-Israel Binational Science Foundation.}\\Department of Mathematics\\University of Iowa\\Iowa City, IA 52242\\e-mail: muhly@math.uiowa.edu
\and Baruch Solel\thanks{Supported in part by the U.S.-Israel Binational Science
Foundation and by the Fund for the Promotion of Research at the Technion.}\\Department of Mathematics\\Technion\\32000 Haifa, Israel\\e-mail: mabaruch@techunix.technion.ac.il}
\date{}
\maketitle
\begin{abstract}
Given a von Neumann algebra $M$ and a $W^{\ast}$-correspondence $E$ over $M$,
we construct an algebra $H^{\infty}(E)$ that we call the Hardy algebra of $E$.
When $M=\mathbb{C}=E$, then $H^{\infty}(E)$ is the classical Hardy space
$H^{\infty}(\mathbb{T})$ of bounded analytic functions on the unit disc. We
show that given any faithful normal representation $\sigma$ of $M$ on a
Hilbert space $H$ there is a natural correspondence $E^{\sigma}$ over the
commutant $\sigma(M)^{\prime}$, called the $\sigma$-dual of $E$, and that
$H^{\infty}(E)$ can be realized in terms of ($B(H)$-valued) functions on the
open unit ball $\mathbb{D}((E^{\sigma})^{\ast})$ in the space of adjoints of
elements in $E^{\sigma}$. We prove analogues of the Nevanlinna-Pick theorem in
this setting and discover other aspects of the value ``distribution theory''
for elements in $H^{\infty}(E)$. We also analyze the ``boundary behavior'' of
elements in $H^{\infty}(E)$ and obtain generalizations of the Sz.-Nagy--Foia\c
{s} functional calculus. The correspondence $E^{\sigma}$ has a dual that is
naturally isomorphic to $E$ and the commutants of certain, so-called induced
representations of $H^{\infty}(E)$ can be viewed as induced representations of
$H^{\infty}(E^{\sigma})$. For these induced representations a double commutant
theorem is proved.
\end{abstract}

\section{Introduction}

This study derives from our work on tensor algebras over $C^{\ast}%
$-correspondences \cite{MS98} in which we introduced a class of
non-self-adjoint operator algebras that contains a large variety of operator
algebras that have received \ substantial attention in the literature in
recent years, as well as other algebras of new, promising interest. Further,
the perspective of the theory gives a common approach to their analysis that
has its roots deeply embedded in the model theory of contraction operators
\cite{szNF70}.

Briefly, the setting of \cite{MS98} and its sequels \cite{MS99, MS00, MSp02}
is this. Let $A$ be a $C^{\ast}$-algebra. A $C^{\ast}$-correspondence over $A$
is a (right) Hilbert $C^{\ast}$-module \cite{cL94} that is made into a
bimodule over $A$ via a $C^{\ast}$-homomorphism $\varphi$ of $A$ into the
algebra $\mathcal{L}(E)$ consisting of all continuous adjointable maps on $E$.
Form the Fock space $\mathcal{F}(E)=A\oplus E\oplus E^{\otimes2}\oplus\cdots$,
which is the direct sum of all the internal tensor powers of $E$. The
\emph{tensor algebra }over $E$, $\mathcal{T}_{+}(E)$, is the norm closed
algebra in $\mathcal{L}(\mathcal{F}(E))$, generated by the creation operators
$\{T_{\xi}\mid\xi\in E\}$, where $T_{\xi}\eta:=\xi\otimes\eta$, $\eta
\in\mathcal{F}(E)$, and the operators $\{\varphi_{\infty}(a)\mid a\in A\}$,
where $\varphi_{\infty}$ is the natural extension of $\varphi$ to a
representation of $A$ in $\mathcal{L}(\mathcal{F}(E))$. When $A$ and $E$
reduce to the complex numbers $\mathbb{C}$, then $\mathcal{F}(E)$ is just
$\ell^{2}(\mathbb{N})$ and $\mathcal{T}_{+}(E)$ is the familiar disc algebra
$A(\mathbb{D})$ realized as analytic Toeplitz matrices on $\ell^{2}%
(\mathbb{N})$. On the other hand, if $A=\mathbb{C}$, but $E=\mathbb{C}^{n}$,
then $\mathcal{F}(E)$ is the full Fock space over $\mathbb{C}^{n}$, and
$\mathcal{T}_{+}(E)$ is the noncommutative disc algebra introduced by Popescu
\cite{gP96a}. There are many more examples and we describe some in the body of
this paper.

One aspect of the theory that is important for our considerations here is the
fact that the completely contractive representations of $\mathcal{T}_{+}(E)$
are determined by pairs $(T,\sigma)$, where $\sigma:A\rightarrow B(H)$ is a
$C^{\ast}$-representation of $A$ on a Hilbert space $H$ and where
$T:E\rightarrow B(H)$ is a completely contractive bimodule map. That is,
$T(a\xi b)=\sigma(a)T(\xi)\sigma(b)$, $a,b\in A$ and $\xi\in E$. The
``complete contractivity'' of $T$ is equivalent to the assertion that the
linear map $\tilde{T}$ defined initially on the balanced algebraic tensor
product $E\otimes H$ by the formula $\tilde{T}(\xi\otimes h):=T(\xi)h$ extends
to an operator of norm at most $1$ on the completion, denoted $E\otimes
_{\sigma}H$. The bimodule property of $T$, then, is equivalent to the equation%
\begin{equation}
\tilde{T}(\varphi(a)\otimes I)=\sigma(a)\tilde{T}\text{,} \label{intertwine0}%
\end{equation}
for all $a\in A$, which means that $\tilde{T}$ intertwines $\sigma$ and the
representation of $A$ on $E\otimes_{\sigma}H$ which is the composition of
Rieffel's induced representation $\sigma^{E}$ of $\mathcal{L}(E)$ and the
$C^{\ast}$-homomorphism $\varphi$ giving the left action of $A$ on $E$. That
is, we can (and shall) write equation (\ref{intertwine0}) as $\tilde{T}%
(\sigma^{E}\circ\varphi)=\sigma\tilde{T}$.

Thus we see that once $\sigma$ is fixed, the representations of $\mathcal{T}%
_{+}(E)$ are parameterized by the elements in the closed unit ball of the
intertwining space $\{\eta\in B(E\otimes_{\sigma}H,H)\mid\eta(\sigma^{E}%
\circ\varphi)=\sigma\eta$ and $\left\|  \eta\right\|  \leq1\}$. Reflecting on
this leads one ineluctably to the functional analyst's imperative: \emph{To
understand an algebra, view it as an algebra of functions on its space of
representations. }In our setting, then, we want to think about $\mathcal{T}%
_{+}(E)$ as a space of functions on this ball. For reasons that will be
revealed in a minute, we prefer to focus on the adjoints of the elements in
this space. Thus we let $E^{\sigma}=\{\eta\in B(H,E\otimes_{\sigma}H)\mid
\eta\sigma=(\sigma^{E}\circ\varphi)\eta\}$ and we write $\mathbb{D}%
((E^{\sigma})^{\ast})$ for the set $\{\eta\in B(E\otimes_{\sigma}H,H)\mid
\eta^{\ast}\in E^{\sigma}$, and $\left\|  \eta\right\|  <1\}$. That is,
$\mathbb{D}((E^{\sigma})^{\ast})$ is the norm-interior of the representation
space consisting of those $(T,\sigma)$ that are ``anchored by $\sigma$''. One
of our interests, then, is to understand the kind of functions that elements
$X$ of $\mathcal{T}_{+}(E)$ determine on $\mathbb{D}((E^{\sigma})^{\ast})$ via
the formula
\[
X(\eta^{\ast})=\sigma\times S(\eta)(X)\text{,}%
\]
where $\sigma\times S(\eta)$ is the so-called integrated form of the
representation determined by $\sigma$ and $\eta^{\ast}$.

In the special case when $A=E=\mathbb{C}$ and $\sigma$ is the one-dimensional
representation of $A$ on $\mathbb{C}$, we will see that $E^{\sigma}$ is also
one-dimensional, so $\mathbb{D}((E^{\sigma})^{\ast})$ is just the open unit
disc in the complex plane and, for $X\in\mathcal{T}_{+}(E)$, $X(\eta^{\ast})$
is the ordinary value of $X$ at the complex number $\bar{\eta}$. On the other
hand, if $A=E=\mathbb{C}$, but $\sigma$ is scalar multiplication on a Hilbert
space $H$ (the only possible representation of $\mathbb{C}$ on $H$), then
$\mathbb{D}((E^{\sigma})^{\ast})$ is the space of \emph{strict }contraction
operators on $H$ and for $\eta^{\ast}\in\overline{\mathbb{D}((E^{\sigma
})^{\ast})}^{\left\|  \cdot\right\|  }$ and $X\in\mathcal{T}_{+}%
(E)=A(\mathbb{D})$, $X(\eta^{\ast})$ is simply the value of $X$ at $\eta
^{\ast}$ defined through the Sz.-Nagy--Foia\c{s} functional calculus
\cite{szNF70}. For another example, if $A=\mathbb{C}$, but $E=\mathbb{C}^{n}$,
and if $\sigma$ is scalar multiplication on a Hilbert space $H$, then
$\mathbb{D}((E^{\sigma})^{\ast})$ is the space of row contractions on $H$,
$(T_{1},T_{2},\cdots,T_{n})$, of norm less than $1$; i.e. $\sum T_{i}^{\ast
}T_{i}\leq rI_{H}$ for some $r<1$. In this case, $X(\eta^{\ast})$ is given by
Popescu's functional calculus \cite{gP95a}.

In addition to paramaterizing certain representations of $\mathcal{T}_{+}(E)$,
$E^{\sigma}$ has another fundamental property: It is itself a $C^{\ast}%
$-correspondence - over \emph{the von Neumann algebra} $\sigma(A)^{\prime}$.
Indeed, it is not difficult to see that $E^{\sigma}$ becomes a bimodule over
$\sigma(A)^{\prime}$ via the formulae: $a\cdot\eta=(I_{E}\otimes a)\eta$ and
$\eta\cdot a=\eta a$, $\eta\in E^{\sigma}$, $a\in\sigma(A)^{\prime}$. Further,
if $\eta$ and $\zeta$ are in $E^{\sigma}$, then the product $\eta^{\ast}\zeta$
lies in the commutant $\sigma(A)^{\prime}$ and defines a $\sigma(A)^{\prime}%
$-valued inner product $\langle\eta,\zeta\rangle$ making $E^{\sigma}$ a
$C^{\ast}$-correspondence. In fact, since $E^{\sigma}$ is a weakly closed
space of operators, it has certain topological properties making it what we
call a $W^{\ast}$-correspondence \cite{MSp02}. It is because $E^{\sigma}$ is a
$W^{\ast}$-correspondence over $\sigma(A)^{\prime}$ that we focus on it, when
studying representations of $\mathcal{T}_{+}(E)$, rather than on its space of
adjoints. While $E^{\sigma}$ plays a fundamental role in our study of quantum
Markov processes \cite{MSp02}, its importance here - besides providing a space
on which to ``evaluate'' elements of $\mathcal{T}_{+}(E)$ - lies in the fact
that a certain natural representation of $E^{\sigma}$ generates the commutant
of the representation of $\mathcal{T}_{+}(E)$ obtained by ``inducing $\sigma$
up to'' $\mathcal{F}(E)$.

More explicitly, form the Hilbert space $K=\mathcal{F}(E)\otimes_{\sigma}H$.
On $K$ we have a natural covariant representation $(V,\rho)$ given by the
formulae $\rho(a):=\varphi_{\infty}(a)\otimes I$ and $V(\xi):=T_{\xi}\otimes
I$, $a\in A$, and $\xi\in E$. We call this representation $(V,\rho)$ the
\emph{induced }covariant representation of $(E,A)$ determined by $\sigma$. As
we elaborate in \cite{MS99}, such induced representations play a role in our
theory that is quite analogous to the role played by unilateral shifts in
model theory. We show here that there is a natural unitary operator $U$ from
$\mathcal{F}(E^{\sigma})\otimes_{\iota}H$ onto $\mathcal{F}(E)\otimes_{\sigma
}H$, where $\iota$ is the identity representation of $\sigma(A)^{\prime}$ on
$H$, such that if $(\tilde{V},\tilde{\rho})$ is induced covariant
representation of $(E^{\sigma},\sigma(A)^{\prime})$ on $\mathcal{F}(E^{\sigma
})\otimes_{\iota}H$ determined by $\iota$, then the ultraweak closure of
$U(\tilde{V}\times\tilde{\rho})(\mathcal{T}_{+}(E^{\sigma}))U^{\ast}$
coincides with the \emph{commutant} of $(V\times\rho)(\mathcal{T}_{+}(E))$.
This, of course, captures the well known fact that the commutatant of
$H^{\infty}(\mathbb{T})$ acting on $H^{2}(\mathbb{T})$ by multiplication is
$H^{\infty}(\mathbb{T})$. It also captures the vector-valued extensions that
play such an important role in model theory. When $E=\mathbb{C}^{n}$, then
this result gives theorems due to Popescu \cite{gP95} and Davidson and Pitts
\cite{DP99}.

It is primarily because of this commutant theorem that we cast our work in
this paper entirely in terms of $W^{\ast}$-correspondences. That is, we work
with von Neumann algebras $M$ and $W^{\ast}$-correspondences $E$ over them. We
still form the Fock space $\mathcal{F}(E)$ and the tensor algebra
$\mathcal{T}_{+}(E)$ over $E$, but because $\mathcal{F}(E)$ is a $W^{\ast}%
$-correspondence over $M$, the space $\mathcal{L}(\mathcal{F}(E))$ is a von
Neumann algebra viewed abstractly, i.e., it is a $W^{\ast}$-algebra. We call
the ultraweak closure of $\mathcal{T}_{+}(E)$ in $\mathcal{L}(\mathcal{F}(E))$
the \emph{Hardy algebra} of $E$ and denote it by $H^{\infty}(E)$. This is our
principal object of study. We will show that there are ``double commutant
theorems'' both at the level of correspondences and at the level of Hardy
algebras. That is, if we take a (faithful)\ normal representation $\sigma$ of
$M$ on a Hilbert space $H$ and form the dual $E^{\sigma}$, getting a $W^{\ast
}$-correspondence over $\sigma(M)^{\prime}$ and if we then form the dual of
$E^{\sigma}$ with respect to the identity representation of $\sigma
(M)^{\prime}$, we get a $W^{\ast}$-correspondence $E^{\sigma,\iota}$ over
$\sigma(M)^{\prime\prime}=\sigma(M)$ that is \emph{naturally }isomorphic to
$E$. (See Theorem \ref{duality}.) Further, $V\times\rho$ and $\tilde{V}%
\times\tilde{\rho}$ are ultraweak homeomorphisms and $(V\times\rho)(H^{\infty
}(E))^{\prime\prime}=(U(\tilde{V}\times\tilde{\rho})(\mathcal{T}_{+}%
(E^{\sigma}))U^{\ast})^{\prime}=(V\times\rho)(H^{\infty}(E))$. (See Corollary
\ref{doublecommutant})

As we will see in Theorem \ref{lem1.10Cor1.11}, given a faithful normal
representation $\sigma$ of $M$ on a Hilbert space $H$, we may evaluate
elements in $H^{\infty}(E)$ at points in $\mathbb{D}((E^{\sigma})^{\ast})$
also. That is, $H^{\infty}(E)$ may be viewed as functions on $\mathbb{D}%
((E^{\sigma})^{\ast})$, also. Further, when $H^{\infty}(E)$ is so represented,
one can study the ``value distribution theory'' of the space of functions. In
this context, we establish two capstone results from function theory: The
first (see Theorem \ref{Theorem3.1}) is a generalization of the
Nevanlinna-Pick interpolation theorem. It asserts that given two $k$-tuples of
operators in $B(H)$ (where $H$ is the representation space of $\sigma$),
$B_{1},B_{2},\cdots,B_{k}$, and $C_{1},C_{2},\cdots,C_{k}$, and given points
$\eta_{1},\eta_{2},\cdots,\eta_{k}$ in $\mathbb{D}((E^{\sigma})^{\ast})$, one
can find an element $X$ of norm at most one such that
\[
B_{i}X(\eta_{i}^{\ast})=C_{i}\text{,}%
\]
for all $i$ if and only if a certain matrix of \emph{maps}, which resembles
the classical Pick matrix, represents a completely positive operator. This
result captures numerous theorems in the literature that go under the name of
generalized Nevanlinna-Pick theorems. Our second capstone result is a
generalization of Schwartz's lemma (see Theorem \ref{Theorem3.5}). It asserts,
among several things, that if $X\in H^{\infty}(E)$ vanishes at the origin and
if $\left\|  X\right\|  \leq1$, then
\[
X(\eta^{\ast})X(\eta^{\ast})^{\ast}\leq\langle\eta,\eta\rangle\text{,}%
\]
where, recall, $\langle\cdot,\cdot\rangle$ is the $\sigma(M)^{\prime}$-valued
inner product on $E^{\sigma}$ defined above. We prove this as a consequence of
our Nevanlinna-Pick theorem.

These two theorems, then, are the main results of the paper. However, the
perspective that leads to them, i.e., viewing $H^{\infty}(E)$ as functions on
$\mathbb{D}((E^{\sigma})^{\ast})$ inspires a host of questions. One of these,
which we begin to address is the ``boundary behavior'' of functions in
$H^{\infty}(E)$. While it is true that given any point $\eta^{\ast}$ in the
boundary of $\mathbb{D}((E^{\sigma})^{\ast})$ it is possible to evaluate every
$X\in\mathcal{T}_{+}(E)$ at $\eta^{\ast}$, it may not be possible to evaluate
every $X\in H^{\infty}(E)$ at $\eta^{\ast}$. This, of course, is true in the
classical setting, where given any point on the unit circle one can find an
$H^{\infty}(\mathbb{T})$-function that cannot be evaluated there. However, in
our setting for \emph{certain }representations $\sigma$ of $M$ and for
\emph{certain }points $\eta^{\ast}$ on the boundary of $\mathbb{D}((E^{\sigma
})^{\ast})$, $X(\eta^{\ast})$ makes sense for every $X\in H^{\infty}(E)$ and,
in fact, the map $X\longrightarrow X(\eta^{\ast})$ is an ultraweakly
continuous, completely contractive, representation of $H^{\infty}(E)$ on $H$.
What we are describing may well be viewed as an extension of the
Sz.-Nagy--Foia\c{s} functional calculus \cite{szNF70}, as we hinted at above.
We are unable to identify all the points on the boundary of $\mathbb{D}%
((E^{\sigma})^{\ast})$ with this extension property, but we show that every
$\eta^{\ast}$ such that $(\eta,\sigma)$ is what we call completely
non-coisometric enjoys this feature. (See Theorem \ref{cncext}.) This notion
also has its roots in model theory and in the noncommutative multivariable
analysis of Popescu \cite{gP89}.

The next section is devoted to the basic properties of $W^{\ast}%
$-correspondences, the definition of $H^{\infty}(E)$, and other background
information that we use throughout the paper. Section \ref{Sect2} is devoted
to the duality theorem for $W^{\ast}$-correspondences, Theorem \ref{duality},
and the commutant theorems, Theorem \ref{commutant} and Corollary
\ref{doublecommutant}. In Section 4 we present a number of examples that
illustrate the duality theorem. Section 5 is devoted to the proofs of the
Nevanlinna-Pick theorem, Theorem \ref{Theorem3.1}, and to the Schwartz lemma,
Theorem \ref{Theorem3.5}. In Section 6, we present a number of applications
and illustrations of the Nevanlinna-Pick theorem. Finally, Section 7 is
devoted to showing how to extend elements of $H^{\infty}(E)$ to points
$\eta^{\ast}$ in the boundary of $\mathbb{D}((E^{\sigma})^{\ast})$ such that
$(\eta,\sigma)$ is completely non-coisometric. We show how to extend both the
Nevanlinna-Pick theorem and the Schwartz lemma to such points (Theorem
\ref{Theorem2.1} and Theorem \ref{Theorem2.2}).

\section{Preliminaries\label{Sect1}}

Throughout this paper, $M$ will be a von Neumann algebra. We will think of $M$
abstractly, i.e. as a $C^{\ast}$-algebra that is a dual space, without a
\emph{preferred} representation, unless otherwise specified. However, some of
the constructs we will associate with $M$ will depend very much on particular
representations and so we shall make these explicit in our analysis and
notation. Also, $E$ will denote a $W^{\ast}$-correspondence over $M$ in the
sense of \cite[Definition 2.2]{MSp02}. This means that $E$ is a self-dual
Hilbert module over $M$ in the sense of Paschke \cite{wP73} and that $E$ is a
left $M$-module through a normal $\ast$-homomorphism $\varphi:M\rightarrow
\mathcal{L}(E)$, where $\mathcal{L}(E)$ denotes the von Neumann algebra of all
bounded (necessarily adjointable) linear maps on $E$. $W^{\ast}$%
-correspondences are also called normal $M$-rigged $M$-modules in
\cite[Definition 5.1]{mR74b}.

Let $\sigma$ be a normal representation of $M$ on a Hilbert space $H$. We form
the so-called \emph{induced representation} \emph{space} $E\otimes_{\sigma}H$
built from $E$ and the representation $\sigma$. It is the Hausdorff completion
of the algebraic tensor product $E\otimes H$ in the positive semi-definite
sesquilinear form defined through the formula $\langle\xi\otimes h,\eta\otimes
k\rangle=\langle h,\sigma(\langle\xi,\eta\rangle)k\rangle$, $\xi\otimes
h,\;\eta\otimes k\in E\otimes H$. The representation of $\mathcal{L}(E)$
\emph{induced} by $\sigma$, denoted $\sigma^{E}$, is defined by the formula
$\sigma^{E}(T)(\xi\otimes h)=(T\xi)\otimes h$, $T\in\mathcal{L}(E)$; i.e.,
$\sigma^{E}(T)=T\otimes I$. Although the following lemma can be extracted from
Theorem 4.8 of \cite{mR74b}, we outline a proof for the sake of completeness.

\begin{lemma}
\label{Top}Let $E$ be a $W^{\ast}$-correspondence over the von Neumann algebra
$M$ and let $\sigma$ be a faithful representation of $M$ on $H$. Then the
induced representation $\sigma^{E}$ mapping $\mathcal{L}(E)$ to $B(E\otimes
_{\sigma}H)$ is a homeomorphism with respect to the ultraweak topologies.
\end{lemma}

\begin{proof}
For a bounded net $\{T_{\alpha}\}$ and an element $T$ in $\mathcal{L}(E)$,
$T_{\alpha}\rightarrow T$ ultraweakly if and only if for every $\xi,\eta\in E$
and every $g$ in the predual of $M$, $g(\langle\eta,T_{\alpha}\xi
\rangle)\rightarrow g(\langle\eta,T\xi\rangle)$ (see \cite[Proof of
Proposition 3.10]{wP73}). Since the net $\{\langle\eta,T_{\alpha}\xi\rangle\}$
is bounded for each fixed $\eta$ and $\xi$, this occurs if and only if for
every $\xi,\eta\in E$ and $h,k\in H$, $\langle\eta\otimes h,\sigma
^{E}(T_{\alpha})(\xi\otimes k)\rangle=\langle h,\langle\eta,T_{\alpha}%
\xi\rangle k\rangle\longrightarrow\langle h,\langle\eta,T\xi\rangle
k\rangle=\langle\eta\otimes h,\sigma^{E}(T)(\xi\otimes k)\rangle$.
\end{proof}

Observe that if $E$ is a $W^{\ast}$-correspondence over a von Neumann algebra
$M$, then each of the tensor powers of $E$, viewed as a $C^{\ast}%
$-correspondence over $M$ in the usual way, is in fact a $W^{\ast}%
$-correspondence over $M$ and so, too, is the full Fock space $\mathcal{F}%
(E)$, which is defined to be the direct sum $M\oplus E\oplus E^{\otimes
2}\oplus\cdots$, with its obvious structure as a right Hilbert module over $M$
and left action given by the map $\varphi_{\infty}$, defined by the formula
$\varphi_{\infty}(a):=diag\{a,\varphi(a),\varphi^{(2)}(a),\varphi
^{(3)}(a),\cdots\}$, where for all $n$, $\varphi^{(n)}(a)(\xi_{1}\otimes
\xi_{2}\otimes\cdots\xi_{n})=(\varphi(a)\xi_{1})\otimes\xi_{2}\otimes\cdots
\xi_{n}$, $\xi_{1}\otimes\xi_{2}\otimes\cdots\xi_{n}\in E^{\otimes n}$. The
\emph{tensor algebra} over $E$, denoted $\mathcal{T}_{+}(E)$, is defined to be
the norm-closed subalgebra of $\mathcal{L}(\mathcal{F}(E))$ generated by
$\varphi_{\infty}(M)$ and the \emph{creation operators} $T_{\xi}$, $\xi\in E$,
defined by the formula $T_{\xi}\eta=\xi\otimes\eta$, $\eta\in\mathcal{F}(E)$.
Although $\varphi$ need not be faithful, $\varphi_{\infty}$ always is, and
consequently it is not difficult to show that each $T_{\xi}$ satisfies
$\left\|  T_{\xi}\right\|  =\left\|  \xi\right\|  $, $\xi\in E$. We refer the
reader to \cite{MS98} for the basic facts about $\mathcal{T}_{+}(E)$.

\begin{definition}
\label{Hinfty}Given a $W^{\ast}$-correspondence $E$ over the von Neumann
algebra $M$, the ultraweak closure of the tensor algebra of $E$,
$\mathcal{T}_{+}(E)$, in $\mathcal{L}(\mathcal{F}(E))$, will be called the
\emph{Hardy Algebra of }$E$, and will be denoted $H^{\infty}(E)$.
\end{definition}

The following corollary is an immediate consequence of Lemma \ref{Top} and the
fact that the kernel of a normal representation of a von Neumann algebra is an
ultraweakly closed, central summand.

\begin{corollary}
Let $E$ be a $W^{\ast}$-correspondence over the von Neumann algebra $M$ and
let $\sigma$ be a normal representation of $M$ on a Hilbert space $H$. Then
the induced representation $\sigma^{\mathcal{F}(E)}:\mathcal{L}(\mathcal{F}%
(E))\rightarrow B(\mathcal{F}(E)\otimes_{\sigma}H)$, when restricted to
$H^{\infty}(E)$ gives an ultraweakly continuous representation of $H^{\infty
}(E)$ onto an ultraweakly closed subalgebra of $B(\mathcal{F}(E)\otimes
_{\sigma}H)$.
\end{corollary}

In most respects, the representation theory of $H^{\infty}(E)$ follows the
lines of the representation theory of $\mathcal{T}_{+}(E)$. However, there are
some differences that will be important here. To help illuminate these, we
need to review some of the basic ideas from \cite{MS98, MS99, MSp02}.

\begin{definition}
\label{Definition1.12}Let $E$ be a $W^{\ast}$-correspondence over a von
Neumann algebra $M$. Then:

\begin{enumerate}
\item  A \emph{completely contractive covariant representation }of $E$ on a
Hilbert space $H$ is a pair $(T,\sigma)$, where

\begin{enumerate}
\item $\sigma$ is a normal $\ast$-representation of $N$ in $B(H)$.

\item $T$ is a linear, completely contractive map from $E$ to $B(H)$ that is
continuous in the $\sigma$-topology of \cite{BDH88} on $E$ and the ultraweak
topology on $B(H).$

\item $T$ is a bimodule map in the sense that $T(S\xi R)=\sigma(S)T(\xi
)\sigma(R)$, $\xi\in E$, and $S,R\in M$.
\end{enumerate}

\item  A completely contractive covariant representation $(T,\sigma)$ of $E$
in $B(H)$ is called \emph{isometric }in case
\begin{equation}
T(\xi)^{\ast}T(\eta)=\sigma(\langle\xi,\eta\rangle)\text{,} \label{isometric}%
\end{equation}
for all $\xi,\eta\in E$.
\end{enumerate}
\end{definition}

It should be noted that the operator space structure of $E$ to which the
Definition \ref{Definition1.12} refers is that which $E$ inherits when viewed
as a subspace of its linking algebra. Also, we shall refer to an isometric,
completely contractive, covariant representation simply as an isometric
covariant representation. There is no problem with doing this because it is
easy to see that if one has a pair $(T,\sigma)$ satisfying all the conditions
of part 1 of Definition \ref{Definition1.12}, except possibly the complete
contractivity assumption, but which is isometric in the sense of equation
(\ref{isometric}), then necessarily $T$ is completely contractive. (See
\cite{MS98}.)

As we noted in the Introduction and developed in \cite[Lemmas 3.4--3.6]{MS98}
and in \cite{MSp02}, if a completely contractive covariant representation,
$(T,\sigma)$, of $E$ in $B(H)$ is given, then it determines a contraction
$\tilde{T}:E\otimes_{\sigma}H\rightarrow H$ defined by the formula $\tilde
{T}(\eta\otimes h):=T(\eta)h$, $\eta\otimes h\in E\otimes_{\sigma}H$. The
operator $\tilde{T}$ and the induced representation $\sigma^{E}$ are related
by the equation%
\begin{equation}
\tilde{T}\sigma^{E}\circ\varphi=\sigma\tilde{T}. \label{covariance}%
\end{equation}
In fact we have the following lemma from \cite[Lemma 2.16]{MSp02}.

\begin{lemma}
\label{CovRep}The map $(T,\sigma)\rightarrow\tilde{T}$ is a bijection between
all completely contractive covariant representations $(T,\sigma)$ of $E$ on
the Hilbert space $H$ and contractive operators $\tilde{T}:E\otimes_{\sigma
}H\rightarrow H$ that satisfy equation (\ref{covariance}). Given such a
$\tilde{T}$ satisfying this equation, $T$, defined by the formula
$T(\xi)h:=\tilde{T}(\xi\otimes h)$, together with $\sigma$ is a completely
contractive covariant representation of $E$ on $H$. Further, $(T,\sigma)$ is
isometric if and only if $\tilde{T}$ is an isometry.
\end{lemma}

\begin{remark}
\label{usefulnote}It is useful to note that this lemma shows that the
ultraweak continuity of $T$ really depends only on the fact that the
representation $\sigma$ is normal.
\end{remark}

The map $\Psi:\mathcal{L}(E)\rightarrow B(H)$ defined, then, by the formula
\[
\Psi(S):=\tilde{T}\sigma^{E}(S)\tilde{T}^{\ast}\text{,}%
\]
$S\in\mathcal{L}(E)$, evidently is completely positive, normal, and contractive.

\begin{definition}
\label{Definition1.12bis}If $(T,\sigma)$ is a completely contractive covariant
representation of $E$ in $B(H),$ the map $\Psi$ is called the \emph{completely
positive extension} of $(T,\sigma)$, and the representation $(T,\sigma)$ is
called \emph{fully coisometric} in case $\Psi(I_{E})=I_{H}$.
\end{definition}

The terminology is inspired by that used in the model theory of a single
contraction. A completely contractive covariant representation $(T,\sigma)$ is
isometric precisely when $\tilde{T}$ is an isometry. Likewise, it is fully
coisometric precisely when $\tilde{T}$ is a coisometry. The map $\Psi$ is a
normal $\ast$-representation precisely when $(T,\sigma)$ is isometric and it
is a unital $\ast$-representation precisely when $(T,\sigma)$ is both
isometric and fully coisometric. (We have, however, resisted the temptation to
call $(T,\sigma)$ unitary in this case.) We state one of the main results of
our theory for future reference. It is proved for $C^{\ast}$-correspondences
in \cite{MS98} and for $W^{\ast}$-correspondences as Theorem and Definition
2.18 of \cite{MSp02}.

\begin{theorem}
\label{Theorem 1.13}Let $E$ be a $W^{\ast}$-correspondence over a von Neumann
algebra $N$ and let $(T,\sigma)$ be a completely contractive covariant
representation of $E$ on the Hilbert space $H$. Then there is a Hilbert space
$K$ containing $H$ and an isometric covariant representation $(V,\rho)$ of $E$
on $K$ such that if $P$ is the projection of $K$ onto $H$, then

\begin{enumerate}
\item $P$ commutes with $\rho(N)$ and $\rho(A)P=\sigma(A)P,$ $A\in N$; and

\item  for all $\eta\in E$, $V(\eta)^{\ast}$ leaves $H$ invariant and
$PV(\eta)P=T(\eta)P$.
\end{enumerate}

The representation $(V,\rho)$ may be chosen so that the smallest subspace $K$
containing $H$ that reduces $(V,\rho)$ is $K$. When this is done, $(V,\rho)$
is unique up to unitary equivalence and is called the \emph{minimal isometric
dilation} of $(T,\sigma)$.

Further, if $(T,\sigma)$ is fully coisometric, the (unique minimal) isometric
dilation $(V,\rho)$ is fully coisometric, too.
\end{theorem}

The significance of condition 2. is that it implies that
\[
PV(\eta_{1})V(\eta_{2})\cdots V(\eta_{n})P=T(\eta_{1})T(\eta_{2})\cdots
T(\eta_{n})P
\]
for all $\eta_{i}\in E$. This, in turn, allows us to extend $T$, or more
accurately the pair $(T,\sigma)$, to a completely contractive representation
of $\mathcal{T}_{+}(E)$, as we shall discuss in a moment.

The minimal isometric dilation $(V,\rho)$ of $(T,\sigma)$ can be written down
explicitly in terms of a natural matrix associated with $\tilde{T}$. We will
need this expression in our analysis, and so we present the essentials here.
For details, consult \cite{MS98} and \cite{MSp02}; see Theorem 3.3 and
Corollary 5.21 in \cite{MS98} and the proof of Theorem 2.18 in \cite{MSp02},
in particular. Let $\Delta=(I-\tilde{T}^{\ast}\tilde{T})^{1/2}$ and let
$\mathcal{D}$ be its range. Then $\Delta$ is an operator on $E\otimes_{\sigma
}H$ and commutes with the representation $\sigma^{E}\circ\varphi$ of $M$, by
equation (\ref{covariance}). Write $\sigma_{1}$ for the restriction of
$\sigma^{E}\circ\varphi$ to $\mathcal{D}$. Form $K=H\oplus\mathcal{F}%
(E)\otimes_{\sigma_{1}}\mathcal{D}$, where the tensor product $\mathcal{F}%
(E)\otimes_{\sigma_{1}}\mathcal{D}$ is balanced over $\sigma_{1}$. The
representation $\rho$ is just $\sigma\oplus\sigma_{1}^{\mathcal{F}(E)}%
\circ\varphi_{\infty}$. For $V$, form $E\otimes_{\rho}K$ and define $\tilde
{V}:E\otimes_{\rho}K\rightarrow K$ matricially as
\begin{equation}
\tilde{V}:=\left[
\begin{array}
[c]{cccccc}%
\tilde{T} & 0 & 0 & \cdots &  & \\
\Delta & 0 & 0 &  & \ddots & \\
0 & I & 0 & \ddots &  & \\
0 & 0 & I & 0 & \ddots & \\
\vdots & 0 & 0 & I & \ddots & \\
&  &  &  & \ddots & \ddots
\end{array}
\right]  , \label{Vtilde}%
\end{equation}
where the identity operators in this matrix are interpreted as the operators
that identify $E\otimes_{\sigma_{n+1}}(E^{\otimes n}\otimes_{\sigma_{1}%
}\mathcal{D})$ with $E^{\otimes(n+1)}\otimes_{\sigma_{1}}\mathcal{D}$, and
where, in turn, $\sigma_{n+1}=\sigma_{1}^{E^{\otimes n}}\circ\varphi^{(n)}$.

It is easily checked that $\tilde{V}$ is an isometry and that the associated
covariant representation $(V,\rho)$ is an isometric dilation $(T,\sigma)$.
Moreover, it is easily checked that $(V,\rho)$ is minimal, i.e., that the
smallest subspace of $K$ containing $K$ and reducing $(V,\rho)$ is $K$.

The proof of the uniqueness of $(V,\rho)$ is the same as in the $C^{\ast}%
$-setting and is given in \cite[Proposition 3.2]{MS98}.

Finally, to see that $V$ is fully coisometric if $T$ is, observe that if $T$
is fully coisometric, then $\widetilde{T}$ is a coisometry as we noted
earlier. Thus $\widetilde{T}\Delta^{2}=0$. This implies that $\widetilde
{T}\Delta=0$. Therefore, from the form of $\widetilde{V}$, we see that
$\widetilde{V}\widetilde{V}^{\ast}=I$, which proves that $V$ is fully coisometric.

An important consequence of Theorem \ref{Theorem 1.13} is a fact that we shall
exploit over and over again: the bijective relationship between the completely
contractive covariant representations of $E$ and (completely contractive)
representations of $\mathcal{T}_{+}(E)$ that is proved in \cite[Theorem
3.10]{MS98} for the $C^{\ast}$-setting. We require a slight variant of it, and
so we present what we need here.

\begin{theorem}
\label{Theorem310MS98}Let $E$ be a $W^{\ast}$-correspondence over a von
Neumann algebra $M$. To every completely contractive covariant representation,
$(T,\sigma)$, of $E$ there is a unique completely contractive representation
$\rho$ of the tensor algebra $\mathcal{T}_{+}(E)$ that satisfies
\[
\rho(T_{\xi})=T(\xi)\text{,\ \ \ }\xi\in E\text{,}%
\]
and%
\[
\rho(\varphi_{\infty}(a))=\sigma(a)\text{,\ \ \ }a\in M\text{.}%
\]
The map $(T,\sigma)\mapsto\rho$ is a bijection between the set of all
completely contractive covariant representations of $E$ and all completely
contractive (algebra) representations of $\mathcal{T}_{+}(E)$ whose
restrictions to $\varphi_{\infty}(M)$ are continuous with respect to the
ultraweak topology on $\mathcal{L}(\mathcal{F}(E))$.
\end{theorem}

\begin{proof}
That every covariant representation $(T,\sigma)$ of $E$ gives rise to a
completely contractive algebra representation $\rho$ of $\mathcal{T}_{+}(E)$
is proved just as in \cite[Theorem 3.10]{MS98}, using Theorem \ref{Theorem
1.13}. The fact that the restriction of $\rho$ to $\varphi_{\infty}(M)$ is
continuous with respect to the ultraweak topology is evident from the fact
that $\varphi_{\infty}$ is a faithful normal representation of $M$ and is,
therefore, a homeomorphism (onto its range) with respect to the ultraweak
topologies involved. If, conversely, $\rho$ is a completely contractive
algebra representation of $\mathcal{T}_{+}(E)$ on a Hilbert space $H$, such
that $\sigma:=\rho\circ\varphi_{\infty}$ is an ultraweakly continuous
representation of $M$ on $H$ and if $T(\xi):=\rho(T_{\xi})$, then of course
$(T,\sigma)$ has all the properties displayed in Definition
\ref{Definition1.12} except, possibly, the property of continuity with respect
to the $\sigma$-topology of \cite{BDH88}. However, because $\sigma$ is normal,
this is automatic, as we noted in Remark \ref{usefulnote}.
\end{proof}

\begin{definition}
\label{integratedform}If $(T,\sigma)$ is a completely contractive covariant
representation of a $W^{\ast}$-correspondence $E$ over a von Neumann algebra
$M$, we call the representation $\rho$ of $\mathcal{T}_{+}(E)$ described in
Theorem \ref{Theorem310MS98} the \emph{integrated form} of $(T,\sigma)$ and
write $\rho=\sigma\times T$.
\end{definition}

\begin{remark}
\label{keyproblem}One of the principal difficulties one faces in dealing with
$\mathcal{T}_{+}(E)$ and $H^{\infty}(E)$ is to decide when the integrated
form, $\sigma\times T$, of a completely contractive covariant representation
$(T,\sigma)$ extends from $\mathcal{T}_{+}(E)$ to $H^{\infty}(E)$. This
problem arises already in the simplest situation, vis. when $M=\mathbb{C}=E$.
In this setting, $T$ is given by a single contraction operator on a Hilbert
space, $\mathcal{T}_{+}(E)$ ``is'' the disc algebra and $H^{\infty}(E)$ ``is''
the space of bounded analytic functions on the disc. The representation
$\sigma\times T$ extends from the disc algebra to $H^{\infty}(E)$ precisely
when there is no singular part to the spectral measure of the minimal unitary
dilation of $T$. There is no comparable global assertion known to us for the
general case, but we will present two sufficient conditions here, conditions
that are very similar to classical results. For the first, see Corollary
\ref{Corollary1.12} that we prove in a minute. The second will be proved in
Section \ref{CNCReps}, Theorem \ref{cncext}.
\end{remark}

A particular class of representations, the \emph{induced representations }of
\cite{MS99}, play a central role in our theory. Suppose $\pi_{0}$ is a normal
representation of the von Neumann algebra $M$ on the Hilbert space $H_{0}$.
Form $\mathcal{F}(E)\otimes_{\pi_{0}}H_{0}$ and the induced representation
$\pi_{0}^{\mathcal{F}(E)}$ of $\mathcal{L}(\mathcal{F}(E))$ on $\mathcal{F}%
(E)\otimes_{\pi_{0}}H_{0}$. Then the restrictions of $\pi_{0}^{\mathcal{F}%
(E)}$ to $\mathcal{T}_{+}(E)$ and $H^{\infty}(E)$ are called the\emph{ induced
representations} of these algebras determined by $\pi_{0}$. The associated
covariant representation of $E$, $(V,\pi)$, is given by the formulae $\pi
:=\pi_{0}^{\mathcal{F}(E)}\circ\varphi_{\infty}$ and $V(\xi):=\pi
_{0}^{\mathcal{F}(E)}(T_{\xi})$, and is called the \emph{induced covariant}
\emph{representation }of $E$. Note that $(V,\pi)$ is isometric. Evidently,
induced covariant representations are analogues of pure isometries. This point
is amplified in \cite{MS99}.

Suppose that $(T,\sigma)$ is an arbitrary completely contractive covariant
representation of the $W^{\ast}$-correspondence $E$ over the von Neumann
algebra $M$ and suppose $(T,\sigma)$ acts on the Hilbert space $H$. Although
it is not possible to form powers of $\tilde{T}$ because $\tilde{T}$ acts
between different spaces, it is possible to form a sequence of operators
$\{\tilde{T}_{n}\}_{n\geq0}$ that provide a serviceable substitute. For $n=0$,
$H$ is naturally isomorphic to $M\otimes_{\sigma}H$ via the map $a\otimes
h\rightarrow\sigma(a)h$ and so we identify $H$ with $M\otimes_{\sigma}H$ and
we set $\tilde{T}_{0}=I$. We set $\tilde{T}_{1}=\tilde{T}$. For $n>1$,
$\tilde{T}_{n}$ is defined inductively: $\tilde{T}_{n+1}:=\tilde{T}_{1}%
(I_{E}\otimes\tilde{T}_{n})$ - mapping $E^{\otimes n+1}\otimes_{\sigma}H$ to
$H$, where, as usual, we identify $E\otimes(E^{\otimes n}\otimes_{\sigma}H)$
with $E^{\otimes n+1}\otimes_{\sigma}H$. The sequence $\{\tilde{T}%
_{n}\}_{n\geq0}$ thus defined satisfies this analogue of the semigroup
property:%
\begin{equation}
\tilde{T}_{n+m}=\tilde{T}_{m}(I_{E^{\otimes m}}\otimes\tilde{T}_{n})=\tilde
{T}_{n}(I_{E^{\otimes n}}\otimes\tilde{T}_{m})\text{,} \label{iterates}%
\end{equation}
where $E^{\otimes m}\otimes(E^{\otimes n}\otimes_{\sigma}H)$ and $E^{\otimes
n}\otimes(E^{\otimes m}\otimes_{\sigma}H)$ are identified with $E^{\otimes
n+m}\otimes_{\sigma}H$ \cite[Section 2]{MS99}. We call $\tilde{T}_{n}$ the
$\emph{n}^{th}$\emph{-generalized power of }$\tilde{T}$.

\begin{remark}
\label{cpendo}We note for future reference the following important facts about
generalized powers. If $(T,\sigma)$ is a completely contractive covariant
representation of $E$ and $M$ on a Hilbert space $H$, then the map $\Theta$
defined on the \emph{commutant }of $\sigma(M)$ by the formula
\[
\Theta(a)=\tilde{T}(I_{E}\otimes a)\tilde{T}^{\ast}\text{,}%
\]
$a\in\sigma(M)^{\prime}$, is a completely positive map with range contained in
$\sigma(M)^{\prime}$. Further, $\Theta$ is an endomorphism of $\sigma
(M)^{\prime}$ if $\tilde{T}$ is an isometry, i.e., if $(T,\sigma)$ is
isometric. (The converse holds subject to a nondegeneracy condition that we
don't need here.) The $n^{th}$ powers of $\Theta$ are given by the formula%
\[
\Theta^{n}(a)=\tilde{T}_{n}(I_{E^{\otimes n}}\otimes a)\tilde{T}_{n}^{\ast
}\text{,}%
\]
$a\in\sigma(M)^{\prime}$, $n\geq0$. For the proofs of these assertions, see
Proposition 2.21 and Theorem 2.24 of \cite{MSp02} and Lemma 2.3 of \cite{MS99}.
\end{remark}

The following result is an analogue of the well known fact that the minimal
isometric dilation of strict contraction is an isometry \cite{szNF70}.

\begin{theorem}
\label{lem1.10Cor1.11}Let $(T,\sigma)$ be a completely contractive covariant
representation of the $W^{\ast}$-correspondence\ $E$ over the von Neumann
algebra $M$ with the property that $\left\|  \tilde{T}\right\|  <1$ and let
$(V,\pi)$ be its minimal isometric dilation. Then $(V,\pi)$ is an induced representation.
\end{theorem}

\begin{proof}
Let the Hilbert space of $(T,\sigma)$ be $H$ and let the Hilbert space of
$(V,\pi)$ be $K$. To show that $(V,\pi)$ is induced, we invoke Corollary 2.10
of \cite{MS99} and show that $\lim_{k\rightarrow\infty}\tilde{V}_{k}\tilde
{V}_{k}^{\ast}=0$ in the strong operator topology, where $\tilde{V}_{k}$ is
the $k^{th}$ generalized power of $\tilde{V}$. Note that each $\tilde{V}%
_{k}\tilde{V}_{k}^{\ast}$ is the projection onto the range of $\tilde{V}_{k}$
and that these form a decreasing sequence, by equation (\ref{iterates}). So it
suffices to show that $\lim_{k\rightarrow\infty}\tilde{V}_{k}\tilde{V}%
_{k}^{\ast}\xi=0$ for $\xi$ running over a spanning subset of $K$. Recall that
the operator $\tilde{V}:E\otimes K\rightarrow K$ associated with $V$ is given
by the matrix (\ref{Vtilde}). It is then a straightforward computation to show
that $\tilde{V}_{k}$ is given by the following matrix%
\[
\tilde{V}_{k}=\left[
\begin{array}
[c]{ccccccccc}%
\tilde{T}_{k} & 0 & 0 & \cdots &  &  &  &  & \\
\Delta(I_{E}\otimes\tilde{T}_{k-1}) & 0 & 0 & \cdots &  &  &  &  & \\
I_{E}\otimes\Delta(I_{E}\otimes\tilde{T}_{k-2}) & 0 & 0 & \cdots &  &  &  &  &
\\
I_{E^{\otimes2}}\otimes\Delta(I_{E}\otimes\tilde{T}_{k-3}) & 0 & 0 & \cdots &
&  &  &  & \\
\vdots & \vdots & \vdots & \vdots &  &  &  &  & \\
I_{E^{\otimes(k-1)}}\otimes\Delta & 0 & 0 & \cdots &  &  &  &  & \\
0 & I & 0 & \ddots &  &  &  &  & \\
0 & 0 & I & \ddots &  &  &  &  & \\
\vdots & \vdots & \vdots & \ddots &  &  &  &  &
\end{array}
\right]
\]
Thus, the $m^{th}$ column of $\tilde{V}_{k}\tilde{V}_{k}^{\ast}$ is%
\[
\left[
\begin{array}
[c]{c}%
\tilde{T}_{k}(I_{E^{\otimes(m-2)}}\otimes\Delta(I_{E}\otimes\tilde{T}%
_{k-m+1}))^{\ast}\\
\vdots\\
(I_{E^{\otimes(k-1)}}\otimes\Delta)(I_{E^{\otimes(m-2)}}\otimes\Delta
(I_{E}\otimes\tilde{T}_{k-m+1}))^{\ast}\\
0\\
0\\
\vdots
\end{array}
\right]
\]
when $2<m\leq k$. Thus the $m^{th}$ column of $\tilde{V}_{k}\tilde{V}%
_{k}^{\ast}$ has $k+1$ nonzero terms, each of norm less than or equal to
$\left\|  \tilde{T}\right\|  ^{k-m+1}$. Hence, for every $\xi\in E^{\otimes
m}\otimes H$, $\left\|  \tilde{V}_{k}\tilde{V}_{k}^{\ast}\xi\right\|
\leq(k+1)\left\|  \tilde{T}\right\|  ^{k-m+1}\left\|  \xi\right\|  $, provided
$k\geq m$. So for such $\xi$, $\lim_{k\rightarrow\infty}\tilde{V}_{k}\tilde
{V}_{k}^{\ast}\xi=0$. Since the vectors in $E^{\otimes m}\otimes H$ span $K$,
the proof is complete.
\end{proof}

\begin{corollary}
\label{Corollary1.12}Let $(T,\sigma)$ be a completely contractive covariant
representation of a $W^{\ast}$-correspondence $E$ over the von Neumann algebra
$M$.\ Suppose that $\left\|  \tilde{T}\right\|  <1$. Then the integrated form
of $(T,\sigma)$, $\sigma\times T$, extends from $\mathcal{T}_{+}(E)$ to an
ultraweakly continuous representation of $H^{\infty}(E)$.
\end{corollary}

\begin{proof}
From Theorem \ref{lem1.10Cor1.11} we know that the minimal isometric dilation
$(V,\pi)$ of $(T,\sigma)$ is an induced representation and so extends to a
representation of $H^{\infty}(E)$, i.e., we know that $\pi\times V=\pi
_{0}^{\mathcal{F}(E)}\mid H^{\infty}(E)$. However, $\sigma\times
T(\cdot)=P_{H}(\pi\times V(\cdot))P_{H}$, and so $\sigma\times T$ extends to
$H^{\infty}(E)$.
\end{proof}

We conclude this section by calling attention to some ``summability'' facts
that we will need in the sequel. Let $E$ be a $W^{\ast}$-correspondence over a
von Neumann algebra $M$. Then the projection $P_{n}$ onto the $n^{th}$ summand
of the Fock space $\mathcal{F}(E)$ belongs to $\mathcal{L}(\mathcal{F}(E))$,
of course. Further, we may form the one-parameter unitary group $\{W_{t}%
\}_{t\in\mathbb{R}}$ in $\mathcal{L}(\mathcal{F}(E))$ defined by the formula%
\begin{equation}
W_{t}:=\sum_{n=0}^{\infty}e^{int}P_{n}\text{.} \label{unitarygroup}%
\end{equation}
One must check that the series representing $W_{t}$ converges in the
weak-$\ast$ topology on $\mathcal{L}(\mathcal{F}(E))$. This is
straightforward, although a bit tedious. An alternate approach is to take a
faithful normal representation $\sigma$ of $M$ on a Hilbert space $H$, say,
and then form the faithful normal representation $\sigma^{\mathcal{F}(E)}$ of
$\mathcal{L}(\mathcal{F}(E))$ on $\mathcal{F}(E)\otimes_{\sigma}H$. Then the
series $\sum_{n\geq0}e^{int}\sigma^{\mathcal{F}(E)}(P_{n})$ clearly converges
in the ultraweak topology on $B(\mathcal{F}(E)\otimes_{\sigma}H)$ and
represents an element in the range of $\sigma^{\mathcal{F}(E)}$. Since
$\sigma^{\mathcal{F}(E)}$ is a homeomorphism between the ultraweak topology
restricted to its range and the weak-$\ast$ topology on $\mathcal{L}%
(\mathcal{F}(E))$, the series representing $W_{t}$ does, indeed, converge in
the weak-$\ast$ topology. Evidently, also, if we set $\gamma_{t}=AdW_{t}$,
then $\{\gamma_{t}\}_{t\in\mathbb{R}}$ is a weak-$\ast$ continuous action of
$\mathbb{R}$ on $\mathcal{L}(\mathcal{F}(E))$, called the \emph{gauge
automorphism group.} It normalizes $\mathcal{T}_{+}(E)$ and $H^{\infty}(E)$
since $\gamma_{t}(T_{\xi})=e^{-it}T_{\xi}$ for all $\xi\in E$ and all
$t\in\mathbb{R}$. Associated with $\{\gamma_{t}\}_{t\in\mathbb{R}}$ is the
sequence of ``Fourier coefficient operators'', $\{\Phi_{j}\}_{j\in\mathbb{Z}}%
$, defined by the formula%
\[
\Phi_{j}(a):=\frac{1}{2\pi}\int_{0}^{2\pi}e^{-int}\gamma_{t}(a)\,dt\text{,}%
\]
for all $a\in\mathcal{L}(\mathcal{F}(E))$, where the integral converges in the
weak-$\ast$ topology on $\mathcal{L}(\mathcal{F}(E))$. It is evident that each
$\Phi_{j}$ leaves invariant any weak-$\ast$ closed subspace of $\mathcal{L}%
(\mathcal{F}(E))$ that is left invariant by $\{\gamma_{t}\}_{t\in\mathbb{R}}$.
In particular, each $\Phi_{j}$ leaves $H^{\infty}(E)$ invariant. It is useful
to note and easy to show that for all $a\in\mathcal{L}(\mathcal{F}(E))$,%
\begin{equation}
\Phi_{j}(a)=\sum_{k}P_{k+j}aP_{k}\text{.} \label{shift}%
\end{equation}
This implies, in particular, that $\Phi_{j}(T_{\xi_{1}}T_{\xi_{2}}\cdots
T_{\xi_{n}})=T_{\xi_{1}}T_{\xi_{2}}\cdots T_{\xi_{n}}$, if $j=n$ and is zero
otherwise. Finally, we write $\Sigma_{k}$ for the ``$k^{th}$ arithmetic mean
operator'':
\[
\Sigma_{k}(a):=\sum_{|j|<k}(1-\frac{|j|}{k})\Phi_{j}(a)\text{,}%
\]
$a\in\mathcal{L}(\mathcal{F}(E))$. Straightforward estimates show that for all
$a\in\mathcal{L}(\mathcal{F}(E))$, $\lim_{k\rightarrow\infty}\Sigma_{k}(a)=a$
in the weak-$\ast$ topology. In fact, it is easy to see that $\lim
_{k\rightarrow\infty}\Sigma_{k}(a)=a$ in the ultrastrong topology of any
normal representation of $\mathcal{L}(\mathcal{F}(E))$.

\section{Duality of correspondences\label{Sect2}}

We continue with our basic setup and suppose that $E$ is a $W^{\ast}%
$-correspondence over $M$ and that $\sigma$ is a (normal) representation of
$M$ on the Hilbert space $H$. Observe that the commutant of $\sigma
^{E}(\mathcal{L}(E))$ is $\{I\otimes X\mid X\in\sigma(M)^{\prime}\}$ by
Theorem 6.23 of \cite{mR74}, where $I\otimes X(\xi\otimes h):=\xi\otimes Xh$,
as the notation suggests. The point is that $I\otimes X$ makes sense as a
bounded operator on $E\otimes_{\sigma}H$ only when $X\in\sigma(M)^{\prime}$.
Evidently, the map $X\rightarrow I\otimes X$ is a normal representation of
$\sigma(M)^{\prime}$ on $E\otimes_{\sigma}H$.

\begin{definition}
Let $\sigma:M\rightarrow B(H)$be a normal representation of the von Neumann
algebra $M$ on the Hilbert space $H$. Then for a $W^{\ast}$-correspondence $E$
over $M$, the $\sigma$\emph{-dual} of $E$, denoted $E^{\sigma}$, is defined to
be
\[
\{T\in B(H,E\otimes_{\sigma}H)\mid T\sigma(a)=\sigma^{E}\circ\varphi
(a)T\text{, }a\in M\}\text{.}%
\]
\end{definition}

\begin{proposition}
With respect to the actions of $\sigma(M)^{\prime}$ and the $\sigma
(M)^{\prime}$-valued inner product defined as follows, $E^{\sigma}$ becomes a
$W^{\ast}$-correspondence over $\sigma(M)^{\prime}$: For $X,Y\in
\sigma(M)^{\prime}$, and $T\in E^{\sigma}$, $X\cdot T\cdot Y:=(I\otimes X)TY$,
and for $T,S\in E^{\sigma}$, $\langle T,S\rangle_{\sigma(M)^{\prime}}%
:=T^{\ast}S$.
\end{proposition}

\begin{proof}
It is straightforward to verify that $E^{\sigma}$ is bimodule over
$\sigma(M)^{\prime}$ and that $\langle\cdot,\cdot\rangle_{\sigma(M)^{\prime}}$
is a bona fide $\sigma(M)^{\prime}$-valued inner product. The fact that
$E^{\sigma}$ is a $W^{\ast}$-correspondence is immediate from Proposition 2.3
of \cite{MSp02}.
\end{proof}

\begin{notation}
When we wish to emphasize module operations over the composition of operators,
we will write the module operations using a `cdot', $\cdot$, as we did in the
preceding proposition. Also, when operator composition needs to be emphasized,
we will denote it, as usual, by a `%
$\backslash$%
circ', $\circ$.
\end{notation}

We have met $E^{\sigma}$ before. By virtue of equation (\ref{covariance}) and
Lemma \ref{CovRep}, we see that the norm-closed unit ball of $(E^{\sigma
})^{\ast}$ - the space of adjoints of the operators in $E^{\sigma}$ -
parameterizes the completely contractive covariant representations of $E$.
More precisely, Lemma \ref{CovRep} can be rephrased as the following theorem,
which we want to record for future reference.

\begin{theorem}
\label{CovRepbis}Let $\sigma:M\rightarrow B(H)$ be a normal representation of
$M$ on the Hilbert space $H$. For $S$ in the closed unit ball of $(E^{\sigma
})^{\ast}$, define $T_{S}:E\rightarrow B(H)$ by $T_{S}(\xi)h=S(\xi\otimes h)$,
$\xi\otimes h\in E\otimes_{\sigma}H$. Then the map $S\rightarrow(T_{S}%
,\sigma)$ is a bijection between the closed unit ball of $(E^{\sigma})^{\ast}$
and the set of all completely contractive covariant representations of $E$
such that the associated representation of $M$ is $\sigma$.
\end{theorem}

Let $\iota$ denote the identity representation of $\sigma(M)^{\prime}$ on $H$.
Then we can form $(E^{\sigma})^{\iota}$, which we shall denote more simply by
$E^{\sigma,\iota}$. More explicitly, $E^{\sigma,\iota}=\{T\in B(H,E^{\sigma
}\otimes_{\iota}H)\mid T\iota(a)=\iota^{E^{\sigma}}\circ\varphi_{E^{\sigma}%
}(a),\;a\in\sigma(M)^{\prime}\}$, where, of course, $\varphi_{E^{\sigma}}$
denotes the left action of $\sigma(M)^{\prime}$ on $E^{\sigma}$. However,
since $\varphi_{E^{\sigma}}$ is given by the formula $\varphi_{E^{\sigma}%
}(a)\xi\otimes h=\xi\otimes\iota(a)h=\xi\otimes ah=(I_{E^{\sigma}}\otimes
a)(\xi\otimes h)$, $a\in\sigma(M)^{\prime}$, $E^{\sigma,\iota}=\{T\in
B(H,E^{\sigma}\otimes_{\iota}H)\mid Ta=(I_{E^{\sigma}}\otimes a)T,\;a\in
\sigma(M)^{\prime}\}$. This is a correspondence over $\sigma(M)^{\prime\prime
}=\sigma(M)$. In fact, it is naturally isomorphic to $E$, assuming $\sigma$ is
faithful. This result, which is fundamental for our analysis, is a consequence
of the following lemma that generalizes Lemma 2.10 of \cite{MSp02}.

\begin{lemma}
\label{Lemma2.10MS02}Let $\pi$ and $\rho$ be two normal $\ast$-representations
of the von Neumann algebra $M$ on Hilbert spaces $H$ and $K$, respectively.
Assume that $\pi$ is injective and let $\mathcal{L}_{M}(H,K)$ be the space of
intertwining operators. Then
\[
\bigvee\{X(H):X\in\mathcal{L}_{M}(H,K)\}=K.
\]
\end{lemma}

\begin{proof}
Write $\mathcal{N}$ for the von Neumann algebra
\[
\mathcal{N}=\{\left(
\begin{array}
[c]{cc}%
\rho(a) & 0\\
0 & \pi(a)
\end{array}
\right)  :\;a\in M\}\subset B(K\oplus H).
\]
Then
\[
\mathcal{N}^{\prime}=\left(
\begin{array}
[c]{cc}%
\rho(M)^{\prime} & \mathcal{L}_{M}(H,K)\\
\mathcal{L}_{M}(K,H) & \pi(M)^{\prime}%
\end{array}
\right)  .
\]
Write $P_{H}$ and $P_{K}$ for the projections onto $H$ and $K$. Both are
projections in $\mathcal{N}^{\prime}$. The central support, $C_{H}$, of
$P_{H}$ is a projection in $\mathcal{N}$ and, so, there is a projection $q$ in
$M$ such that
\[
C_{H}=\left(
\begin{array}
[c]{cc}%
\rho(q) & 0\\
0 & \pi(q)
\end{array}
\right)  .
\]
It follows, then, that $0=(I-C_{H})P_{H}=\left(
\begin{array}
[c]{cc}%
0 & 0\\
0 & \pi(I-q)
\end{array}
\right)  $. Thus $q=I$, since $\pi$ is injective. We conclude that $C_{H}=I$.
Consider the projection onto the subspace $[\mathcal{N}^{\prime}(H)]$. It lies
in the center of $\mathcal{N}$ and dominates $P_{H}$. Hence it is $I$.
Therefore $K=[P_{K}\mathcal{N}^{\prime}(H)]=\bigvee\{X(H):X\in\mathcal{L}%
_{M}(H,K)\}$, which was to be shown.
\end{proof}

To define the isomorphism $W:E\rightarrow E^{\sigma,\iota}$ that we want to
use, observe that $E^{\sigma,\iota}$ consists of bounded linear operators $S$,
say, from $H$ to $E^{\sigma}\otimes_{\iota}H$ such that $Sa=(I_{E^{\sigma}%
}\otimes a)S$ for all $a\in\sigma(M)^{\prime}$. For $\xi\in E$, we define
$W(\xi)\in E^{\sigma,\iota}$, by defining the adjoint $W(\xi)^{\ast}$, which
is a map from $E^{\sigma}\otimes_{\iota}H$ to $H$. For this purpose, we note
that $\xi\in E$ defines a bounded linear map $L_{\xi}$ from $H$ to
$E\otimes_{\sigma}H$ via the formula%
\[
L_{\xi}(h)=\xi\otimes h\text{.}%
\]
Indeed, $\left\|  L_{\xi}h\right\|  ^{2}=\langle\xi\otimes h,\xi\otimes
h\rangle=\langle h,\sigma(\langle\xi,\xi\rangle)h\rangle\leq\left\|
\xi\right\|  ^{2}\left\|  h\right\|  ^{2}$. Then, $W(\xi)^{\ast}$ is defined
by the formula
\begin{equation}
W(\xi)^{\ast}(T\otimes h)=L_{\xi}^{\ast}(Th), \label{wxs}%
\end{equation}
$T\otimes h\in E^{\sigma}\otimes H$.

\begin{theorem}
\label{duality}If the representation $\sigma$ of $M$ on $H$ is faithful, then
the map $W:E\rightarrow E^{\sigma,\iota}$ just defined is an isomorphism of
$W^{\ast}$-correspondences.
\end{theorem}

\begin{proof}
If $\xi\in E$, $h_{i}\in H$, and $T_{i}\in E^{\sigma}$, $i=1,2,\cdots,n$,
then
\begin{align*}
\left\|  \sum_{i=1}^{n}L_{\xi}^{\ast}(T_{i}h_{i})\right\|  ^{2}\leq\left\|
\xi\right\|  ^{2}\left\|  \sum_{i=1}^{n}T_{i}h_{i}\right\|  ^{2}=\left\|
\xi\right\|  ^{2}\sum_{i,j=1}^{n}\langle T_{i}h_{i},T_{j}h_{j}\rangle\\
=\left\|  \xi\right\|  ^{2}\sum_{i,j=1}^{n}\langle h_{i},\langle T_{i}%
,T_{j}\rangle h_{j}\rangle=\left\|  \xi\right\|  ^{2}\left\|  \sum_{i=1}%
^{n}T_{i}\otimes h_{i}\right\|  ^{2},
\end{align*}
which proves that $W(\xi)$ is a bounded linear map with $\left\|
W(\xi)\right\|  \leq\left\|  \xi\right\|  $. Also, for $a\in\sigma(M)^{\prime
}$, $T\in E^{\sigma}$, and $h\in H$, the equation $W(\xi)^{\ast}(a\otimes
I)\cdot(T\otimes h)=W(\xi)^{\ast}(a\cdot T\otimes h)=L_{\xi}^{\ast}((a\cdot
T)h)=L_{\xi}^{\ast}((I\otimes a)(Th))=a(L_{\xi}^{\ast}Th)=aW(\xi)^{\ast
}(T\otimes h)$, shows that $W(\xi)$ lies in $E^{\sigma,\iota}$.

Write $\mu$ for the map from $E^{\sigma}\otimes_{\iota}H\rightarrow
E\otimes_{\sigma}H$ defined by the formula $\mu(T\otimes h)=Th$. Then
$\left\|  \mu(\sum T_{i}\otimes h_{i})\right\|  ^{2}=\sum_{i,j}\langle
T_{i}h_{i},T_{j}h_{j}\rangle=\sum_{i,j}\langle h_{i},T_{i}^{\ast}T_{j}%
h_{j}\rangle=\left\|  \sum T_{i}\otimes h_{i}\right\|  ^{2}$. Thus $\mu$ is a
well-defined isometry from $E^{\sigma}\otimes_{\iota(\sigma(M)^{\prime})}H$ to
$E\otimes_{\sigma}H$. However, Lemma \ref{Lemma2.10MS02} shows that
$E\otimes_{\sigma}H=\bigvee\{Th\mid T\in E^{\sigma}$, $h\in H\}$ and this
means that $\mu$ maps onto $E\otimes_{\sigma}H$. Thus, $\mu$ is a Hilbert
space isomorphism.

Now observe that $W(\xi)^{\ast}=L_{\xi}^{\ast}\circ\mu$. Consequently,
$W(\xi)=\mu^{\ast}\circ L_{\xi}$ and for $\xi_{1},\xi_{2}\in E$, we have
$W(\xi_{1})^{\ast}W(\xi_{2})=(L_{\xi_{1}}^{\ast}\circ\mu)(\mu^{\ast}\circ
L_{\xi_{2}})=L_{\xi_{1}}^{\ast}L_{\xi_{2}}=\sigma(\langle\xi_{1},\xi
_{2}\rangle)$. Thus, $W$ is an isometry from $E$ into $E^{\sigma,\iota}$. On
the other hand, $W$ is a bimodule map because for $a,b\in M$, $\xi\in E$ and
$T\otimes h\in E^{\sigma}\otimes_{\iota(\sigma(M)^{\prime})}H$, $W(a\cdot
\xi\cdot b)^{\ast}(T\otimes h)=L_{a\cdot\xi\cdot b}^{\ast}(Th)$, by definition
and for $\eta\otimes k\in E\otimes_{\sigma}H$, $L_{a\cdot\xi\cdot b}^{\ast
}(\eta\otimes k)=\sigma(\langle a\cdot\xi\cdot b,\eta\rangle)k=\sigma
(b)^{\ast}\sigma(\langle\xi,a^{\ast}\cdot\eta\rangle)k=\sigma(b)^{\ast}L_{\xi
}^{\ast}(a^{\ast}\otimes I)(\eta\otimes k)$. Thus $W(a\cdot\xi\cdot b)^{\ast
}(T\otimes h)=L_{a\cdot\xi\cdot b}^{\ast}(Th)=\sigma(b)^{\ast}L_{\xi}^{\ast
}(a^{\ast}\otimes I)(Th)=\sigma(b)^{\ast}L_{\xi}^{\ast}(T\sigma(a^{\ast
})h)=\sigma(b)^{\ast}W(\xi)^{\ast}(I\otimes\sigma(a^{\ast}))(T\otimes h)$.
Hence $W(a\cdot\xi\cdot b)=I\otimes\sigma(a)W(\xi)\sigma(b)$, as we wanted to show.

It remains to show that $W$ is surjective. To this end, fix $S\in
E^{\sigma,\iota}$. Then $S$ is a bounded linear operator from $H$ to
$E^{\sigma}\otimes_{\iota}H$ such that $Sa=(I_{E^{\sigma}}\otimes a)S$ for all
$a\in\sigma(M)^{\prime}$. We set $F:=\mu\circ S$, obtaining a bounded linear
transformation from $H$ to $E\otimes_{\sigma}H$. We claim that $F$ satisfies
the equation $Fa=(I_{E}\otimes a)F$ for all $a\in\sigma(M)^{\prime}$. Indeed,
if $T\otimes h$ lies in $E^{\sigma}\otimes_{\iota}H$, then $\mu(a\otimes
I)(T\otimes h)=\mu(a\cdot T\otimes h)=\mu(((I\otimes a)T)\otimes h)=(I\otimes
a)Th=(I\otimes a)\mu(T\otimes h)$. Hence, $Fa=\mu Sa=\mu(I_{E^{\sigma}}\otimes
a)S=(I_{E^{\sigma}}\otimes a)\mu S=(I_{E^{\sigma}}\otimes a)F$, as we claimed.
Next, we define a $B(H)$-valued function $B$ on $E$ by the formula
$B(\xi):=F^{\ast}\circ L_{\xi}$. Then, for $a\in\sigma(M)^{\prime}$,
$aB(\xi)=aF^{\ast}L_{\xi}=F^{\ast}(I_{E^{\sigma}}\otimes a)L_{\xi}=F^{\ast
}L_{\xi}a=B(\xi)a$. Consequently, $B$ takes values in $\sigma(M)$. Further,
for $b\in M$, $B(\xi b)=F^{\ast}L_{\xi b}=F^{\ast}L_{\xi}\sigma(b)=B(\xi
)\sigma(b)$. Thus, $B$ is a continuous, right module map from $E$ into
$\sigma(M)$. Since $E$ is self-dual and since we may identify $M$ with
$\sigma(M)$, because $\sigma$ is assumed to be faithful, we conclude that
there is an $\eta\in E$ such that $B(\xi)=\langle\eta,\xi\rangle$, for all
$\xi\in E$. Thus $F^{\ast}L_{\xi}=L_{\eta}^{\ast}L_{\xi}$, for all $\xi\in E$.
Consequently, $F=L_{\eta}$, so $S=\mu^{\ast}F=\mu^{\ast}L_{\eta}=W(\eta)$.
\end{proof}

Our next goal is to show that if $\sigma:M\rightarrow B(H)$ is a faithful
normal representation of $M$ on a Hilbert space $H$, then there is a canonical
isometric covariant representation of $E^{\sigma}$ on $\mathcal{F}%
(E)\otimes_{\sigma}H$ that generates the commutant of $\sigma^{\mathcal{F}%
(E)}(H^{\infty}(E))$. So let $\sigma:M\rightarrow B(H)$ be fixed and define
$\pi:\sigma(M)^{\prime}\rightarrow\mathcal{F}(E)\otimes_{\sigma}H$ by the
formula $\pi(a)=I_{\mathcal{F}(E)}\otimes a$, $a\in\sigma(M)^{\prime}$. As we
have mentioned before, $\pi$ is a bona fide normal representation of
$\sigma(M)^{\prime}$ mapping onto $\sigma^{\mathcal{F}(E)}(\mathcal{L}%
(\mathcal{F}(E)))^{\prime}$, by Theorem 6.23 of \cite{mR74}. To define the map
$\Psi:E^{\sigma}\rightarrow\mathcal{F}(E)\otimes_{\sigma}H$, note that since
$E^{\sigma}=\{T:H\rightarrow E\otimes_{\sigma}H\mid\sigma^{E}\circ
\varphi(a)T=T\sigma(a),\;a\in M\}$, each $T\in E^{\sigma}$ defines an operator
$T^{(k)}:E^{\otimes n}\otimes_{\sigma}H\rightarrow E^{\otimes n+1}%
\otimes_{\sigma}H$ via the formula $T^{(k)}(\xi\otimes h)=\xi\otimes Th$,
where we have identified $E^{\otimes n+1}\otimes_{\sigma}H$ with $E^{\otimes
n}\otimes_{\sigma^{E}\circ\varphi}(E\otimes_{\sigma}H)$ in the obvious
fashion. The norms of all the $T^{(k)}$ are all bounded by $\left\|
T\right\|  $ and so we may define $\Psi(T):\mathcal{F}(E)\otimes_{\sigma
}H\rightarrow\mathcal{F}(E)\otimes_{\sigma}H$ as $\Psi(T):=\sum_{n\geq
0}^{\oplus}T^{(k)}$. Thus we may think of $\Psi(T)$ as $I\otimes T$, but it
must be remembered that $T$ maps $H$ to $E\otimes_{\sigma}H$. To check that
$(\Psi,\pi)$ is an isometric covariant representation of $(\sigma(M)^{\prime
},E^{\sigma})$ on $\mathcal{F}(E)\otimes_{\sigma}H$, observe that $\Psi$ is a
bimodule map, i.e., $\Psi$ satisfies the equation $\Psi(a\cdot T\cdot
b)=\pi(a)\Psi(T)\pi(b)$, $a,b\in\sigma(M)^{\prime}$, $T\in E^{\sigma}$, almost
by definition. As for the inner product, observe that if $T_{1},T_{2}\in
E^{\sigma}$, and if $\xi\otimes h,\;\eta\otimes k\in E^{\otimes n}%
\otimes_{\sigma}H$, then
\begin{multline*}
\langle\Psi(T_{1})(\xi\otimes h),\Psi(T_{2})(\eta\otimes k)\rangle=\langle
\xi\otimes T_{1}h,\eta\otimes T_{2}k\rangle=\langle T_{1}h,\sigma^{E}%
(\varphi(\langle\xi,\eta\rangle))T_{2}k\rangle\\
=\langle T_{1}h,T_{2}(\sigma(\langle\xi,\eta\rangle)k)\rangle=\langle
h,T_{1}^{\ast}T_{2}(\sigma(\langle\xi,\eta\rangle)k)\rangle\\
=\langle h,\sigma(\langle\xi,\eta\rangle)T_{1}^{\ast}T_{2}k\rangle=\langle
\xi\otimes h,\eta\otimes T_{1}^{\ast}T_{2}k\rangle=\langle\xi\otimes
h,\pi(T_{1}^{\ast}T_{2})(\eta\otimes k)\rangle.
\end{multline*}
Thus $(\Psi,\pi)$ is an isometric covariant representation of $E^{\sigma}$
viewed as a $C^{\ast}$-correspondence (See \cite[Definition 2.11]{MS98}.)
However, as we noted above in Remark \ref{usefulnote}, the normality of $\pi$
implies that $\Psi$ has the desired continuity properties. Thus $(\Psi,\pi)$
is an isometric covariant representation of $E^{\sigma}$ viewed as a $W^{\ast
}$-correspondence.

\begin{lemma}
\label{Lemma1.6}Suppose that for $i=1,2,$ $E_{i}$ is a $W^{\ast}%
$-correspondence over the von Neumann algebra $M$ and suppose that
$\sigma:M\rightarrow B(H)$ is a faithful normal representation of $M$ on $H$.
Then, under the map $\Lambda$ defined by the formula
\begin{equation}
\Lambda(\eta_{1}\otimes\eta_{2})=(I_{E_{2}}\otimes\eta_{1})\eta_{2},
\label{dualproduct}%
\end{equation}
$E_{1}^{\sigma}\otimes E_{2}^{\sigma}$ is isomorphic to $(E_{2}\otimes
E_{1})^{\sigma}$.
\end{lemma}

Observe that if $\eta_{i}\in E_{i}^{\sigma}$, then $\eta_{i}$ is a map from
$H$ to $E_{i}\otimes_{\sigma}H$. Therefore $(I_{E_{2}}\otimes\eta_{1})\eta
_{2}$ is a map from $H$ to $E_{2}\otimes E_{1}\otimes_{\sigma}H$, i.e., we see
that $(I_{E_{2}}\otimes\eta_{1})\eta_{2}$ can be represented schematically as%
\[
H\overset{\eta_{2}}{\rightarrow}E_{2}\otimes_{\sigma}H\overset{I_{E_{2}%
}\otimes\eta_{1}}{\rightarrow}E_{2}\otimes E_{1}\otimes_{\sigma}H\text{.}%
\]
Thus the proposed isomorphism $\Lambda$ makes good sense \emph{a priori}.

\begin{proof}
Note first that for $\eta_{i}\in E_{i}$, $i=1,2$, the map $(I_{E_{2}}%
\otimes\eta_{1})\eta_{2}$ satisfies the equation $(I_{E_{2}}\otimes\eta
_{1})\eta_{2}\sigma(a)=(I_{E_{2}}\otimes\eta_{1})(\varphi_{2}(a)\otimes
I_{H})\eta_{2}=(\varphi_{2}(a)\otimes\eta_{1})\eta_{2}=(\varphi_{2}(a)\otimes
I_{E_{1}}\otimes I_{H})(I_{E_{2}}\otimes\eta_{1})\eta_{2}$, where $a\in M$ and
$\varphi_{2}$ denotes the left action of $M$ on $E_{2}$.

Let $\eta_{i}$ and $\theta_{i}$ be elements of $E_{i}^{\sigma}$, $i=1,2$ and
compute: $\langle\eta_{1}\otimes\eta_{2},\theta_{1}\otimes\theta_{2}%
\rangle=\langle\eta_{2},\langle\eta_{1},\theta_{1}\rangle\cdot\theta
_{2}\rangle=\langle\eta_{2},(I\otimes\eta_{1}^{\ast}\theta_{1})\theta
_{2}\rangle=\eta_{2}^{\ast}(I_{E_{2}}\otimes\eta_{1})^{\ast}(I_{E_{2}}%
\otimes\theta_{1})\theta_{2}$. But this last expression is the inner product
$\langle(I_{E_{2}}\otimes\eta_{1})\eta_{2},(I_{E_{2}}\otimes\theta_{1}%
)\theta_{2}\rangle$ in $(E_{2}\otimes E_{1})^{\sigma}$. Thus, $\Lambda$
preserves inner products.

Next, we show that $\Lambda$ is a bimodule map. To this end, for $S\in
\sigma(M)^{\prime}$, calculate to find that on the one hand $\Lambda
(S\cdot(\eta_{1}\otimes\eta_{2}))=\Lambda((I\otimes S)\eta_{1}\otimes\eta
_{2})=(I_{E_{2}}\otimes(I_{E_{1}}\otimes S)\eta_{1})\eta_{2}=(I_{E_{2}\otimes
E_{1}}\otimes S)(I_{E_{2}}\otimes\eta_{1})\eta_{2}=(I_{E_{2}\otimes E_{1}%
}\otimes S)\Lambda(\eta_{1}\otimes\eta_{2})=S\cdot\Lambda(\eta_{1}\otimes
\eta_{2})$. While on the other, $\Lambda((\eta_{1}\otimes\eta_{2})\cdot
S)=\Lambda((\eta_{1}\otimes\eta_{2}S))=(I_{E_{2}}\otimes\eta_{1})\eta
_{2}S=\Lambda((\eta_{1}\otimes\eta_{2})\cdot S$. Hence, $\Lambda$ is a
bimodule isomorphism from $E_{1}^{\sigma}\otimes E_{2}^{\sigma}$ onto a closed
submodule of $(E_{2}\otimes E_{1})^{\sigma}$.

To show that $\Lambda$ is surjective, it will suffice to show that
$(\operatorname{Im}(\Lambda))^{\perp}=\{0\}$, since $(E_{2}\otimes
E_{1})^{\sigma}$ is selfdual. So fix $T\in(\operatorname{Im}(\Lambda))^{\perp
}$. Then for every $\eta_{i}\in E_{i}^{\sigma}$, $i=1,2$, we have $\eta
_{2}^{\ast}(I_{E_{2}}\otimes\eta_{1})^{\ast}T=0$. However, $(I_{E_{2}}%
\otimes\eta_{1})^{\ast}T$ maps $H$ into $E_{2}\otimes_{\sigma}H$ and satisfies
the equation $(I_{E_{2}}\otimes\eta_{1})^{\ast}T\sigma(a)=(I_{E_{2}}%
\otimes\eta_{1})^{\ast}(\varphi_{2}(a)\otimes I_{E_{1}}\otimes I_{H}%
)T=(\varphi_{2}(a)\otimes I_{H})(I_{E_{2}}\otimes\eta_{1})^{\ast}T$ for all
$a\in M$. Thus $(I_{E_{2}}\otimes\eta_{1})^{\ast}T$ lies in $E_{2}^{\sigma}$
and is orthogonal to every $\eta_{2}$ in $E_{2}^{\sigma}$. Thus $(I_{E_{2}%
}\otimes\eta_{1})^{\ast}T=0$. This, then, means that the range of $T$ is
orthogonal to $\bigvee\{(I_{E_{2}}\otimes\eta_{1})(E_{2}\otimes_{\sigma}%
H)\mid\eta_{1}\in E_{1}\}$. However, since $\bigvee\{\eta_{1}(H)\mid\eta
_{1}\in E_{1}\}=E_{1}\otimes_{\sigma}H$ by Lemma \ref{Lemma2.10MS02}, we find
that $T=0$.
\end{proof}

\begin{lemma}
\label{Lemma1.7}Let $E$ be a $W^{\ast}$-correspondence over the von Neumann
algebra $M$ and let $\sigma:M\rightarrow B(H)$ be a faithful normal
representation of $M$ on $H$. Then:

\begin{enumerate}
\item  For each $n\geq1$, the correspondences $(E^{\otimes n})^{\sigma}$ and
$(E^{\sigma})^{\otimes n}$ are isomorphic under the map $\Lambda
_{n}:(E^{\sigma})^{\otimes n}\rightarrow(E^{\otimes n})^{\sigma}$ defined by
the formula
\[
\Lambda_{n}(\eta_{1}\otimes\eta_{2}\otimes\cdots\otimes\eta_{n}%
)=(I_{E^{\otimes n-1}}\otimes\eta_{1})(I_{E^{\otimes n-2}}\otimes\eta
_{2})\cdots(I_{E}\otimes\eta_{n-1})\eta_{n}\text{,}%
\]
$\eta_{1}\otimes\eta_{2}\otimes\cdots\otimes\eta_{n}\in(E^{\sigma})^{\otimes
n}$.

\item  The map $U_{n}:(E^{\sigma})^{\otimes n}\otimes_{\iota}H\rightarrow
E^{\otimes n}\otimes_{\sigma}H$ defined by the formula%
\[
U_{n}(\eta\otimes h)=\Lambda_{n}(\eta)h\text{,}%
\]
$\eta\in(E^{\sigma})^{\otimes n}$, is a Hilbert space isomorphism mapping
$(E^{\sigma})^{\otimes n}\otimes_{\iota}H$ onto $E^{\otimes n}\otimes_{\sigma
}H$.

\item  The Hilbert space direct sum $U:=\sum^{\oplus}U_{n}$ maps
$\mathcal{F}(E^{\sigma})\otimes_{\iota}H$ isometrically onto $\mathcal{F}%
(E)\otimes_{\sigma}H$ and satisfies the equations%
\[
U\iota^{\mathcal{F}(E^{\sigma})}(T_{\eta})U^{\ast}=\Psi(\eta)
\]
and%
\[
U\iota^{\mathcal{F}(E^{\sigma})}(\varphi_{\infty}(a))U^{\ast}=\pi(a)\text{,}%
\]
for all $\eta\in E^{\sigma}$, $a\in\sigma(M)^{\prime}$, where $(\Psi,\pi)$ is
the covariant representation constructed in the paragraph before Lemma
\ref{Lemma1.6}.
\end{enumerate}
\end{lemma}

We note in passing that $U_{1}$ is just the map $\mu$ of Theorem \ref{duality}
and $U_{n}$ is just the natural ``inflation'' of $U_{1}$ to higher tensor
powers of $E^{\sigma}$.

\begin{proof}
The first assertion follows from Lemma \ref{Lemma1.6} and an easy induction.
Also, the second assertion is an immediate consequence of the first. Moreover,
of course, $U$ is a Hilbert space isomorphism from $\mathcal{F}(E^{\sigma
})\otimes_{\iota}H$ onto $\mathcal{F}(E)\otimes_{\sigma}H$. So, we need only
check that the indicated equations are satisfied. If $a\in\sigma(M)^{\prime}$,
then $\pi(a)U_{n}(\eta_{1}\otimes\eta_{2}\otimes\cdots\otimes\eta_{n}\otimes
h)=(I_{E^{\otimes n}}\otimes a)(I_{E^{\otimes n-1}}\otimes\eta_{1}%
)(I_{E^{\otimes n-2}}\otimes\eta_{2})\cdots(I_{E}\otimes\eta_{n-1})\eta
_{n}h=(I_{E^{\otimes n-1}}\otimes a\cdot\eta_{1})(I_{E^{\otimes n-2}}%
\otimes\eta_{2})\cdots(I_{E}\otimes\eta_{n-1})\eta_{n}h=U_{n}(a\cdot\eta
_{1}\otimes\eta_{2}\otimes\cdots\otimes\eta_{n}\otimes h)=U_{n}\sigma
^{\mathcal{F}(E_{\sigma})}(\varphi_{\infty}(a))(\eta_{1}\otimes\eta_{2}%
\otimes\cdots\otimes\eta_{n}\otimes h)$ and, if $\eta\in E^{\sigma}$, then
$\Psi(\eta)U_{n}(\eta_{1}\otimes\eta_{2}\otimes\cdots\otimes\eta_{n}\otimes
h)=(I_{E^{\otimes n}}\otimes\eta)(I_{E^{\otimes n-1}}\otimes\eta
_{1})(I_{E^{\otimes n-2}}\otimes\eta_{2})\cdots(I_{E}\otimes\eta_{n-1}%
)\eta_{n}h=U_{n+1}(\eta\otimes\eta_{1}\otimes\eta_{2}\otimes\cdots\otimes
\eta_{n}\otimes h)=U_{n+1}\iota^{\mathcal{F}(E^{\sigma})}(T_{\eta})(\eta
_{1}\otimes\eta_{2}\otimes\cdots\otimes\eta_{n}\otimes h)$.
\end{proof}

We want to emphasize for the sake of reference that the formula for
$U:\mathcal{F}(E^{\sigma})\otimes_{\iota}H\rightarrow\mathcal{F}%
(E)\otimes_{\sigma}H$, expressed on decomposable tensors, is%
\begin{equation}
U(\eta_{1}\otimes\eta_{2}\cdots\eta_{n}\otimes h)=(I_{E^{\otimes n-1}}%
\otimes\eta_{1})(I_{E^{\otimes n-2}}\otimes\eta_{2})\cdots(I_{E}\otimes
\eta_{n-1})\eta_{n}h\text{.} \label{Uformula}%
\end{equation}

\begin{theorem}
\label{commutant}Let $E$ be a $W^{\ast}$-correspondence over the von Neumann
algebra $M$ and let $\sigma:M\rightarrow B(H)$ be a faithful normal
representation of $M$ on $H$. With the notation of the preceding lemma, let
$\rho$ be the representation of $H^{\infty}(E^{\sigma})$ defined by the
formula%
\begin{equation}
\rho(X)=U\iota^{\mathcal{F}(E^{\sigma})}(X)U^{\ast}\text{,}
\label{dualityisomorphism}%
\end{equation}
$X\in H^{\infty}(E^{\sigma})$. Then $\rho$ is an ultraweakly continuous,
completely isometric representation of $H^{\infty}(E^{\sigma})$ on
$\mathcal{F}(E)\otimes_{\sigma}H$ that extends the representation $\pi
\times\Psi$ of $\mathcal{T}_{+}(E^{\sigma})$. Moreover, the range of $\rho$ is
the commutant of $\sigma^{\mathcal{F}(E)}(H^{\infty}(E))$.
\end{theorem}

\begin{proof}
The first assertion is immediate from Lemma \ref{Lemma1.7}. To see that
$\rho(H^{\infty}(E^{\sigma}))\subseteq\sigma^{\mathcal{F}(E)}(H^{\infty
}(E))^{\prime}$, it is enough to show that $\Psi(\eta)$, $\eta\in E^{\sigma},$
and $\pi\left(  a\right)  $, $a\in\sigma(M)^{\prime}$, commute with
$\sigma^{\mathcal{F}(E)}(T_{\xi})$ and $\sigma^{\mathcal{F}(E)}(\varphi
_{\infty}(b))$ for all $\xi\in E$ and $b\in M$. Evaluating on vectors of the
form $\xi_{1}\otimes\xi_{2}\otimes\cdots\xi_{n}\otimes h\in E^{\otimes
n}\otimes_{\sigma}H$, we find that $\Psi(\eta)\sigma^{\mathcal{F}(E)}(T_{\xi
})(\xi_{1}\otimes\xi_{2}\otimes\cdots\xi_{n}\otimes h)=\Psi(\eta)(\xi
\otimes\xi_{1}\otimes\xi_{2}\otimes\cdots\xi_{n}\otimes h)=\xi\otimes\xi
_{1}\otimes\xi_{2}\otimes\cdots\xi_{n}\otimes\eta h=\sigma^{\mathcal{F}%
(E)}(T_{\xi})\Psi(\eta)(\xi_{1}\otimes\xi_{2}\otimes\cdots\xi_{n}\otimes h)$.
The proof of the other relations are similar.

For the reverse inclusion, $\sigma^{\mathcal{F}(E)}(H^{\infty}(E))^{\prime
}\subseteq\rho(H^{\infty}(E^{\sigma}))$, recall the unitary group
$\{W_{t}\}_{t\in\mathbb{R}}$ in $\mathcal{L}(\mathcal{F}(E))$ defined by
formula (\ref{unitarygroup}). We write $W_{t}^{\sigma}$ for $\sigma
^{\mathcal{F}(E)}(W_{t})$ and we write $P_{n}^{\sigma}$ for $\sigma
^{\mathcal{F}(E)}(P_{n})$. Since $\gamma_{t}=AdW_{t}$ leaves $H^{\infty}(E)$
invariant, $AdW_{t}^{\sigma}$ leaves $\sigma^{\mathcal{F}(E)}(H^{\infty}(E))$
\emph{and }$\sigma^{\mathcal{F}(E)}(H^{\infty}(E))^{\prime}$ invariant. Write
$\Phi_{j}^{\sigma}$ and $\Sigma_{k}^{\sigma}$ for the operators on
$B(\mathcal{F}(E)\otimes_{\sigma}H)$ that are defined in a fashion analogous
to the definitions given at the end of the preceding section: $\Phi
_{j}^{\sigma}:=\frac{1}{2\pi}\int_{0}^{2\pi}e^{-itj}AdW_{t}^{\sigma}\,dt$ and
$\Sigma_{k}^{\sigma}:=\sum_{|j|<k}(1-\frac{|j|}{k})\Phi_{j}^{\sigma}$.
(Alternatively, $\Phi_{j}^{\sigma}=\sigma^{\mathcal{F}(E)}\circ\Phi_{j}$ and
$\Sigma_{k}^{\sigma}=\sigma^{\mathcal{F}(E)}\circ\Sigma_{k}$.) Since
$\Sigma_{k}^{\sigma}(a)\rightarrow a$ ultrastrongly for all $a\in
B(\mathcal{F}(E)\otimes_{\sigma}H))$, it clearly suffices to show that for all
$a\in\sigma^{\mathcal{F}(E)}(H^{\infty}(E))^{\prime}$, the operator $\Phi
_{j}^{\sigma}(a)$ belongs to $\rho(H^{\infty}(E^{\sigma}))$ for all $j$. To
this end, suppose first that $j<0$ and that $h$ lies in $H$, viewed as the
zero$^{th}$ summand of $\mathcal{F}(E)\otimes_{\sigma}H$. Then $\Phi
_{j}^{\sigma}(a)h=\frac{1}{2\pi}\int_{0}^{2\pi}e^{-itj}W_{t}^{\sigma}%
aW_{t}^{\sigma\ast}h\,dt=(\frac{1}{2\pi}\int_{0}^{2\pi}e^{-itj}W_{t}^{\sigma
})ah\,dt=0$. On the other hand, since $\Phi_{j}^{\sigma}(a)$ commutes with
$\sigma^{\mathcal{F}(E)}(H^{\infty}(E))$ and since $\sigma^{\mathcal{F}%
(E)}(H^{\infty}(E))H$ is dense in $\mathcal{F}(E)\otimes_{\sigma}H$, we see
that $\Phi_{j}^{\sigma}(a)$ must be the zero operator. So now fix $j\geq0$ and
consider the operator $b:=P_{j}^{\sigma}\Phi_{j}^{\sigma}(a)P_{0}^{\sigma}$,
viewed as an operator from $H$ to $E^{\otimes j}\otimes_{\sigma}H$. Then for
all $x\in M$, we have $b\sigma(x)=P_{j}^{\sigma}\Phi_{j}^{\sigma}%
(a)P_{0}^{\sigma}\sigma(x)=P_{j}^{\sigma}\Phi_{j}^{\sigma}(a)\sigma
^{\mathcal{F}(E)}(\varphi_{\infty}(x))P_{0}^{\sigma}=P_{j}^{\sigma}%
\sigma^{\mathcal{F}(E)}(\varphi_{\infty}(x))\Phi_{j}^{\sigma}(a)P_{0}^{\sigma
}=(\varphi^{(j)}(x)\otimes I_{H})P_{j}^{\sigma}\Phi_{j}^{\sigma}%
(a)P_{0}^{\sigma}$, which shows that $b\in(E^{\otimes j})^{\sigma}$. But for
$h\in H$ we have $\rho(b)h=bh=P_{j}^{\sigma}\Phi_{j}^{\sigma}(a)P_{0}^{\sigma
}h=P_{j}^{\sigma}\Phi_{j}^{\sigma}(a)h=\Phi_{j}^{\sigma}(a)h$ because
$\Phi_{j}^{\sigma}(a)$ maps the range of $P_{0}^{\sigma}$ to the range of
$P_{j}^{\sigma}$. So, given $\xi\otimes h\in E^{\otimes n}\otimes H$, we have
$\rho(b)(\xi\otimes h)=\xi\otimes bh=\sigma^{\mathcal{F}(E)}(T_{\xi})\Phi
_{j}^{\sigma}(a)h=\Phi_{j}^{\sigma}(a)\sigma^{\mathcal{F}(E)}(T_{\xi}%
)h=\Phi_{j}^{\sigma}(a)(\xi\otimes h)$. Thus, $\Phi_{j}^{\sigma}(a)=\rho(b)$,
an element in $\rho(H^{\infty}(E^{\sigma}))$.
\end{proof}

As we noted in the introduction, this result leads to the following corollary
which is a double commutant theorem for induced representations of Hardy algebras.

\begin{corollary}
\label{doublecommutant}If $\sigma$ is a faithful normal representation of $M$
on a Hilbert space $H$, then $\sigma^{\mathcal{F}(E)}(H^{\infty}%
(E))^{\prime\prime}=\sigma^{\mathcal{F}(E)}(H^{\infty}(E))$.
\end{corollary}

\begin{proof}
The proof rests on Theorem \ref{commutant} and on Theorem \ref{duality} and is
a matter of tracing through a string of definitions. First recall from Lemma
\ref{Lemma1.7} that we have a Hilbert space isomorphism $U:\mathcal{F}%
(E^{\sigma})\otimes_{\iota}H\rightarrow\mathcal{F}(E)\otimes_{\sigma}H$. Then
by Theorem \ref{commutant} $U\iota^{\mathcal{F}(E^{\sigma})}(H^{\infty
}(E^{\sigma}))U^{\ast}$ is the commutant of $\sigma^{\mathcal{F}(E)}%
(H^{\infty}(E))$. On the other hand, there is also a similarly defined Hilbert
space isomorphism $V:\mathcal{F}(E^{\sigma,\iota})\otimes_{\sigma}%
H\rightarrow\mathcal{F}(E^{\sigma})\otimes_{\iota}H$ so that $V\sigma
^{\mathcal{F}(E^{\sigma,\iota})}(H^{\infty}(E^{\sigma,\iota}))V^{\ast}$ is the
commutant of $\iota^{\mathcal{F}(E^{\sigma})}(H^{\infty}(E^{\sigma}))$.
(Recall that $E^{\sigma,\iota}$ is a correspondence over $M$, so this makes
sense.) Hence we see that the double commutant of $\sigma^{\mathcal{F}%
(E)}(H^{\infty}(E))$ is $UV\sigma^{\mathcal{F}(E^{\sigma,\iota})}(H^{\infty
}(E^{\sigma,\iota}))V^{\ast}U^{\ast}$. But also, the correspondence
isomorphism $W$ of Theorem \ref{duality} induces a Hilbert space isomorphism
$\tilde{W}$ from $\mathcal{F}(E)\otimes_{\sigma}H$ onto $\mathcal{F}%
(E^{\sigma,\iota})\otimes_{\sigma}H$ via the formula $\tilde{W}(\xi_{1}%
\otimes\xi_{2}\otimes\cdots\xi_{n}\otimes h)=W(\xi_{1})\otimes W(\xi
_{2})\otimes\cdots W(\xi_{n})\otimes h$, and it is a straightforward matter to
check that $\tilde{W}\sigma^{\mathcal{F}(E)}(H^{\infty}(E))\tilde{W}^{\ast
}=\sigma^{\mathcal{F}(E^{\sigma,\iota})}(H^{\infty}(E^{\sigma,\iota}))$. Thus
the double commutant of $\sigma^{\mathcal{F}(E)}(H^{\infty}(E))$ is
$(UV\tilde{W})\sigma^{\mathcal{F}(E)}(H^{\infty}(E))(UV\tilde{W})^{\ast}$, and
so it suffices to prove that $UV\tilde{W}$ is the identity.

Since each of $U$, $V$ and $\tilde{W}$ is determined by how it acts on
homogeneous spaces of tensors, we may focus our attention on the restriction
of $UV\tilde{W}$ to $E^{\otimes n}\otimes_{\sigma}H$. We shall subscript these
operators with $n$ to indicate their restrictions to the appropriate spaces of
$n$-homogeneous tensors. First note that the definition of $V$ is given by
equation (\ref{Uformula}), but with the appropriate modifications:
\[
V_{n}(\eta_{1}\otimes\eta_{2}\cdots\eta_{n}\otimes h)=(I_{(E^{\sigma
})^{\otimes n-1}}\otimes\eta_{1})(I_{(E^{\sigma})^{\otimes n-2}}\otimes
\eta_{2})\cdots(I_{(E^{\sigma})}\otimes\eta_{n-1})\eta_{n}h\text{,}%
\]
where $\eta_{1}\otimes\eta_{2}\cdots\eta_{n}\otimes h\in(E^{\sigma,\iota
})^{\otimes n}\otimes_{\sigma}H$. Using this, and the definition of $\tilde
{W}$, we see that for $\xi_{1}\otimes\xi_{2}\cdots\xi_{n}\otimes h\in
E^{\otimes n}\otimes_{\sigma}H$
\begin{multline*}
V_{n}\tilde{W}_{n}(\xi_{1}\otimes\xi_{2}\cdots\xi_{n}\otimes h)\\
=(I_{(E^{\sigma})^{\otimes n-1}}\otimes W(\xi_{1}))(I_{(E^{\sigma})^{\otimes
n-2}}\otimes W(\xi_{2}))\cdots(I_{(E^{\sigma})}\otimes W(\xi_{n-1}))W(\xi
_{n})h\text{.}%
\end{multline*}
Thus, we must prove that
\begin{align}
U_{n}((I_{(E^{\sigma})^{\otimes n-1}}\otimes W(\xi_{1}))(I_{(E^{\sigma
})^{\otimes n-2}}\otimes W(\xi_{2}))\cdots(I_{(E^{\sigma})}\otimes W(\xi
_{n-1}))W(\xi_{n})h)\nonumber\\
=\xi_{1}\otimes\xi_{2}\cdots\xi_{n}\otimes h\text{.} \label{BasicEquation}%
\end{align}
We will show in detail that this equation is valid when $n=1,2$, but we shall
set things up so that checking the equation for larger $n$ will be a
straightforward, albeit tedious, matter.

The problem with managing the formulae is that elements in $E^{\sigma}$ and in
$E^{\sigma,\iota}$ play two roles: sometimes they are viewed as operators and
sometimes they are viewed as elements of a space to which other operators are
applied. So the key is to express various equations as equations among
operators. For this purpose observe that we may express vectors of the form
$\xi_{1}\otimes\xi_{2}\cdots\xi_{n}\otimes h\in E^{\otimes n}\otimes_{\sigma
}H$ as $L_{\xi_{1}}^{(n)}L_{\xi_{2}}^{(n-1)}\cdots L_{\xi_{n}}^{(1)}h$, where
for $\zeta\in E$, $L_{\zeta}^{(k)}:E^{\otimes k-1}\otimes_{\sigma}H\rightarrow
E^{\otimes k}\otimes_{\sigma}H$ is defined by the formula $L_{\zeta}^{(k)}%
\eta=\zeta\otimes\eta$. We have met operators of this form before. However,
what is important here is to keep track of the indices and to know how to
calculate their adjoints. But the adjoints are easy: For $\xi_{1}\otimes
\xi_{2}\cdots\xi_{k}\otimes h\in E^{\otimes k}\otimes_{\sigma}H$ and $\zeta\in
E$,
\begin{multline*}
(L_{\zeta}^{(k)})^{\ast}\xi_{1}\otimes\xi_{2}\cdots\xi_{k}\otimes
h=(\sigma^{E^{\otimes k-1}}\circ\varphi(\langle\zeta,\xi_{1}\rangle))(\xi
_{2}\otimes\xi_{3}\cdots\xi_{k}\otimes h)\\
=((\varphi_{k-1}(\langle\zeta,\xi_{1}\rangle)(\xi_{2}\otimes\xi_{3}\cdots
\xi_{k}))\otimes h\\
=(\varphi(\langle\zeta,\xi_{1}\rangle)\xi_{2})\otimes\xi_{3}\cdots\xi
_{k}\otimes h\text{.}%
\end{multline*}
Also, we need the following computation. Let $\eta\in E^{\sigma}$ and remember
that $\eta$ is a map from $H$ to $E\otimes_{\sigma}H$. Then for $\zeta\in E$,
the following equation holds between maps on $E^{\otimes k}\otimes_{\sigma}H$%
\begin{equation}
(I_{E^{\otimes k-1}}\otimes\eta)(L_{\zeta}^{(k)})^{\ast}=(L_{\zeta}%
^{(k+1)})^{\ast}(I_{E^{\otimes k}}\otimes\eta)\text{.} \label{adjoint1}%
\end{equation}
Turning to the verification of equation (\ref{BasicEquation}), we see that
when expressed entirely in terms of operators, what we need to prove is
\begin{align}
U_{n}((I_{(E^{\sigma})^{\otimes n-1}}\otimes W(\xi_{1}))(I_{(E^{\sigma
})^{\otimes n-2}}\otimes W(\xi_{2}))\cdots(I_{(E^{\sigma})}\otimes W(\xi
_{n-1}))W(\xi_{n}))\nonumber\\
=L_{\xi_{1}}^{(n)}L_{\xi_{2}}^{(n-1)}\cdots L_{\xi_{n}}^{(1)}\text{,}
\label{OperatorEquat1}%
\end{align}
as an operator from $H$ to $E^{\otimes n}\otimes_{\sigma}H$. We shall do this
for $n=1,2$.

By definition, $W(\xi)=U_{1}^{\ast}L_{\xi}^{(1)}$ (See equation (\ref{wxs})
and the note just before the proof of Lemma \ref{Lemma1.7}.) So equation
(\ref{OperatorEquat1}) is satisfied in this case. To verify
(\ref{OperatorEquat1}) when $n=2$, it seems preferable to compute adjoints.
That is, we want to show that%
\begin{equation}
W(\xi_{2})^{\ast}(I_{E^{\sigma}}\otimes W(\xi_{1}))^{\ast}U_{2}^{\ast}%
=(L_{\xi_{2}}^{(1)})^{\ast}(L_{\xi_{1}}^{(2)})^{\ast}\text{.} \label{adjoint2}%
\end{equation}
We begin by applying the left hand side to a vector of the form $(I_{E}%
\otimes\eta_{1})\eta_{2}h\in E^{\otimes2}\otimes_{\sigma}H$. These span
$E^{\otimes2}\otimes_{\sigma}H$. By (\ref{Uformula}), $U_{2}^{\ast}%
(I_{E}\otimes\eta_{1})\eta_{2}h=\eta_{1}\otimes\eta_{2}\otimes h$ in
$(E^{\sigma})^{\otimes2}\otimes_{\iota}H$. But then
\begin{align}
W(\xi_{2})^{\ast}(I_{E^{\sigma}}\otimes W(\xi_{1}))^{\ast}U_{2}^{\ast}%
(I_{E}\otimes\eta_{1})\eta_{2}h  &  =W(\xi_{2})^{\ast}(I_{E^{\sigma}}\otimes
W(\xi_{1}))^{\ast}(\eta_{1}\otimes\eta_{2}\otimes h)\nonumber\\
&  =W(\xi_{2})^{\ast}(\eta_{1}\otimes(W(\xi_{1})^{\ast}(\eta_{2}\otimes
h)))\text{.} \label{adjoint3}%
\end{align}
By definition, $W(\xi_{1})=U_{1}^{\ast}L_{\xi_{1}}^{(1)}$, so
\[
W(\xi_{1})^{\ast}(\eta_{2}\otimes h)=(L_{\xi_{1}}^{(1)})^{\ast}U_{1}(\eta
_{2}\otimes h)=(L_{\xi_{1}}^{(1)})^{\ast}\eta_{2}h\text{.}%
\]
Hence continuing with (\ref{adjoint3}), we have%
\begin{align}
W(\xi_{2})^{\ast}(\eta_{1}\otimes(W(\xi_{1})^{\ast}(\eta_{2}\otimes h)))  &
=W(\xi_{2})^{\ast}(\eta_{1}\otimes((L_{\xi_{1}}^{(1)})^{\ast}\eta
_{2}h))\nonumber\\
&  =(L_{\xi_{2}}^{(1)})^{\ast}\eta_{1}(L_{\xi_{1}}^{(1)})^{\ast}\eta
_{2}h\text{.} \label{adjoint4}%
\end{align}
However, by equation (\ref{adjoint1}), the last term in (\ref{adjoint4})
equals $(L_{\xi_{2}}^{(1)})^{\ast}(L_{\xi_{1}}^{(2)})^{\ast}(I_{E}\otimes
\eta_{1})\eta_{2}h$. Inserting this expression into equation (\ref{adjoint4})
and working backwards through equation (\ref{adjoint3}), we find that, indeed,
equation (\ref{adjoint2}) is satisfied, which is what we wanted to show.
\end{proof}

\section{Examples of Correspondence Duals}

In this section we collect a number of examples of duals of various correspondences.

\begin{example}
\label{Example2.1}Let $M=\mathbb{C}$ and let $E=\mathbb{C}$ to get what is the
simplest correspondence over the simplest von Neumann algebra. Let $H$ be an
arbitrary Hilbert space and let $\sigma$ be the representation of $M$ on $H$
(the only representation of $M$ on $H$) given by multiplication: $\sigma
(a)\xi=a\xi$. Then $\sigma(M)^{\prime}=B(H)$, of course, and $E^{\sigma
}=\{\eta:H\rightarrow\mathbb{C}\otimes_{\sigma}H\mid\eta a=(a\otimes
I_{H})\eta$, $a\in\mathbb{C\}}$. Evidently, $E^{\sigma}$ may be identified
with $B(H)$ as a $B(H)$-bimodule, where the left and right actions are given
by operator multiplication. The $B(H)$-valued inner product on $E^{\sigma}$ is
given by the formula $\langle\eta_{1},\eta_{2}\rangle=\eta_{1}^{\ast}\eta_{2}$.
\end{example}

The Fock space in this example is $\mathcal{F}(E)$ is $\ell^{2}(\mathbb{Z}%
_{+})$, which, via the Fourier transform, may be identified with the Hardy
space $H^{2}(\mathbb{T)}$. The Hardy algebra $H^{\infty}(E)$, as defined, is
the collection of all lower triangular Toeplitz matrices, which, of course, is
(completely) isometrically isomorphic to $H^{\infty}(\mathbb{T})$. The induced
space, $\mathcal{F}(E)\otimes_{\sigma}H$, is just the Hilbert space direct sum
of copies of $H$ indexed by $\mathbb{Z}_{+}$ and $\sigma^{\mathcal{F}(E)}$ is
simply the representation of $H^{\infty}(\mathbb{T})$ given matricially by the
formula%
\[
\sigma^{\mathcal{F}(E)}(f)=\left[
\begin{array}
[c]{cccc}%
{\hat{f}}_{0} & 0 & 0 & \\
{\hat{f}}_{1} & {\hat{f}}_{0} & 0 & \ddots\\
{\hat{f}}_{2} & {\hat{f}}_{1} & {\hat{f}}_{0} & \ddots\\
& \ddots & \ddots & \ddots
\end{array}
\right]  \text{,}%
\]
where ${\hat{f}}_{k}$ denotes the $k^{th}$ Taylor coefficient of $f$ times the
identity operator on $H$, for $f\in H^{\infty}(\mathbb{T})$. On the other
hand, the $n$-fold tensor product of $E^{\sigma}$ with itself is naturally
isomorphic to $B(H)$, for each $n$, so the Fock space $\mathcal{F}(E^{\sigma
})$ is just the direct sum of infinitely many copies of $B(H)$, with the
column space inner product. The covariant representation $(\Psi,\pi)$ is
realized as follows. First note that for $S\in\sigma(M)^{\prime}=B(H)$ and
$\eta\in E^{\sigma}=B(H)$,
\[
\sigma^{\mathcal{F}(E)}(\varphi_{\infty}(S))(T_{0}\otimes h_{0},T_{1}\otimes
h_{1},\cdots)=((ST_{0})\otimes h_{0},(ST_{1})\otimes h_{1},\cdots)
\]
and%
\[
\sigma^{\mathcal{F}(E)}(T_{\eta})(T_{0}\otimes h_{0},T_{1}\otimes h_{1}%
,\cdots)=(0,(\eta T_{0})\otimes h_{0},(\eta T_{1})\otimes h_{1},\cdots
)\text{.}%
\]
The operator $U:\mathcal{F}(E^{\sigma})\otimes H\rightarrow\mathcal{F}%
(E)\otimes H$ of Lemma \ref{Lemma1.7} maps $T_{i}\otimes h_{i}$ to $T_{i}%
h_{i}$. Hence $\pi(S)=U\sigma^{\mathcal{F}(E^{\sigma})}(S)U^{\ast}$ and
$\Psi(\eta)=U\sigma^{\mathcal{F}(E^{\sigma})}(T_{\eta})U^{\ast}$ are given by
the formulae%
\[
\pi(S)(h_{0},h_{1},\cdots)=(Sh_{0},Sh_{1},\cdots)
\]
and
\[
\Psi(\eta)(h_{0},h_{1},\cdots)=(0,\eta h_{0},\eta h_{1},\cdots)\text{.}%
\]
The ultraweakly closed algebra generated by these operators is manifestly the
commutant of $\sigma^{\mathcal{F}(E)}(H^{\infty}(E))$.

\begin{example}
\label{Example2.2}Let $M=\mathbb{C}$ again, but let $E=\mathbb{C}^{n}$. Also,
as in Example \ref{Example2.1}, let $H$ be an arbitrary Hilbert space, with
the same representation $\sigma$ of $M$. Then $M^{\prime}=B(H)$ and
$E^{\sigma}=\{\eta:H\rightarrow\mathbb{C}^{n}\otimes H\mid\eta a=(a\otimes
I_{H})\eta$, $a\in\mathbb{C\}=}B(H,\mathbb{C}^{n}\otimes H)$. The right action
of $M^{\prime}=B(H)$ on $E^{\sigma}$ is via multiplication, while the left
action of $S\in B(H)$ on $E^{\sigma}$ is given by the formula $S\cdot
\eta=(I\otimes S)\eta$. The inner product is $\langle\eta_{1},\eta_{2}%
\rangle=\eta_{1}^{\ast}\eta_{2}$. Hence, $E^{\sigma}$ is isomorphic to the
$n$-fold column space over $B(H)$, $C_{n}(B(H))$, with the left action
$S\cdot(T_{i})=(ST_{i})$.
\end{example}

Of course $\mathcal{F}(E)=\mathcal{F}(\mathbb{C}^{n})$ and, thus, $H^{\infty
}(E)$ is the algebra that Davidson and Pitts denoted by $\mathfrak{L}_{n}$ in
\cite{DP98a}. It is the ultraweak closure of Popescu's noncommutative disc
algebra (in the Fock representation) \cite{gP96a}. The Fock space
$\mathcal{F}(E^{\sigma})$ is the direct sum $\sum_{k\geq0}C_{n}(B(H))^{\otimes
k}$ and the operator $U:\mathcal{F}(E^{\sigma})\otimes H\rightarrow
\mathcal{F}(E)\otimes H$ of Lemma \ref{Lemma1.7} maps a typical decomposable
element $(T_{i}^{(1)})\otimes(T_{i}^{(2)})\otimes\cdots\otimes(T_{i}%
^{(k)})\otimes h$ in $(E^{\sigma})^{\otimes k}\otimes_{\sigma}H$ to $\sum
e_{i_{k}}\otimes e_{i_{k-1}}\otimes\cdots\otimes e_{i_{1}}\otimes(T_{i_{1}%
}^{(1)}T_{i_{2}}^{(2)}\cdots T_{i_{k}}^{(k)})h$ in $E^{\otimes k}%
\otimes_{\sigma}H$, where $\{e_{i}\}_{i=1}^{n}$ is the standard basis for
$\mathbb{C}^{n}$.

To see how $\pi$ and $\Psi$ work for this example, fix $\eta=(\eta_{i})\in
E^{\sigma}=\mathbb{C}^{n}(B(H))$ and $S\in M^{\prime}=B(H)$. Then if
$\alpha_{1}\otimes\alpha_{2}\otimes\cdots\otimes\alpha_{k}\otimes
h\in(\mathbb{C}^{n})^{\otimes k}\otimes H\subseteq\mathcal{F}(E)\otimes
_{\sigma}H$ is a typical decomposable tensor, we have
\[
\pi(S)(\alpha_{1}\otimes\alpha_{2}\otimes\cdots\otimes\alpha_{k}\otimes
h)=\alpha_{1}\otimes\alpha_{2}\otimes\cdots\otimes\alpha_{k}\otimes Sh
\]
and%
\[
\Psi(\eta)(\alpha_{1}\otimes\alpha_{2}\otimes\cdots\otimes\alpha_{k}\otimes
h)=\sum_{j=1}^{n}\alpha_{1}\otimes\alpha_{2}\otimes\cdots\otimes\alpha
_{k}\otimes e_{j}\otimes\eta h\text{.}%
\]
The ultraweakly closed algebra generated by these operators is $\rho
(H^{\infty}(E^{\sigma}))$, which is completely isometrically isomorphic and
ultraweakly homeomorphic to $H^{\infty}(E^{\sigma})$, and coincides with the
commutant of $\mathfrak{L}_{n}\otimes I_{H}=\sigma^{\mathcal{F}(E)}%
(\mathfrak{L}_{n})$. When $\dim H=1$, this result was obtained by Popescu in
\cite[Corollary 1.3]{gP95} and by Davidson and Pitts in \cite[Theorem
1.2]{DP99}.

\begin{example}
\label{Example2.3} (Quiver algebras or ``non-self-adjoint Cuntz-Krieger
algebras'') We follow the notation and terminology from \cite[Section 5]%
{MS99}. A (finite) quiver is simply a directed graph with $n$ vertices, say,
$\{v_{1},v_{2},\cdots,v_{n}\}$ and $c_{ij}$ edges or arrows from $v_{j}$ to
$v_{i}$ - $c_{ij}$ being a nonnegative integer. We write $C$ for the $n\times
n$ matrix $\left[  c_{ij}\right]  $. Our von Neumann algebra $M$ in this
setting is simply $D_{n}$, the algebra of $n\times n$ diagonal matrices. The
correspondence that we construct from $C$, $E(C)$, is the direct sum of
Hilbert spaces, $\sum_{i,j=1}^{n}\oplus H_{ij}$, where $\dim H_{ij}=c_{ij}$.
For each $(i,j)$ such that $c_{ij}>0$, we fix an orthonormal basis
$\{e_{ij}^{(k)}\mid1\leq k\leq c_{ij}\}$. This choice determines a $D_{n}%
$-$D_{n}$ bimodule structure on $E(C)$ via the formulae:%
\[
e_{ij}^{(k)}e_{ll}=\delta_{jl}e_{ij}^{(k)}%
\]
and%
\[
\varphi(e_{ll})e_{ij}^{(k)}=\delta_{il}e_{ij}^{(k)}\text{,}%
\]
where $\delta_{kl}$ denotes the Kronecker delta and where $\varphi$, as usual,
denotes the left action. Then with respect to the $D_{n}$-valued inner product
defined by the formula%
\[
\langle e_{ij}^{(k)},e_{lm}^{(p)}\rangle=\delta_{kp}\delta_{il}\delta
_{jm}e_{jj}\text{,}%
\]
$E(C)$ becomes a $W^{\ast}$-correspondence over $D_{n}$. As shown in
\cite[Corollary 5.2]{MS99}, $E(C)^{\otimes n}$ is naturally isomorphic to
$E(C^{n})$ as $W^{\ast}$-correspondences over $D_{n}$, and so these spaces
will be identified. Note, in particular, that since $C^{0}$ is the identity
matrix, this identification may be made even when $n=0$. We have then, of
course, the equation $\mathcal{F}(E(C))=\sum_{n\geq0}E(C^{n})$.

If $\sigma:D_{n}\rightarrow B(H)$ is a representation, then $\sigma$ is
determined by the multiplicities, $m_{k}$, of the one-dimensional
representations $\delta_{k}$ given by the formulae $\delta_{k}(a)=a_{kk}$,
$a\in D_{n}$. The representation $\sigma$ is completely determined up to
unitary equivalence by the vector $\mathbf{m}=(m_{1},m_{2},\cdots,m_{n})$ of
multiplicities. Note that $0\leq m_{i}\leq\infty$, and that $\sigma$ is
faithful (which is the only situation we consider here) precisely when each
$m_{i}$ is positive. The Hilbert space $H$, then is naturally written as
$H=\sum_{k=1}^{n}\mathbb{C}^{m_{k}}$, where $\mathbb{C}^{\infty}$ is
interpreted as $\ell^{2}(\mathbb{Z}_{+})$. Evidently, $\sigma(M)^{\prime}$ is
$\sum_{k=1}^{n}B(\mathbb{C}^{m_{k}})$. By definition the correspondence
$E^{\sigma}(C)$ is $\{\eta:H\rightarrow E(C)\otimes_{\sigma}H\mid\eta
a=(a\otimes I)\eta$, $a\in D_{n}\}$. This space may be written as a family of
operator matrices because for each $e_{kk}\in D_{n}$ and each $\eta\in
E^{\sigma}(C)$, $\eta e_{kk}=(e_{kk}\otimes I)\eta$. Since $\sigma
(e_{kk})H=\mathbb{C}^{m_{k}}$, while $\sigma^{E(C)}\circ\varphi(e_{kk}\otimes
I)E(C)\otimes H=\sum_{i=1}^{n}H_{ki}\otimes\mathbb{C}^{m_{i}}$. Thus
$E^{\sigma}(C)=\{\eta=\left(  \eta_{ij}\right)  \mid\eta_{ij}\in
B(\mathbb{C}^{m_{j}},H_{ji}\otimes\mathbb{C}^{m_{i}})\}$.
\end{example}

In the special case when the graph has one vertex and $n$ edges, Theorem
\ref{commutant} is Theorem 1.2 of \cite{DP99}. For a general (finite) graph,
Theorem \ref{commutant} generalizes Proposition 5.4 of \cite{MS99}. A bit more
specifically, in \cite{MS99}, we considered representations $\sigma$ of
uniform multiplicity one, i.e., we assumed that $m_{i}=1$ for all $i$. In this
case, $B(\mathbb{C}^{m_{j}},H_{ji}\otimes\mathbb{C}^{m_{i}})\}$ is naturally
isomorphic to $H_{ji}$ and we see that $E^{\sigma}=E(C^{t})$. From this,
Proposition 5.4 of \cite{MS99} is immediate. We note, too, that Theorem 5.8 of
\cite{KPp02}, which extends Proposition 5.4 of \cite{MS99} beyond the
multiplicity one case, also follows from Theorem \ref{commutant}.

\begin{example}
\label{Example2.4}In this example, we fix a von Neumann algebra $M$ acting on
a Hilbert space $H$ and we let $P:M\rightarrow M$ be a unital, normal
completely positive map. Following \cite[Section 5]{wP73} (see \cite[p.
48]{cL94} also) we use $P$ to build a $W^{\ast}$-correspondence $E$ over $M$
as follows. Form the algebraic tensor product $M\otimes M$ and endow it with
the inner product $\langle x_{1}\otimes y_{1},x_{2}\otimes y_{2}\rangle
=y_{1}^{\ast}P(x_{1}^{\ast}x_{2})y_{2}$. The complete positivity of $P$
guarantees that this inner product is positive semidefinite. We define
$E=M\otimes_{P}M$ to be the Hilbert space completion of $M\otimes M$ in this
inner product (with the vectors of length zero modded out). The left and right
actions of $M$ on $E$ are the obvious ones and with respect to them, $E$ is a
$W^{\ast}$-correspondence over $M$. The dual correspondence is $E^{\sigma
}=\{\eta:H\rightarrow M\otimes_{P}M\otimes_{\iota}H\mid\eta a=(a\otimes
I_{H})\eta$, $a\in M\}$, where $\iota$ denotes the identity representation of
$M$ on $H$. However, a moment's reflection reveals that $M\otimes_{P}%
M\otimes_{\iota}H=M\otimes_{P}H$, where $M\otimes_{P}H$ denotes the
Stinespring dilation space for the map $P$. Further, $\iota^{M\otimes_{P}%
M}\circ\varphi=\pi_{P}$, where $\varphi$ denotes the left action of $M$ on $E$
and $\pi_{P}$ denotes the Stinespring representation that dilates $P$. Thus,
in this case, we see that $E^{\sigma}$ is the correspondence that we denoted
$\mathcal{L}_{\mathcal{M}}(H,\mathcal{M}\otimes_{P}H)$ in \cite{MSp02}.
\end{example}

It is easy to see that $E^{\otimes k}=(M\otimes_{P}M)^{\otimes k}$ is
naturally isomorphic to the $k+1$-fold tensor product $\overset{k+1}%
{\overbrace{M\otimes_{P}M\otimes_{P}\cdots\otimes_{P}M}}$. We find, then, that
$\mathcal{F}(E)\otimes H=\sum_{k\geq0}\overset{k+1}{\overbrace{M\otimes
_{P}M\otimes_{P}\cdots\otimes_{P}M}}\otimes_{\iota}H$ and that for $S\in
M^{\prime}$, and $\eta\in E^{\sigma}$, $\pi(S)(a_{1}\otimes a_{2}\otimes\cdots
a_{k+1}\otimes h)=a_{1}\otimes a_{2}\otimes\cdots a_{k+1}\otimes Sh$, while
$\Psi(\eta)(a_{1}\otimes a_{2}\otimes\cdots a_{k+1}\otimes h)=a_{1}\otimes
a_{2}\otimes\cdots a_{k+1}\otimes\eta(h)\in E^{\otimes(k+1)}\otimes_{\iota}H$.
The ultraweakly closed algebra generated by $\{\pi(S),\Psi(\eta)\mid S\in
M^{\prime}$, $\eta\in E^{\sigma}\}$ is the commutant of $\sigma^{\mathcal{F}%
(E)}(H^{\infty}(E))$ by Theorem \ref{commutant}.

In the special case when $P$ is a $\ast$-endomorphism of $M$ (in which event
we shall write $P=\alpha$), we see that $E=M\otimes_{\alpha}M$ may be
identified with $M$ via the map $a\otimes b\mapsto\alpha(a)b$. (In this
situation, $E$ is usually written as $_{\alpha}M$.) When this identification
is made, the left action of $M$ on $E=\,_{\alpha}M$ is given by $\alpha$,
i.e., $\varphi(a)b=\alpha(a)b$. Also, $M\otimes_{\alpha}H$ may be identified
with $H$ via the map $a\otimes h\mapsto\alpha(a)h$ and when this
identification is made, then $E^{\sigma}$ becomes the intertwining space of
$\alpha$, $\{\eta\in B(H)\mid\eta a=\alpha(a)\eta$, $a\in M\}$. More
generally, we may identify $E^{\otimes k}\otimes_{\iota}H=\overset
{k+1}{\overbrace{M\otimes_{\alpha}M\otimes_{\alpha}\cdots\otimes_{\alpha}M}%
}\otimes_{\iota}H$ with $H$ via the map $t_{1}\otimes t_{2}\otimes
\cdots\otimes t_{k+1}\otimes h\mapsto\alpha^{k+1}(t_{1})\alpha^{k-1}%
(t_{2})\cdots\alpha(t_{1})h$. When this identification is made, $\mathcal{F}%
(E)\otimes_{\iota}H$ is just the direct sum of infinitely many copies of $H$,
$\sum_{k\geq0}H$, and we see that
\[
\iota^{\mathcal{F}(E)}\circ\varphi_{\infty}(a)(h_{0},h_{1},h_{2}%
,\cdots)=(ah_{0},\alpha(a)h_{1},\alpha^{2}(a)h_{2},\cdots)
\]
while
\[
\iota^{\mathcal{F}(E)}(T_{\xi})(h_{0},h_{1},h_{2},\cdots)=(0,\xi h_{0}%
,\alpha(\xi)h_{1},\alpha^{2}(\xi)h_{2},\cdots)\text{.}%
\]
Since $\iota$ is the identity representation of $M$, it is natural to identify
the ultraweakly closed algebra generated by $\{\iota^{\mathcal{F}(E)}%
\circ\varphi_{\infty}(a)$, $\iota^{\mathcal{F}(E)}(T_{\xi})\mid a\in M$,
$\xi\in E=\,_{\alpha}M\}$ with $H^{\infty}(E)$. However, we shall continue to
write this algebra as $\iota^{\mathcal{F}(E)}(H^{\infty}(E))$. By Theorem
\ref{commutant}, its commutant is generated by $\{\pi(s),\Psi(\eta)\mid s\in
M^{\prime}$, $\eta\in E^{\sigma}\}$, where $\pi$ and $\Psi$ are given by the
formulae%
\[
\pi(s)(h_{0},h_{1},h_{2},\cdots)=(sh_{0},sh_{1},sh_{2},\cdots),\;\;\;\;s\in
M^{\prime}\text{,}%
\]
and
\[
\Psi(\eta)(h_{0},h_{1},h_{2},\cdots)=(0,\eta h_{0},\eta h_{1},\eta
h_{2},\cdots)\text{,\ \ \ \ }\eta\in E^{\sigma}\text{.}%
\]
Of course, these formulae are familiar from the theory of self-adjoint crossed products.

\begin{example}
\label{Example2.5}Suppose that $M$ is a finite type $II$ factor, acting via
the left regular representation $\lambda$ on $H=L^{2}(M,\tau)$, where $\tau$
is a faithful normal tracial state. Suppose also that $N$ is a subfactor of
$M$ and that $\Phi:M\rightarrow N$ is the unique trace-preserving conditional
expectation of $M$ onto $N$. Then, as is well known, $\Phi$ is a completely
positive map on $M$ and so we may form $E:=M\otimes_{\Phi}M$. According to
Example \ref{Example2.4}, the dual of $E$ with respect to $\lambda$ is
$E^{\lambda}=\{\eta:H\rightarrow M\otimes_{\Phi}H\mid\eta\lambda(a)=(a\otimes
I_{H})\eta$, $a\in M\}$. We are going to show that $E^{\lambda}$ is naturally
isomorphic to $M^{op}\otimes_{\tilde{\Phi}}M^{op}$, where $M^{op}$ denotes the
opposite von Neumann algebra of $M$ and where $\tilde{\Phi}$ is the
conditional expectation on $M^{op}$ defined by the formula $\tilde{\Phi
}(x^{op})=\Phi(x)^{op}$. Here and throughout, we write $x\rightarrow x^{op}$
for the identity map regarded as an anti-isomorphism from $M$ onto $M^{op}$.
\end{example}

As is customary (see \cite{JS97}), let $J$ be the canonical antiunitary map on
$H=L^{2}(M,\tau)$ given by the formula $J(\lambda(a)\Omega)=\lambda(a)^{\ast
}\Omega$, where $\Omega$ denotes the identity operator in $M$ viewed as the
canonical cyclic and separating vector in $L^{2}(M,\tau)$. Then $J\lambda
(M)J=\rho(M)=\lambda(M)^{\prime}$, where $\rho$ denotes the right regular
representation of $M$. Usually, $\rho$ is treated as an antirepresentation of
$M$. However, we want to think of $\rho$ as a (faithful) representation of
$M^{op}$. Let $e$ denote the operator on $L^{2}(M,\tau)$ determined by the
conditional expectation $\Phi$ via the formula $e\lambda(a)\Omega=\lambda
(\Phi(a))\Omega$, $a\in M$. Then $e$ is a projection onto the closure of
$N\Omega$ in $L^{2}(M,\tau)$, and $e\lambda(x)e=\lambda(\Phi(x))e$ for all
$x\in M$. Note that $eJ=Je$ and, for $x\in M$, $\rho(x)=J\lambda(x^{*})J$. It
follows that $e\rho(x)e=eJ\lambda(x^{*})Je=Je\lambda(x^{*})eJ=J\lambda
(\Phi(x^{*}))eJ=\rho(\tilde{\Phi}(x^{*}))e$. We write $M_{1}$ for the von
Neumann algebra in $B(H)$ generated by $\lambda(M)$ and $e$. Then, as is well
known, $M_{1}=J\lambda(N)^{\prime}J$, so that $M_{1}^{\prime}=J\lambda(N)J$ is
contained in $\lambda(M)^{\prime}=\rho(M^{op})$.

For $b$ and $d$ in $M^{op}$, define $S(b,d)$ on $M\otimes_{\Phi}H$ $\ $mapping
to $H$ via the formula%
\[
S(b,d)(a\otimes h):=\rho(d)^{\ast}\lambda(a)e\rho(b)^{\ast}h\text{.}%
\]
We claim $S(b,d)$ is a well defined, bounded operator from $M\otimes_{\Phi}H$
to $H$ whose adjoint lies in $E^{\lambda}$. Once this is verified, we will
show that the map $\Psi:M^{op}\otimes_{\tilde{\Phi}}M^{op}\rightarrow
E^{\lambda}$ defined by $\Psi(b\otimes d)=S(b,d)^{\ast}$ is an isomorphism of
correspondences, where in $E^{\lambda}$ we identify $\rho(M^{op})$ with
$M^{op}$.

We have
\[
\left\|  \sum_{i}\rho(d)^{\ast}\lambda(a_{i})e\rho(b)^{\ast}h_{i}\right\|
^{2}\leq\left\|  \rho(d)^{\ast}\right\|  ^{2}\left\|  \sum_{i}\lambda
(a_{i})e\rho(b)^{\ast}h_{i}\right\|  ^{2},
\]
of course, and
\begin{align}
\left\|  \sum_{i}\lambda(a_{i})e\rho(b)^{\ast}h_{i}\right\|  ^{2}  &
=\sum_{i,j}\langle\lambda(a_{i})e\rho(b^{\ast})h_{i},\lambda(a_{j}%
)e\rho(b)^{\ast}h_{j}\rangle\label{DefS}\\
&  =\sum_{i,j}\langle h_{i},\rho(b)e\lambda(a_{i}^{\ast}a_{j})e\rho(b)^{\ast
}h_{j}\rangle\allowbreak\nonumber\\
&  =\sum_{i,j}\langle h_{i},\rho(b)e\lambda(\Phi(a_{i}^{\ast}a_{j}%
))e\rho(b)^{\ast}h_{j}\rangle\allowbreak.\nonumber
\end{align}
However, $\lambda(\Phi(a_{i}^{\ast}a_{j}))$ lies in $\lambda(N)$ while
$e\rho(b)^{\ast}$ and $\rho(b)e$ lie in $\lambda(N)^{\prime}$. Hence the
matrix inequality%
\[
\left(  \rho(b)e\Phi(a_{i}^{\ast}a_{j}\right)  e\rho(b^{\ast}))\leq\left\|
b\right\|  ^{2}\left(  \Phi(a_{i}^{\ast}a_{j}\right)  )
\]
is satisfied. This shows that the last expression in equation (\ref{DefS}) is
dominated by $\left\|  b\right\|  ^{2}\sum_{i,j}\langle h_{i},\Phi(a_{i}%
^{\ast}a_{j})h_{j}\rangle=\left\|  b\right\|  ^{2}\left\|  \sum_{i}%
a_{i}\otimes h_{i}\right\|  ^{2}$, which in turn shows that $S(b,d)$ is indeed
well defined and that $\left\|  S(b,d)\right\|  \leq\left\|  b\right\|
\left\|  d\right\|  $. To see that $S(b,d)^{\ast}\in E^{\lambda}$, observe
that if $c\in M$, then for $a\otimes h\in M\otimes_{\Phi}H$, we have
$S(b,d)(c\otimes I)(a\otimes h)=S(b,d)(ca\otimes h)\allowbreak=\rho(d)^{\ast
}\lambda(ca)e\rho(b)^{\ast}h=\lambda(c)\rho(d)^{\ast}\lambda(a)e\rho(b)^{\ast
}h\allowbreak=\lambda(c)S(b,d)(a\otimes h)$.

Now $\Psi:M^{op}\otimes_{\tilde{\Phi}}M^{op}\rightarrow E^{\lambda}$ is
defined by the formula $\Psi(b,d)=S(b,d)^{\ast}$, $b\otimes d\in M^{op}%
\otimes_{\tilde{\Phi}}M^{op}$. To show that $\Psi$ is well-defined and
bounded, we show that, in fact, $\Psi$ is isometric, i.e., that
\begin{align*}
\langle b_{1}\otimes d_{1},b_{2}\otimes d_{2}\rangle_{M^{op}\otimes
_{\tilde{\Phi}}M^{op}}  &  =\langle S(b_{1},d_{1})^{\ast},S(b_{2},d_{2}%
)^{\ast}\rangle\\
&  =S(b_{1},d_{1})S(b_{2},d_{2})^{\ast}\text{.}%
\end{align*}
To this end, fix $h,k\in H$ and $a\in M$. Then $\langle S(b_{2},d_{2})^{\ast
}(k),a\otimes h\rangle\allowbreak=\langle k,\rho(d_{2})^{\ast}\lambda
(a)e\rho(b_{2})^{\ast}h\rangle\allowbreak\allowbreak=\langle\rho
(d_{2})k,\lambda(a)e\rho(b_{2})^{\ast}h\rangle$. Choose a sequence $\{c_{n}\}$
in $M$ so that $\lambda(c_{n})\Omega\rightarrow\rho(d_{2})k$ in the norm of
$H$. Then
\begin{multline*}
\langle S(b_{2},d_{2})^{\ast}(k),a\otimes h\rangle=\lim_{n}\langle
\lambda(c_{n})\Omega,\lambda(a)e\rho(b_{2})^{\ast} h\rangle\\
=\lim_{n}\langle e\lambda(a^{\ast}c_{n})\Omega,\rho(b_{2})^{\ast}h\rangle
=\lim_{n}\langle e\lambda(a^{\ast}c_{n})e\Omega,e\rho(b_{2})^{\ast}h\rangle\\
=\lim_{n}\langle\lambda(\Phi(a^{\ast}c_{n}))\Omega,e\rho(b_{2})^{\ast}%
h\rangle=\lim_{n}\langle\rho(b_{2})e\Omega,\lambda(\Phi(c_{n}^{\ast
}a))h\rangle\\
=\lim_{n}\langle c_{n}\otimes\rho(b_{2})e\Omega,a\otimes h\rangle\text{.}%
\end{multline*}
Thus, for $k,g \in H$,
\begin{multline*}
\langle S(b_{1},d_{1})S(b_{2},d_{2})^{*}k,g\rangle=\langle S(b_{2},d_{2}%
)^{*}k,S(b_{1},d_{1})^{*}g\rangle\\
=lim_{n}\langle c_{n}\otimes\rho(b_{2})\Omega,S(b_{1},d_{1})^{*}%
g\rangle=lim_{n}\langle\rho(d_{1})^{*}\lambda(c_{n})e\rho(b_{1})^{*}\rho
(b_{2})\Omega,g\rangle\\
= lim_{n}\langle\rho(d_{1})^{*}\lambda(c_{n})e\rho(b_{2}b_{1}^{*}%
)\Omega,g\rangle=lim_{n}\langle\rho(d_{1})^{*}\lambda(c_{n})\rho(\tilde{\Phi
}(b_{2}b_{1}^{*}))\Omega,g\rangle\\
= lim_{n}\langle\rho(d_{1})^{*}\rho(\tilde{\Phi}(b_{2}b_{1}^{*}))\lambda
(c_{n})\Omega,g\rangle=\langle\rho(d_{1})^{*}\rho(\tilde{\Phi}(b_{2}b_{1}%
^{*}))\rho(d_{2})k,g\rangle.
\end{multline*}

Hence $S(b_{1},d_{1})S(b_{2},d_{2})^{\ast}=\rho(d_{2}\tilde{\Phi}(b_{2}%
b_{1}^{\ast})d_{1}^{\ast})$. Since $d_{2}\tilde{\Phi}(b_{2}b_{1}^{\ast}%
)d_{1}^{\ast}$ is the inner product $\langle b_{1}\otimes d_{1},b_{2}\otimes
d_{2}\rangle$ in $M^{op}\otimes_{\tilde{\Phi}}M^{op}$, this shows that $\Psi$
is isometric (and, in particular, well defined). Now compute, for
$b,d,c_{1},c_{2}$ in $M^{op}$, $a\in M$ and $h\in H$,
\begin{multline*}
S(bc_{1},c_{2}d)(a\otimes h)=\rho(c_{2}d)^{\ast}\lambda(a)e\rho(bc_{1})^{\ast
}h\\
=\rho(c_{2})^{\ast}\rho(d)^{\ast}\lambda(a)e\rho(b)^{\ast}\rho(c_{1})^{\ast
}h\\
=\rho(c_{2})^{\ast}S(b,d)(I\otimes\rho(c_{1}))^{\ast}(a\otimes h).
\end{multline*}
Thus $\Psi(bc_{1},c_{2}d)=(I\otimes\rho(c_{1}))\Psi(b,d)\rho(c_{2})$
completing the proof that $\Psi$ is an isomorphism of correspondences with
range in $E^{\lambda}$. It is left to show that $\Psi$ is surjective.

In order to show this, we claim first that%

\begin{equation}
\bigcap Ker\{S(b,I):\;b\in M^{op}\}=\{0\}. \label{ker}%
\end{equation}

For this, let us first write $S$ for $S(I,I)$ and note that $S(b,I)=S(I\otimes
\rho( b)^{*})$ for $b\in M^{op}$. From the fact that $\Psi$ is an isometry, we
know that $SS^{*}$ is the inner product of $I\otimes I$ with itself (in
$M^{op}\otimes_{\tilde{\Phi}}M^{op}$). Consequently, $SS^{*}=I$. Hence
$S^{*}S$ is a projection. Since $S(a\otimes h)=\lambda(a)eh$, it follows that
$S=S(I\otimes e)$ (note that, since $e\in\lambda(N)^{\prime}$, $I\otimes e$ is
a well defined projection on $M\otimes_{\Phi}H$). Hence $S^{*}S\leq I\otimes
e$. We also have, for $a_{1},a_{2}$ in $M$ and $h_{1},h_{2}$ in $H$,%

\[
\langle S(a_{1}\otimes h_{1}),S(a_{2}\otimes h_{2})\rangle=\langle
\lambda(a_{1})eh_{1},\lambda(a_{2})eh_{2}\rangle=
\]
\[
=\langle h_{1},e\Phi(a_{1}^{*}a_{2})eh_{2} \rangle= \langle a_{1}\otimes
eh_{1},a_{2}\otimes eh_{2}\rangle.
\]

Hence $S^{\ast}S=I\otimes e$ and, for every unitary operator $u\in M^{op}$,
\[
S(u,I)^{\ast}S(u,I)=(I\otimes\rho(u))(I\otimes e)(I\otimes\rho(u)^{\ast
})=I\otimes\rho(u)e\rho(u)^{\ast}.
\]
But $\bigvee\{\rho(u)e\rho(u)^{\ast}\}$ (where $u$ runs over all unitaries in
$M^{op}$) is $I$ since $\Omega=e\Omega$ and $\rho(M^{op})\Omega$ is dense in
$H$. Hence $\bigvee\{S(u,I)^{\ast}S(u,I)\}=I_{M\otimes H}$ and the claim
(\ref{ker}) follows.

To complete the proof of the surjectivity of $\Psi$, fix $\eta\in E^{\lambda}$
that is orthogonal to the range of $\Psi$. Then, for every $b,d \in M^{op}$,
$S(b,d)\eta=0$ and it follows from the claim that $\eta=0$.

Our last example is inspired in part by a correspondence that arises in the
theory of wavelets. See \cite{BJ02}.

\begin{example}
\label{Example2.6}Let $M$ be a von Neumann algebra and let $\sigma$ be a
representation of $M$ on a Hilbert space $H$. Assume that $\sigma(M)$ has a
cyclic vector $\Omega$. Also, let $\alpha$ be an endomorphism of $M$ with the
property that the state determined by $\Omega$ is invariant under $\alpha$.
Form the correspondence $E=M\otimes_{\alpha}M$, which as we saw in the
discussion after Example \ref{Example2.4}, is $_{\alpha}M$. Then we claim that
$E^{\sigma}=\{mS\mid m\in\sigma(\alpha(M))^{\prime}\}$, where $S$ is the
\emph{isometry} defined by the formula $S(\sigma(a)\Omega):=\sigma
(\alpha(a))\Omega$.
\end{example}

First note that $S$ is indeed an isometry: $\left\|  S(\sigma(a)\Omega
)\right\|  ^{2}=\left\|  \sigma(\alpha(a))\Omega\right\|  ^{2}=\langle
\sigma(\alpha(a^{\ast}a))\Omega,\Omega\rangle=\langle\sigma(a^{\ast}%
a)\Omega,\Omega\rangle=\left\|  \sigma(a)\Omega\right\|  ^{2}$ because of the
invariance of the state determined by $\Omega$, and because of the cyclicity
of $\Omega$, $S$ extends to all of $H$. The principal virtue of $S $, from our
perspective, is that it satisfies the equation $S\sigma(a)=\sigma(\alpha
(a))S$, which may be readily verified from its definition. So, if $\eta=mS$,
for some $m\in\sigma(\alpha(M))^{\prime}$, then $\eta\sigma(a)=mS\sigma
(a)=m\sigma(\alpha(a))S=\sigma(\alpha(a))mS=\eta\sigma(\alpha(a))$. Thus
$\eta\in E^{\sigma}$. Conversely, if $\eta\in E^{\sigma}$, then $\eta S^{\ast
}\in\sigma(\alpha(M))^{\prime}$, since $\sigma(a)S^{\ast}=S^{\ast}%
\sigma(\alpha(a))$. Consequently $\eta=\eta S^{\ast}S\in\{mS\mid m\in
\sigma(\alpha(M))^{\prime}\}$, proving our claim.

The operator $S$ is closely related to what is called a transfer operator or a
Ruelle operator in \cite{BJ02} and \cite{rEp00}. Indeed, in \cite{BJ02}, the
von Neumann algebra $M$ is $L^{\infty}(\mathbb{T)}$ acting by multiplication
operators on $L^{2}(\mathbb{T)}$. The endomorphism $\alpha$ is given by the
formula $\alpha(f)(z)=f(z^{N})$, $f\in L^{\infty}(\mathbb{T)}$, for a positive
integer $N$. The vector $\Omega$ is the constant function $1$ in
$L^{2}(\mathbb{T)}$ and the isometry $S$ is given by the formula $S\xi
(z)=\xi(z^{N})$. The adjoint of $S$ is the transfer operator.

\section{Nevanlinna-Pick Theorem}

Our goal in this section is to prove a generalization of the Nevanlinna-Pick
interpolation theorem. We begin by setting up notation. Throughout this
section, we fix a faithful normal representation $\sigma$ of $M$ on a Hilbert
space $H$. \ Since we will be considering representations of both $H^{\infty
}(E)$ and $H^{\infty}(E^{\sigma})$, we will distinguish between the various
constructs for these two algebras, when necessary, by indexing them with $E$
or $E^{\sigma}$. \ So, for example, we will write $\varphi_{\infty}^{E}$ for
the left action of $M$ on $\mathcal{F}(E)$, while we will write $\varphi
_{\infty}^{E^{\sigma}}$ for the left action of $\sigma(M)^{\prime}$ on
$\mathcal{F}(E^{\sigma})$. We write $\iota$ for the identity representation of
$\sigma(M)^{\prime}$ on $H$, but we will abbreviate the induced representation
of $\mathcal{L}(\mathcal{F}(E^{\sigma}))$ on $\mathcal{F}(E^{\sigma}%
)\otimes_{\iota}H$, $\iota^{\mathcal{F}(E^{\sigma})}$, by $\tau$. \ That is,
$\tau=\iota^{\mathcal{F}(E^{\sigma})}$, is given by the formula $\tau
(T)\zeta\otimes h=(T\zeta)\otimes h$, $\zeta\otimes h\in\mathcal{F}(E^{\sigma
})$. \ Given $\eta_{1}$ and $\eta_{2}$ in $E^{\sigma}$ with $\left\|  \eta
_{i}\right\|  <1$, $i=1,2$, we write $\theta_{\eta_{1},\eta_{2}}$ for the map
on $\sigma(M)^{\prime}$ defined by the formula
\[
\theta_{\eta_{1},\eta_{2}}(a):=\langle\eta_{1},a\eta_{2}\rangle_{M^{\prime}%
}=\eta_{1}^{\ast}(I_{E}\otimes a)\eta_{2}\text{,}%
\]
where, recall, each $\eta_{i}$ is a map from $H$ to $E\otimes_{\sigma}H$, and
the $\sigma(M)^{\prime}$-valued inner product is given by the formula
$\langle\eta_{1},\eta_{2}\rangle=\eta_{1}^{\ast}\eta_{2}$. Of course we have
the inequality
\[
\left\|  \theta_{\eta_{1},\eta_{2}}\right\|  \leq\left\|  \eta_{1}\right\|
\left\|  \eta_{2}\right\|  <1\text{.}%
\]
Consequently, we may form $(id-\theta_{\eta_{1},\eta_{2}})^{-1}=id+\theta
_{\eta_{1},\eta_{2}}+\theta_{\eta_{1},\eta_{2}}^{2}+\cdots$, obtaining a
completely bounded map on $\sigma(M)^{\prime}$ that we denote by $t_{\eta
_{1},\eta_{2}}$. When $\eta_{1}=\eta_{2}:=\eta$, $\theta_{\eta_{1},\eta_{2}}$
and $t_{\eta_{1},\eta_{2}}$ are both completely positive and we denote them by
$\theta_{\eta}$ and $t_{\eta}$, respectively.

Given operators $C_{1}$ and $C_{2}$ on the Hilbert space $H$ we shall write
$Ad(C_{1},C_{2})$ for the so-called elementary operator on $B(H)$ defined by
the formula
\[
Ad(C_{1},C_{2})(A)=C_{1}AC_{2}^{\ast}\text{,}%
\]
$A\in B(H)$.

Using Lemma \ref{Lemma1.7} and Theorem \ref{commutant}, we interchange the
roles of $E$ and $E^{\sigma}$ and let $\rho$ be the representation of
$H^{\infty}(E)$ on $\mathcal{F}(E^{\sigma})\otimes_{\iota}H$ given by the
formulae
\begin{equation}
\rho(\varphi_{\infty}^{E}(a))(\eta\otimes h)=\eta\otimes\sigma(a)h
\label{rhophiinfty}%
\end{equation}
and
\begin{equation}
\rho(T_{\xi}^{E})(\eta\otimes h)=\eta\otimes W(\xi)h\text{,} \label{rhoTxi}%
\end{equation}
where $\eta\otimes h\in\mathcal{F}(E^{\sigma})\otimes_{\iota}H$, $a\in M$,
$\xi\in E$, and where $W$ is defined by equation (\ref{wxs}) before Theorem
\ref{duality}.

By Theorem \ref{commutant}, $\rho(H^{\infty}(E))=\tau(H^{\infty}(E^{\sigma
}))^{\prime}$.

Given $\eta\in E^{\sigma}$, with $\left\|  \eta\right\|  \leq1$, we define a
covariant representation of $E$ that we denote by $(S(\eta),\sigma)$ through
the formula%
\[
S(\eta)(\xi)h:=\eta^{\ast}(\xi\otimes h)\text{,}%
\]
$\xi\in E$, $h\in H$. That is, $\widetilde{S(\eta)}=\eta^{\ast}$. The
integrated form of this representation, $\sigma\times S(\eta)$ gives a
representation of $\mathcal{T}_{+}(E)$, and if $\left\|  \eta\right\|  \lneq
1$, then, as we saw in Corollary \ref{Corollary1.12}, $\sigma\times S(\eta)$
extends to all of $H^{\infty}(E)$. In this event, for $X\in H^{\infty}(E)$, we
write $X(\eta^{\ast})$ for $\sigma\times S(\eta)(X)$. Thus we like to think of
the open unit ball in $(E^{\sigma})^{\ast}$ as a generalization of the open
unit disc, $\mathbb{D}$, in the complex plane and for $\eta^{\ast}$ in the
open unit ball of $(E^{\sigma})^{\ast}$ and $X\in H^{\infty}(E)$,
$X(\eta^{\ast})$ is the value of $X$ at $\eta^{\ast}$. \ Note that when $M$,
$H$ and $E$ are all one-dimensional, so is $E^{\sigma}$, $H^{\infty}(E)$ is
the space $H^{\infty}(\mathbb{T})$, and for $X\in H^{\infty}(E)=H^{\infty
}(\mathbb{T})$, $X(\eta^{\ast})=X(\overline{\eta})$ is the value of $X$ at
$\bar{\eta}$ in the ordinary sense. For this reason, we shall write
$\mathbb{D((}E^{\sigma})^{\ast})$ for the collection of all $\eta^{\ast}%
\in(E^{\sigma})^{\ast}$ with norm strictly less than one. \ We will give other
examples of how this formula arises in the next section. \ We recommend, in
particular, Example \ref{Example4.7}.

First, we want to note that in general, $H^{\infty}(E)$ does \emph{not}
separate the points of $\mathbb{D((}E^{\sigma})^{\ast})$. The simplest
example, perhaps, is afforded by Proposition 2.4 in Davidson's and Pitts's
article \cite{DP98a}. In our terminology, if $M$ is $\mathbb{C}$, if
$E=\mathbb{C}^{n}$, $n>1$ and if $\sigma$ is the identity representation, then
$H^{\infty}(E)$ is their free semigroup algebra and their proposition states
that the collection of all elements of $H^{\infty}(E)$ that vanish on
$\mathbb{D((}E^{\sigma})^{\ast})$ is the ultraweak closure of the commutator
ideal in $H^{\infty}(E)$. What the situation is in other settings remains mysterious.

The value of $X$ at an $\eta^{\ast}\in\mathbb{D((}E^{\sigma})^{\ast})$ can
also be calculated in terms of the representation $\rho$ and a map $L_{\eta}$
from $H$ to $\mathcal{F}(E^{\sigma})\otimes_{i}H$ that hints of the Cauchy
kernel. \ For this reason, we call it the \emph{Cauchy transform. }It is given
by the formula%
\[
L_{\eta}h=h\oplus\eta\otimes h\oplus\eta^{\otimes2}\otimes h\oplus
\cdots\text{,}%
\]
$h\in H$, where for $n>0$, $\eta^{\otimes n}=\eta\otimes\eta\otimes
\cdots\otimes\eta$, $n$ times and where $h$ should be viewed as $1_{\sigma
(M)^{\prime}}\otimes h$. (That is, we identify $\sigma(M)^{\prime}%
\otimes_{\iota}H$ with $H$ and, likewise, $M\otimes_{\sigma}H$ is identified
with $H$.) Observe that $L_{\eta}$ is, indeed, a bounded operator, with norm
dominated by $\frac{1}{1-\left\|  \eta\right\|  }$. \ In the sequel, we will
have to consider the natural imbedding $\iota_{H}$ from $H$ into
$\mathcal{F}(E^{\sigma})\otimes_{\iota}H$ that is given by the formula
$\iota_{H}(h)=h\oplus0\oplus0\oplus\cdots$. \ Abusing notation slightly, we
shall write $P_{H}$ for the adjoint of $\iota_{H}$.

\begin{proposition}
\label{CauchyInt}\emph{If }$\eta^{\ast}\in\mathbb{D((}E^{\sigma})^{\ast})$,
\emph{then }%
\begin{equation}
X(\eta^{\ast})=L_{\eta}^{\ast}\rho(X)\iota_{H} \label{Evaluate}%
\end{equation}
\emph{ for all }$X\in H^{\infty}(E)$\emph{.\medskip}
\end{proposition}

\begin{proof}
It clearly suffices to prove the formula for the two cases when $X=\varphi
_{\infty}^{E}(a)$, $a\in M$, and when $X=T_{\xi}^{E}$, $\xi\in E^{\otimes n}$,
$n>0$. \ Note first that if $\zeta\otimes h\in(E^{\sigma})^{\otimes n}%
\otimes_{\iota}H$, regarded as a summand of $\mathcal{F}(E^{\sigma}%
)\otimes_{\iota}H$, then $L_{\eta}^{\ast}(\zeta\otimes h)=L_{\eta}^{\ast
}(0\oplus\cdots0\oplus\zeta\otimes h\oplus0\oplus\cdots)=\langle\eta^{\otimes
n},\zeta\rangle h$. Suppose that $X=\varphi_{\infty}^{E}(a)$. \ Then
$X(\eta^{\ast})=\sigma(a)$ and for every $h\in H$,
\begin{multline*}
L_{\eta}^{\ast}\rho(\varphi_{\infty}^{E}(a))\iota_{H}(h)=L_{\eta}^{\ast}%
\rho(\varphi_{\infty}^{E}(a))(h\oplus0\oplus\cdots)\\
=L_{\eta}^{\ast}(\sigma(a)h\oplus0\oplus\cdots)=\sigma(a)h\text{.}%
\end{multline*}
\ Thus equation (\ref{Evaluate}) holds when $X=\varphi_{\infty}^{E}(a)$. \ 

To verify equation (\ref{Evaluate}) when $X=T_{\xi}^{E}$, $\xi\in E^{\otimes
n}$, we need two preliminary calculations. For the first, let $y$ be in
$E^{\sigma}\otimes_{\iota}H$ and let $\theta$ be in $(E^{\sigma})^{\otimes k}%
$, then $\theta\otimes y$ lies in $(E^{\sigma})^{\otimes k+1}\otimes_{\iota}%
H$, of course, and when $y$ and $\theta\otimes y$ are viewed as elements of
$\mathcal{F}(E^{\sigma})\otimes_{\iota}H$,
\[
L_{\eta}^{\ast}(\theta\otimes y)=L_{\eta}^{\ast}(\varphi_{\infty}^{E^{\sigma}%
}(\langle\eta^{\otimes k},\theta\rangle)\otimes I_{H})y\text{.}%
\]
To see this, assume, as we may, that $y$ has the form $y=\zeta\otimes g\in
E^{\sigma}\otimes_{\iota}H$. Then
\begin{multline*}
L_{\eta}^{\ast}(\varphi_{\infty}^{E^{\sigma}}(\langle\eta^{\otimes k}%
,\theta\rangle)\otimes I_{H})(\zeta\otimes g)=L_{\eta}^{\ast}(\varphi_{\infty
}^{E^{\sigma}}(\langle\eta^{\otimes k},\theta\rangle)\zeta\otimes g)=L_{\eta
}^{\ast}(\langle\eta^{\otimes k},\theta\rangle\zeta\otimes g)\\
=\langle\eta,\langle\eta^{\otimes k},\theta\rangle\zeta\rangle g=\langle
\eta\otimes\eta^{\otimes k},\theta\otimes\zeta\rangle g=\langle\eta^{\otimes
k+1},\theta\otimes\zeta\rangle g\\
=L_{\eta}^{\ast}(\theta\otimes\zeta\otimes g)\text{,}%
\end{multline*}
as asserted.

For the second, we claim that if $\xi\in E$, then
\begin{equation}
L_{\eta}^{\ast}\rho(T_{\xi}^{E})=L_{\eta}^{\ast}\rho(T_{\xi}^{E})\iota
_{H}L_{\eta}^{\ast}. \label{letarhotxiequal}%
\end{equation}
To see this, recall the definition of $\rho(T_{\xi}^{E})$ in equation
(\ref{rhoTxi}) and the fact that $W(\xi)$, defined in equation (\ref{wxs}),
maps $H$ to $E^{\sigma}\otimes_{\iota}H$. Thus, if $\theta\otimes k$ lies in
$(E^{\sigma})^{\otimes k}\otimes_{\iota}H$, then $L_{\eta}^{\ast}\rho(T_{\xi
}^{E})(\theta\otimes k)=L_{\eta}^{\ast}(\theta\otimes W(\xi)k)$. By the
calculation of the preceding paragraph, $L_{\eta}^{\ast}(\theta\otimes
W(\xi)k)=L_{\eta}^{\ast}(\varphi_{\infty}^{E^{\sigma}}(\langle\eta^{\otimes
k},\theta\rangle)\otimes I_{H})(W(\xi)k)$. However, by the covariance
properties of $W(\xi)$, this expression is $L_{\eta}^{\ast}(W(\xi)\langle
\eta^{\otimes k},\theta\rangle k)$. \ Since $\langle\eta^{\otimes k}%
,\theta\rangle k=L_{\eta}^{\ast}(\theta\otimes k)$, $L_{\eta}^{\ast}%
(W(\xi)\langle\eta^{\otimes k},\theta\rangle k)=L_{\eta}^{\ast}(W(\xi)L_{\eta
}^{\ast}(\theta\otimes k))=L_{\eta}^{\ast}(\rho(T_{\xi}^{E})\iota_{H}L_{\eta
}^{\ast})(\theta\otimes k)$. This proves that $L_{\eta}^{\ast}\rho(T_{\xi}%
^{E})=L_{\eta}^{\ast}\rho(T_{\xi}^{E})\iota_{H}L_{\eta}^{\ast}$ on each
summand $(E^{\sigma})^{\otimes k}\otimes_{i}H$ of $\mathcal{F}(E^{\sigma
})\otimes_{i}H$, and so $L_{\eta}^{\ast}\rho(T_{\xi}^{E})=L_{\eta}^{\ast}%
\rho(T_{\xi}^{E})\iota_{H}L_{\eta}^{\ast}$ on all of $\mathcal{F}(E^{\sigma
})\otimes_{i}H$.

Next observe that if $\xi\in E$ and if $h$ and $k$ are in $H$, then
\begin{align*}
\langle L_{\eta}^{\ast}\rho(T_{\xi})\iota_{H}h,k\rangle=\langle W(\xi
)h,L_{\eta}k\rangle=\langle W(\xi)h,\eta\otimes k\rangle\\
=\langle\xi\otimes h,\eta k\rangle=\langle\eta^{\ast}(\xi\otimes h),k\rangle\\
=\langle T_{\xi}(\eta^{\ast})h,k\rangle\text{.}%
\end{align*}
That is, equation (\ref{Evaluate}) is satisfied when $X$ is of the form
$T_{\xi}$, $\xi\in E$. Since the map $X\rightarrow X(\eta^{\ast})$ is
multiplicative, we conclude from equation (\ref{letarhotxiequal}) that if
$\xi=\xi_{1}\otimes\xi_{2}\otimes\cdots\otimes\xi_{k}$ is in $E^{\otimes k},$
then
\begin{align*}
L_{\eta}^{\ast}\rho(T_{\xi})\iota_{H}=L_{\eta}^{\ast}\rho(T_{\xi_{1}}%
)\rho(T_{\xi_{2}})\cdots\rho(T_{\xi_{k}})\iota_{H}\\
=L_{\eta}^{\ast}\rho(T_{\xi_{1}})\iota_{H}L_{\eta}^{\ast}\rho(T_{\xi_{2}%
})\cdots\rho(T_{\xi_{k}})\iota_{H}\\
=L_{\eta}^{\ast}\rho(T_{\xi_{1}})\iota_{H}L_{\eta}^{\ast}\rho(T_{\xi_{2}%
})\iota_{H}L_{\eta}^{\ast}\cdots\rho(T_{\xi_{k}})\iota_{H}\\
\cdots\\
=L_{\eta}^{\ast}\rho(T_{\xi_{1}})\iota_{H}L_{\eta}^{\ast}\rho(T_{\xi_{2}%
})\iota_{H}L_{\eta}^{\ast}\cdots L_{\eta}^{\ast}\rho(T_{\xi_{k}})\iota
_{H}L_{\eta}^{\ast}\iota_{H}\\
=L_{\eta}^{\ast}\rho(T_{\xi_{1}})\iota_{H}L_{\eta}^{\ast}\rho(T_{\xi_{2}%
})\iota_{H}L_{\eta}^{\ast}\cdots L_{\eta}^{\ast}\rho(T_{\xi_{k}})\iota_{H}\\
=T_{\xi_{1}}(\eta^{\ast})T_{\xi_{2}}(\eta^{\ast})\cdots T_{\xi_{k}}(\eta
^{\ast})\\
=T_{\xi}(\eta^{\ast})\text{.}%
\end{align*}
Taking linear combinations and limits, we conclude that equation
(\ref{Evaluate}) is valid.
\end{proof}

\begin{remark}
\label{letarhotxiequalbis}It will be important for a later calculation to
observe that the argument just given establishes this generalization of
equation (\ref{letarhotxiequal}):
\[
L_{\eta}^{\ast}\rho(X)=L_{\eta}^{\ast}\rho(X)\iota_{H}L_{\eta}^{\ast},
\]
for all $X\in H^{\infty}(E)$.
\end{remark}

The following is our generalization of the Nevanlinna-Pick interpolation Theorem.

\begin{theorem}
\label{Theorem3.1}Let $E$ be a $W^{\ast}$-correspondence over the von Neumann
algebra $M$ and let $\sigma$ be a representation of $M$ on $H$. Given
$\eta_{1}^{\ast},\eta_{2}^{\ast},\cdots,\eta_{k}^{\ast}\in\mathbb{D((}%
E^{\sigma})^{\ast})$ and two $k$-tuples of operators in $B(H)$, $B_{1}%
,B_{2},\cdots,B_{k}$ and $C_{1},C_{2},\cdots,C_{k}$, then there is an element
$X\in H^{\infty}(E)$ such that $\left\|  X\right\|  \leq1$ and such that%
\begin{equation}
B_{i}X(\eta_{i}^{\ast})=C_{i}\text{,} \label{Interpequat}%
\end{equation}
$i=1,2,\cdots k$, if and only if the map from $M_{k}(\sigma(M)^{\prime})$ to
$M_{k}(B(H))$ defined by the $k\times k$ matrix in $M_{k}(B(\sigma(M)^{\prime
},B(H)))$,%
\begin{equation}
\left(  \left(  Ad(B_{i},B_{j})-Ad(C_{i},C_{j})\right)  \circ\left(
id-\theta_{\eta_{i},\eta_{j}}\right)  ^{-1}\right)  \label{positequat}%
\end{equation}
is completely positive.

Moreover, if $N$ is a von Neumann subalgebra of $M$ and if for each $i$,
$i=1,2,\cdots,k$, $C_{i}$ and $B_{i}$ lie in the von Neumann algebra generated
by $\sigma(N)$ and $\sigma(M)^{\prime}$, $\sigma(N)\bigvee\sigma(M)^{\prime}$,
then whenever the matrix (\ref{positequat}) represents a completely positive
operator, a solution $X$ to the interpolation equation (\ref{Interpequat}) can
be found that commutes with $\varphi_{\infty}(M\cap N^{\prime})^{\prime}$.
\end{theorem}

\begin{remark}
\label{Remark3.2}If $M$ is $B(H)$ and if $\sigma$ is the identity
representation, then $\sigma(M)^{\prime}=\mathbb{C}$ and the matrix
(\ref{positequat}) really defines a map on $M_{k}(\mathbb{C)}$. We shall
denote this map by $\Psi$. To tell if $\Psi$ is completely positive, one can
invoke a theorem of Choi \cite{mC75}, which asserts that a map $\Phi$ from
$M_{k}(\mathbb{C)}$ into $B(K)$, for some Hilbert space $K$, is completely
positive if and only if the matrix of operators $(\Phi(e_{ij}))$ is a positive
element in $M_{k}(B(K)\mathbb{)}$, where $\{e_{ij}\}_{i,j=1}^{k}$ are the
matrix units for $M_{k}(\mathbb{C})$. So, if we form $(\Psi(e_{ij}))$, we get
the matrix whose $i,j$ entry is the $k\times k$ matrix in $M_{k}(B(H))$ whose
$i,j$ entry is
\[
\left(  Ad(B_{i},B_{j})-Ad(C_{i},C_{j})\right)  \circ\left(  id-\theta
_{\eta_{1},\eta_{2}}\right)  ^{-1}(1)=\frac{1}{1-\langle\eta_{i},\eta
_{j}\rangle}(B_{i}B_{j}^{\ast}-C_{i}C_{j}^{\ast})
\]
and whose other entries are zero. \ It follows that when $M=B(H)$ and $\sigma$
is the identity representation, the complete positivity of the map represented
by the matrix (\ref{positequat}) is equivalent to the positivity of the
operator matrix%
\begin{equation}
\left(  \frac{B_{i}B_{j}^{\ast}-C_{i}C_{j}^{\ast}}{1-\langle\eta_{i},\eta
_{j}\rangle}\right)  \text{.} \label{MatricialNPThm}%
\end{equation}
In particular, if $M$, $E$ and $H$ are all one dimensional, if we take
$B_{i}=1$, $i=1,2,\cdots k$, and if we write $z_{i}$ for $\eta_{i}^{\ast}$ and
$w_{i}$ for $C_{i}$, we see that the theorem asserts that we can find a
bounded analytic function $f$ of \textrm{sup-}norm at most one such that
$f(z_{i})=w_{i}$, for all $i$, if and only if the $k\times k$ matrix%
\[
\left(  \frac{1-w_{i}\bar{w}_{j}}{1-z_{i}\bar{z}_{j}}\right)
\]
is non-negative. \ This, of course, is the classical Nevanlinna-Pick
interpolation theorem.
\end{remark}

\begin{proof}
(of Theorem \ref{Theorem3.1}) The proof we give is a descendant of Sarason's
proof of the classical Nevanlinna-Pick theorem \cite{dS67} because it rests on
the commutant lifting theorem. Our organization of the proof was inspired by
the presentation in \cite{RR97}.

Suppose the map represented by the matrix (\ref{positequat}) is completely
positive. \ Recall that we are writing $\tau$ for the induced representation
$\iota^{\mathcal{F}(E^{\sigma})}$ and for each $i$, $1\leq i\leq k$, write
$L_{i}$ for the Cauchy transform $L_{\eta_{i}}$ mapping from $H$ to
$\mathcal{F}(E^{\sigma})\otimes_{\iota}H$. Then for $\eta\in E^{\sigma}$ and
$h\in H$, we have
\begin{align}
\tau(T_{\eta})^{\ast}L_{i}h=\langle\eta,\eta_{i}\rangle h+\varphi(\langle
\eta,\eta_{i}\rangle)\eta_{i}\otimes h+\cdots\label{invariant}\\
=\tau(\varphi_{\infty}(\langle\eta,\eta_{i}\rangle))L_{i}h\text{,}\nonumber
\end{align}
which lies in $\tau(\varphi_{\infty}(M^{\prime}))L_{i}(H)$. Let $\mathfrak
{M}_{i}$ be the closure of $\tau(\varphi_{\infty}(M^{\prime}))L_{i}%
(B_{i}^{\ast}H)$ and set $\mathfrak{M}=\overline{\bigvee\mathfrak{M}_{i}}$. We
shall show that if the matrix in (\ref{positequat}) is completely positive,
then there is an operator $R:\mathfrak{M}\rightarrow\mathcal{F}(E^{\sigma
})\otimes_{\iota}H$ of norm at most one satisfying the equation%
\begin{equation}
R\tau(\varphi_{\infty}(a))L_{i}B_{i}^{\ast}h=\tau(\varphi_{\infty}%
(a))L_{i}C_{i}^{\ast}h \label{2'}%
\end{equation}
for all $a\in\sigma(M)^{\prime}$ and all $h\in H$.

A general element of $\mathfrak{M}$ is a limit of sums of the form $\sum
\tau(\varphi_{\infty}(a_{i,j}))L_{i}B_{i}^{\ast}h_{i,j}$, where $i$ runs from
$1$ to $k$, $j$ runs over a finite set of indices, the $a_{i,j}$ all lie in
$\sigma(M)^{\prime}$ and the $h_{i,j}$ are in $H$. So a contractive operator
$R:\mathfrak{M}\rightarrow H$ exists satisfying equation (\ref{2'}) if and
only if for all such sums the inequality%
\begin{equation}
\left\|  \sum\tau(\varphi_{\infty}(a_{i,j}))L_{i}C_{i}^{\ast}h_{i,j}\right\|
^{2}\leq\left\|  \sum\tau(\varphi_{\infty}(a_{i,j}))L_{i}B_{i}^{\ast}%
h_{i,j}\right\|  ^{2} \label{i}%
\end{equation}
is satisfied. So, if we expand this inequality in terms of inner products, we
will see how the complete positivity of the matrix (\ref{positequat})
intervenes. To this end, we require a preliminary calculation.

For $a,b\in\sigma(M)^{\prime}$, and $h,k\in H$, we find that for all $i$ and
$m$, $1\leq i,m\leq k$,%
\begin{multline*}
\langle\tau(\varphi_{\infty}(a))L_{i}h,\tau(\varphi_{\infty}(b))L_{m}%
k\rangle\\
=\sum_{j,p}\langle a\overset{j\;\mathrm{times}}{\overbrace{\eta_{i}\otimes
\eta_{i}\cdots\otimes\eta_{i}}}\otimes h,b\overset{p\;\mathrm{times}%
}{\overbrace{\eta_{m}\otimes\eta_{m}\cdots\otimes\eta_{m}}}\otimes k\rangle\\
=\sum_{j=0}^{\infty}\langle h,\langle\overset{j\;\mathrm{times}}%
{\overbrace{\eta_{i}\otimes\eta_{i}\cdots\otimes\eta_{i}}},a^{\ast}%
b\overset{j\;\mathrm{times}}{\overbrace{\eta_{m}\otimes\eta_{m}\cdots
\otimes\eta_{m}}}\rangle_{(E^{\sigma})^{\otimes j}}k\rangle\\
=\langle h,\sum_{j=0}^{\infty}\theta_{\eta_{i},\eta_{m}}^{j}(a^{\ast
}b)k\rangle\\
=\langle h,t_{\eta_{i},\eta_{m}}(a^{\ast}b)k\rangle\text{,}%
\end{multline*}
where, recall, $\theta_{\eta_{i},\eta_{m}}(a)=\langle\eta_{i},a\eta_{m}%
\rangle$ and where $t_{\eta_{i},\eta_{m}}$ is the completely bounded map
$\sum_{j=0}^{\infty}\theta_{\eta_{i},\eta_{m}}^{j}$.

So, we may write%
\begin{multline*}
\left\|  \sum\tau(\varphi_{\infty}(a_{i,j}))L_{i}C_{i}^{\ast}h_{i,j}\right\|
^{2}\\
=\sum_{(i,j),(m,p)}\langle\tau(\varphi_{\infty}(a_{i,j}))L_{i}C_{i}^{\ast
}h_{i,j},\tau(\varphi_{\infty}(a_{m,p}))L_{m}C_{m}^{\ast}h_{m,p}\rangle\\
=\sum_{(i,j),(m,p)}\langle C_{i}^{\ast}h_{i,j},t_{\eta_{i},\eta_{m}}%
(a_{i,j}^{\ast}a_{m,p})C_{m}^{\ast}h_{m,p}\rangle\\
=\sum_{(i,j),(m,p)}\langle h_{i,j},C_{i}t_{\eta_{i},\eta_{m}}(a_{i,j}^{\ast
}a_{m,p})C_{m}^{\ast}h_{m,p}\rangle\text{.}%
\end{multline*}
In a like manner, we may write%
\[
\left\|  \sum\tau(\varphi_{\infty}(a_{i,j}))L_{i}B_{i}^{\ast}h_{i,j}\right\|
^{2}=\sum_{(i,j),(m,p)}\langle h_{i,j},B_{i}t_{\eta_{i},\eta_{m}}%
(a_{i,j}^{\ast}a_{m,p})B_{m}^{\ast}h_{m,p}\rangle\text{.}%
\]
So, to verify the inequality (\ref{i}), we need to show that the difference
\[
\sum_{(i,j),(m,p)}\langle h_{i,j},B_{i}t_{\eta_{i},\eta_{m}}(a_{i,j}^{\ast
}a_{m,p})B_{m}^{\ast}h_{m,p}\rangle-\sum_{(i,j),(m,p)}\langle h_{i,j}%
,C_{i}t_{\eta_{i},\eta_{m}}(a_{i,j}^{\ast}a_{m,p})C_{m}^{\ast}h_{m,p}\rangle
\]
is non negative. However, this difference is%
\begin{align}
&  \sum_{(i,j),(m,p)}\langle h_{i,j},B_{i}t_{\eta_{i},\eta_{m}}(a_{i,j}^{\ast
}a_{m,p})B_{m}^{\ast}h_{m,p}-C_{i}t_{\eta_{i},\eta_{m}}(a_{i,j}^{\ast}%
a_{m,p})C_{m}^{\ast}h_{m,p}\rangle\label{i'}\\
&  =\sum_{(i,j),(m,p)}\langle h_{i,j},\Psi_{i,m}(a_{i,j}^{\ast}a_{m,p}%
)h_{m,p}\rangle\text{,}\nonumber
\end{align}
where%
\[
\Psi_{i,m}=\left(  Ad(B_{i},B_{m})-Ad(C_{i},C_{m})\right)  \circ\left(
id-\theta_{\eta_{i},\eta_{m}}\right)  ^{-1}\text{.}%
\]
Our assumption is that the map represented by the matrix (\ref{positequat}) is
completely positive. This means, then, that the operator matrix $(\Psi
_{i,m}(a_{i,j}^{\ast}a_{m,p}))$ is positive and, therefore, that the quantity
in equation (\ref{i'}) is nonnegative, which is what we were out to prove.
Thus there is indeed a contractive map $R:\mathfrak{M}\rightarrow
\mathcal{F}(E^{\sigma})\otimes_{\iota}H$ satisfying equation (\ref{2'}).

Observe that by definition, in fact, $R$ carries $\mathfrak{M}$ into
$\mathfrak{M}$. On the other hand, by equation (\ref{invariant}),
$\tau(\mathcal{T}_{+}(E^{\sigma}))^{\ast}$ leaves $\mathfrak{M}$ invariant. We
claim that $R$ commutes with the restriction of $\tau(\mathcal{T}%
_{+}(E^{\sigma}))^{\ast}$ to $\mathfrak{M}$. To see this, fix an $\eta\in
E^{\sigma}$ and a sum of the form $\sum\tau(\varphi_{\infty}(a_{i,j}%
))L_{i}B_{i}^{\ast}h_{i,j}$ in $\mathfrak{M}$. Then%
\begin{multline*}
\tau(T_{\eta})^{\ast}R(\sum_{(i,j)}\tau(\varphi_{\infty}(a_{i,j}))L_{i}%
B_{i}^{\ast}h_{i,j})=\tau(T_{\eta})^{\ast}\sum_{(i,j)}\tau(\varphi_{\infty
}(a_{i,j}))L_{i}C_{i}^{\ast}h_{i,j}\\
=\sum_{\substack{(i,j)\\p}}\tau(T_{\eta})^{\ast}\tau(\varphi_{\infty}%
(a_{i,j}))\eta_{i}^{\otimes p}\otimes C_{i}^{\ast}h_{i,j}\\
=\sum_{\substack{(i,j)\\p}}\tau(\varphi_{\infty}(\langle\eta,a_{i,j}\eta
_{i}\rangle)\eta_{i}^{\otimes(p-1)}\otimes C_{i}^{\ast}h_{i,j}\\
=\sum_{(i,j)}\tau(\varphi_{\infty}(\langle\eta,a_{i,j}\eta_{i}\rangle
)L_{i}C_{i}^{\ast}h_{i,j}\\
=R\sum_{(i,j)}\tau(\varphi_{\infty}(\langle\eta,a_{i,j}\eta_{i}\rangle
)L_{i}B_{i}^{\ast}h_{i,j}\\
=R\tau(T_{\eta})^{\ast}\sum_{(i,j)}\tau(\varphi_{\infty}(a_{i,j}))L_{i}%
B_{i}^{\ast}h_{i,j}\text{.}%
\end{multline*}
Hence $R\tau(T_{\eta})^{\ast}P_{\mathfrak{M}}=\tau(T_{\eta})^{\ast
}RP_{\mathfrak{M}}$ and, thus, $P_{\mathfrak{M}}\tau(T_{\eta})R^{\ast
}=P_{\mathfrak{M}}R^{\ast}\tau(T_{\eta})$ for all $\eta\in E^{\sigma}$. Since
$R$ commutes with the action of $\tau(\varphi_{\infty}(M^{\prime}))$ on
$\mathfrak{M}$, we conclude that
\[
P_{\mathfrak{M}}\tau(X)R^{\ast}=P_{\mathfrak{M}}R^{\ast}\tau(X)
\]
for all $X\in H^{\infty}(E^{\sigma})$.

At this point, we could apply the commutant lifting theorem \cite[Theorem
4.4]{MS98} to obtain an operator on all of $\mathcal{F}(E^{\sigma}%
)\otimes_{\iota}H$ that dilates $R^{\ast}$ and commutes with $\tau(H^{\infty
}(E))$ to complete the proof, but before doing so, we need to attend to the
details of what happens if we assume that the $C_{i}$ and $B_{i}$ all lie in
$\sigma(N)\bigvee\sigma(M)^{\prime}$, where, recall, $N$ is a von Neumann
algebra contained in $M$. Note that $(\sigma(N)\bigvee\sigma(M)^{\prime
})^{\prime}=\sigma(N)^{\prime}\bigwedge\sigma(M)$, so we may define a faithful
normal representation $\pi$ of $(\sigma(N)\bigvee\sigma(M)^{\prime})^{\prime}$
on $\mathcal{F}(E^{\sigma})\otimes_{\iota}H$ by the formula $\pi
(a):=I\otimes\sigma(a)$, $a\in(\sigma(N)\bigvee\sigma(M)^{\prime})^{\prime}$.
\ It is evident from the definitions of $\pi$ and $R$ that $R$ commutes with
the von Neumann algebra $\pi((\sigma(N)\bigvee\sigma(M)^{\prime})^{\prime})$.
Observe also that from its definition, it is easy to see that $\mathfrak{M}$
is invariant under $\pi((\sigma(N)\bigvee\sigma(M)^{\prime})^{\prime})$,
assuming, as we are, that the $C_{i}$ and $B_{i}$ all lie in $\sigma
(N)\bigvee\sigma(M)^{\prime}$. We let $\tilde{M}$ denote the von Neumann
subalgebra of $B(\mathcal{F}(E^{\sigma})\otimes_{\iota}H)$ generated by
$\tau(\varphi_{\infty}(\sigma(M)^{\prime}))$ and $\pi((\sigma(N)\bigvee
\sigma(M)^{\prime})^{\prime})$ and we let $F$ denote the ultraweak closure of
the span $\{\tau(T_{\eta})\pi(a)\mid a\in(\sigma(N)\bigvee\sigma(M)^{\prime
})^{\prime}$, $\eta\in E^{\sigma}\}$. Then since $\pi((\sigma(N)\bigvee
\sigma(M)^{\prime})^{\prime})$ commutes with all $\tau(T_{\eta})$, $\eta\in
E^{\sigma}$, $F$ is a bimodule over $\tilde{M}$. Also, $F$ is naturally a
$W^{\ast}$-correspondence over $\tilde{M},$ via the $\tilde{M}$-valued inner
product given on generators by the formula%
\begin{align*}
\langle\tau(T_{\eta_{1}})\pi(a_{1}),\tau(T_{\eta_{2}})\pi(a_{2})\rangle
:=(\tau(T_{\eta_{1}})\pi(a_{1}))^{\ast}\tau(T_{\eta_{2}})\pi(a_{2})\\
=\tau(\varphi_{\infty}(\langle\eta_{1},\eta_{2}\rangle))\pi(a_{1}^{\ast}%
a_{2})\text{,}%
\end{align*}
$a_{i}\in\tilde{M}$ and $\eta_{i}\in E^{\sigma}$, $i=1,2$. Since $\mathfrak
{M}$ reduces $\pi$ and since $T_{\eta}^{\ast}$ leaves $\mathfrak{M}$ invariant
for all $\eta\in E^{\sigma}$, we conclude that compressing $F$ and $\tilde{M}$
to $\mathfrak{M}$ yields a completely contractive covariant representation of
the pair $(F,\tilde{M})$. Also, of course, $R^{\ast}$ commutes with the range
of this representation. Since the identity representation of $(F,\tilde{M})$
on $\mathcal{F}(E^{\sigma})\otimes_{\iota}H$ is an isometric dilation of the
representation on $\mathfrak{M}$ (in the sense of \cite[Definition 3.1]%
{MS98}), we may apply Theorem 4.4 of \cite{MS98} to conclude that there is an
operator $Y\in B(\mathcal{F}(E^{\sigma})\otimes_{\iota}H)$ with norm at most
one such that $Y^{\ast}$ leaves $\mathfrak{M}$ invariant, $P_{\mathfrak{M}%
}Y|\mathfrak{M}=R^{\ast}$, and $Y$ commutes with $\tilde{M}$ and with $F$.
Hence, $Y$ commutes with $\tau(\mathcal{T}_{+}(E^{\sigma}))$ and with
$\pi((\sigma(N)\bigvee\sigma(M)^{\prime})^{\prime})$. By Theorem
\ref{commutant}, we conclude that there is an $X\in H^{\infty}(E)$, $\left\|
X\right\|  \leq1$, such that $\rho(X)=Y.$ To show that $X$ commutes with
$\varphi_{\infty}(M\wedge N^{\prime})$ in $\mathcal{L}(\mathcal{F}(E))$, we
note that if $U$ is the unitary operator from $\mathcal{F}(E^{\sigma}%
)\otimes_{\iota}H$ to $\mathcal{F}(E)\otimes_{\sigma}H$ defined in Lemma
\ref{Lemma1.7}, then $U\rho(X)U^{-1}=\sigma^{\mathcal{F}(E)}(X)$, while
$U\pi(a)U^{-1}=\sigma^{\mathcal{F}(E)}(\varphi_{\infty}(a))$. Since $\rho(X)$
commutes with $\pi((\sigma(N)\bigvee\sigma(M)^{\prime})^{\prime})$, we see
that $\sigma^{\mathcal{F}(E)}(X)$ commutes with $\sigma^{\mathcal{F}%
(E)}(\varphi_{\infty}((M\wedge N^{\prime}))$. \ But then, of course, $X$
commutes with $\varphi_{\infty}(M\wedge N^{\prime})$.

Our goal now is to show that $X$ satisfies equation (\ref{Interpequat}). By
virtue of the fact that $\rho(X)^{\ast}=Y^{\ast}$ leaves $\mathfrak{M}$
invariant and $\rho(X)^{\ast}|\mathfrak{M}=R^{\ast}$, we conclude from the
definition of $R$ and $\mathfrak{M}$ that
\begin{equation}
\rho(X)^{\ast}L_{i}B_{i}^{\ast}h=L_{i}C_{i}^{\ast}h \label{*}%
\end{equation}
for all $h\in H$ and for all $i$, $1\leq i\leq k$. Observe that by definition
of $L_{i}$, $P_{H}L_{i}h=h$ for all $h\in H$ where, recall, we are viewing the
projection of $\mathcal{F}(E^{\sigma})\otimes_{\iota}H$ onto the zeroth
summand as the adjoint of the imbedding, $\iota_{H}$, of $H$ into
$\mathcal{F}(E^{\sigma})\otimes_{\iota}H$. Hence equation (\ref{*}) gives
$C_{i}^{\ast}h=P_{H}L_{i}C_{i}^{\ast}h=P_{H}\rho(X)^{\ast}L_{i}B_{i}^{\ast
}h=(L_{i}^{\ast}\rho(X)\iota_{H})^{\ast}B_{i}^{\ast}h$ for all $h\in H$ and
all $i$, $1\leq i\leq k$. That is, we have $C_{i}=B_{i}(L_{i}^{\ast}%
\rho(X)\iota_{H})$ for all $i$, $1\leq i\leq k$. However, by Proposition
\ref{CauchyInt}, $L_{i}^{\ast}\rho(X)\iota_{H}=X(\eta_{i}^{\ast})$. Thus, we
find that equation (\ref{Interpequat}) is satisfied, which completes one half
of the proof.

For the converse direction, suppose $X$ is an element of $H^{\infty}(E)$ with
norm at most $1$ that satisfies the equations $B_{i}X(\eta_{i}^{\ast})=C_{i}$,
$i=1,2\ldots k.$ Then, by Proposition \ref{CauchyInt}, $B_{i}L_{i}^{\ast}%
\rho(X)\iota_{H}=C_{i}$, $i=1,2\ldots k.$ Multiply these equations on the
right by $L_{i}^{\ast}$ to conclude that $B_{i}L_{i}^{\ast}\rho(X)\iota
_{H}L_{i}^{\ast}=C_{i}L_{i}^{\ast}$, $i=1,2\ldots k,$ and recall from Remark
\ref{letarhotxiequalbis} that $L_{i}^{\ast}\rho(X)\iota_{H}L_{i}^{\ast}%
=L_{i}^{\ast}\rho(X)$. Thus we conclude that $B_{i}L_{i}^{\ast}\rho
(X)=C_{i}L_{i}^{\ast}$, $i=1,2\ldots k$. Taking adjoints allows us to write
$\rho(X)^{\ast}L_{i}B_{i}^{\ast}=L_{i}C_{i}^{\ast}$, $i=1,2\ldots k$. Since
the images of the representations $\rho$ and $\tau$ commute, and since the
restriction of $\tau$ to $\varphi_{\infty}^{E^{\sigma}}(\sigma(M)^{\prime})$
is a $\ast$-representation, we conclude that for all $a\in\sigma(M)^{\prime}$,
and all $i$, $i=1,2\ldots k$,%
\[
\rho(X)^{\ast}\tau(\varphi_{\infty}^{E^{\sigma}}(a))L_{i}B_{i}^{\ast}%
=\tau(\varphi_{\infty}^{E^{\sigma}}(a))L_{i}C_{i}^{\ast}\text{.}%
\]
So if $a_{i,j}$ are elements of $\sigma(M)^{\prime}$, where $i=1,2\ldots k$,
and the $j$ run over some finite set, and if $h_{i,j}$ are similarly indexed
elements of $H$, then
\begin{align*}
\left\|  \sum\tau(\varphi_{\infty}^{E^{\sigma}}(a_{i,j}))L_{i}C_{i}^{\ast
}h_{i,j}\right\|  ^{2}=\left\|  \sum\rho(X)^{\ast}\tau(\varphi_{\infty
}^{E^{\sigma}}(a_{i,j}))L_{i}B_{i}^{\ast}h_{i,j}\right\|  ^{2}\\
\leq\left\|  \sum\tau(\varphi_{\infty}^{E^{\sigma}}(a_{i,j}))L_{i}B_{i}^{\ast
}h_{i,j}\right\|  ^{2}\text{.}%
\end{align*}
That is, inequality (\ref{i}) is satisfied. However, as we observed in the
course of developing inequality (\ref{i}), its validity for all choices of
$a_{i,j}\in\sigma(M)^{\prime}$ and $h_{i,j}\in H$ is equivalent to the
complete positivity of the map represented by the matrix (\ref{positequat}).
Thus the map is completely positive and the proof is complete. \textbf{ }\ 
\end{proof}

Specializing Theorem \ref{Theorem3.1}, we obtain the following corollary,
which illuminates the ``value distribution theory" for functions in
$H^{\infty}(E)$ at points in the disc $\mathbb{D} ((E^{\sigma})^{\ast})$.

\begin{corollary}
\label{Corollary0.2}Given $\eta^{\ast}$ in $\mathbb{D}((E^{\sigma})^{\ast})$,
then an operator $C\in B(H)$ lies in the image $\sigma\times S(\eta
)(H^{\infty}(E))$ of the representation $\sigma\times S(\eta)$, if and only if
$I\otimes C$ is a bounded operator on the Stinespring dilation space
$\sigma(M)^{\prime}\otimes_{t_{\eta,\eta}}H$ for the completely positive map
${t_{\eta,\eta}}$.
\end{corollary}

\begin{proof}
The operator $I\otimes C$ is bounded with norm less or equal to $1$ if and
only if for every $n$-tuple of elements $a_{1},\ldots,a_{n}$ in $\sigma
(M)^{\prime}$ and every $n$-tuple $h_{1},\ldots,h_{n}$ of elements in $H$,
\[
\sum_{i,j}\langle a_{i}\otimes Ch_{i},a_{j}\otimes Ch_{j}\rangle\leq\sum
_{i,j}\langle a_{i}\otimes h_{i},a_{j}\otimes h_{j}\rangle.
\]
Since $\langle a_{i}\otimes Ch_{i},a_{j}\otimes Ch_{j}\rangle=\langle
h_{i},C^{\ast}t_{\eta,\eta}(a_{i}^{\ast}a_{j})Ch_{j}\rangle$, the inequality
is equivalent to the complete positivity of the map $(id-Ad(C,C))\circ
t_{\eta,\eta}$. Theorem \ref{Theorem3.1}, with $k=1$, gives the result.
\end{proof}

Another consequence of Theorem \ref{Theorem3.1}, which contributes to the
``value distribution theory'' of elements in $H^{\infty}(E)$ is the following
generalization of Schwartz's lemma.

\begin{theorem}
\label{Theorem3.5}Suppose an element $X$ of $H^{\infty}(E)$ has norm at most
one and satisfies the equation $X(0)=0$. Then for every $\eta^{\ast}%
\in\mathbb{D}((E^{\sigma})^{\ast})$ the following assertions are valid:

\begin{enumerate}
\item  If $a$ is a nonnegative element in $\sigma(M)^{\prime}$, and if
$\eta^{\ast}(I_{E}\otimes a)\eta=\langle\eta,a\cdot\eta\rangle\leq a$, then
\[
X(\eta^{\ast})aX(\eta^{\ast})^{\ast}\leq\langle\eta,a\cdot\eta\rangle.
\]

\item  If $\eta^{\otimes k}$ denotes the element $\eta\otimes\eta\otimes
\cdots\otimes\eta\in E^{\otimes k}$, then
\[
X(\eta^{\ast})\langle\eta^{\otimes k},\eta^{\otimes k}\rangle X(\eta^{\ast
})^{\ast}\leq\langle\eta^{\otimes k+1},\eta^{\otimes k+1}\rangle\text{.}%
\]

\item $X(\eta^{\ast})X(\eta^{\ast})^{\ast}\leq\langle\eta,\eta\rangle$.
\end{enumerate}
\end{theorem}

Beside Theorem \ref{Theorem3.1}, the proof requires two initial observations.
The first is certainly known, but we could not find the exact statement in the
literature. Therefore, for completeness, we supply a proof. The second is a
consequence of the first that is the key to our proof of Theorem
\ref{Theorem3.5}.

\begin{lemma}
\label{Lemma3.3}Let $A$ be a unital $C^{\ast}$-algebra. Then an Hermitian
matrix
\[
\left(
\begin{array}
[c]{cc}%
a & b\\
b^{\ast} & c
\end{array}
\right)  \in M_{2}(A)
\]
is nonnegative if and only if $a$ and $c$ are nonnegative in $A$ and for every
$\epsilon>0$,%
\[
b^{\ast}(a+\epsilon)^{-1}b\leq c+\epsilon\text{.}%
\]
\end{lemma}

Here, of course, $\epsilon$ stands for $\epsilon$ times the identity
$\mathbf{1}_{A}$ in $A$.

\begin{proof}
The proof we give is inspired by an argument on page 100 of \cite{vP86}. The
conditions $a\geq0$ and $c\geq0$ are clearly necessary for the matrix to be
nonnegative. Therefore, we shall assume this throughout our proof. Of course,
then $a+\epsilon$ and $c+\epsilon$ are both invertible and we may write%
\begin{multline*}
\left(
\begin{array}
[c]{cc}%
\mathbf{1}_{A} & (a+\epsilon)^{-1/2}b(c+\epsilon)^{-1/2}\\
(c+\epsilon)^{-1/2}b^{\ast}(a+\epsilon)^{-1/2} & \mathbf{1}_{A}%
\end{array}
\right) \\
=\left(
\begin{array}
[c]{cc}%
(a+\epsilon)^{-1/2} & 0\\
0 & (c+\epsilon)^{-1/2}%
\end{array}
\right)  \left(
\begin{array}
[c]{cc}%
a+\epsilon &  b\\
b^{\ast} & c+\epsilon
\end{array}
\right)  \left(
\begin{array}
[c]{cc}%
(a+\epsilon)^{-1/2} & 0\\
0 & (c+\epsilon)^{-1/2}%
\end{array}
\right)
\end{multline*}
Since the matrix $\left(
\begin{array}
[c]{cc}%
a & b\\
b^{\ast} & c
\end{array}
\right)  $ is nonnegative if and only if $\left(
\begin{array}
[c]{cc}%
a+\epsilon &  b\\
b^{\ast} & c+\epsilon
\end{array}
\right)  $ is nonnegative for every positive $\epsilon$, we conclude from this
equation that $\left(
\begin{array}
[c]{cc}%
a & b\\
b^{\ast} & c
\end{array}
\right)  $ is nonnegative if and only if the matrix%
\[
\left(
\begin{array}
[c]{cc}%
\mathbf{1}_{A} & (a+\epsilon)^{-1/2}b(c+\epsilon)^{-1/2}\\
(c+\epsilon)^{-1/2}b^{\ast}(a+\epsilon)^{-1/2} & \mathbf{1}_{A}%
\end{array}
\right)
\]
is nonnegative. However, by Lemma 3.1 in \cite{vP86}, this matrix is
nonnegative if and only if $\left\|  (c+\epsilon)^{-1/2}b^{\ast}%
(a+\epsilon)^{-1/2}\right\|  \leq1$. Since this happens if and only if
$(c+\epsilon)^{-1/2}b^{\ast}(a+\epsilon)^{-1}b(c+\epsilon)^{-1/2}%
\leq\mathbf{1}_{A}$, we conclude that our original matrix is nonnegative if
and only if $b^{\ast}(a+\epsilon)^{-1}b\leq c+\epsilon$.
\end{proof}

\begin{corollary}
\label{Corollary3.4}Suppose that $A$ is a unital subalgebra of a $C^{\ast}%
$-algebra $B$ with the unit of $A$ serving as the unit of $B$. Write $j$ for
the inclusion map of $A$ into $B$. Let $\Psi$ be a unital linear map from $A$
into $B$ and write $\tau:M_{2}(A)\rightarrow M_{2}(B)$ for the matrix of maps%
\[
\left(
\begin{array}
[c]{cc}%
j & j\\
j & \Psi
\end{array}
\right)  \text{.}%
\]
Then $\tau$ is completely positive if and only if $\Psi-j$ is completely positive.
\end{corollary}

\begin{proof}
Suppose first that $\tau$ is positive. Lemma \ref{Lemma3.3} implies that the
matrix $\left(
\begin{array}
[c]{cc}%
c^{-1} & \mathbf{1}_{A}\\
\mathbf{1}_{A} & c
\end{array}
\right)  $ is nonnegative for all nonnegative invertible elements $c\in A$.
Applying $\tau$ we conclude that $\left(
\begin{array}
[c]{cc}%
c^{-1} & \mathbf{1}_{A}\\
\mathbf{1}_{A} & \Psi(c)
\end{array}
\right)  $ is nonnegative for all nonnegative invertible elements $c\in A$.
Hence, by that lemma, we conclude that if $c\geq0$ in $A$ and if $\epsilon>0$,
then $(c^{-1}+\epsilon)^{-1}\leq(\Psi(c)+\epsilon)$. Letting $\epsilon$ tend
to zero gives $c\leq\Psi(c)$ for all $c\in A$, $c\geq0$. That is, $\Psi-j$ is
a positive linear map. Next, suppose the $\Psi-j$ is a positive linear map,
and let $\left(
\begin{array}
[c]{cc}%
a & b\\
b^{\ast} & c
\end{array}
\right)  $ be a nonnegative element in $M_{2}(A)$. Then again by Lemma
\ref{Lemma3.3}, we have $b^{\ast}(a+\epsilon)^{-1}b\leq c+\epsilon$ for all
$\epsilon>0$. Since $\Psi-j$ is assumed to be positive, we conclude that
$b^{\ast}(a+\epsilon)^{-1}b\leq\Psi(c)+\epsilon$ for all $\epsilon>0$. Thus
once more by Lemma \ref{Lemma3.3} we conclude that $\tau\left(
\begin{array}
[c]{cc}%
a & b\\
b^{\ast} & c
\end{array}
\right)  \geq0$ in $M_{2}(A)$. Thus, we conclude that $\tau$ is positive if
and only if $\Psi-j$ is positive. To conclude that $\tau$ is completely
positive if and only if $\Psi-j$ is completely positive simply note that the
inflated map $\tau_{n}$ (acting on $M_{n}(M_{2}(A))$) is, after reshuffling,
just $\left(
\begin{array}
[c]{cc}%
j^{(n)} & j^{(n)}\\
j^{(n)} & \Psi_{n}%
\end{array}
\right)  $ where $j^{(n)}$ is the identity map embedding $M_{n}(A)$ in
$M_{n}(B)$.
\end{proof}

\begin{proof}
(of Theorem \ref{Theorem3.5}) Since condition (2) is obtained from (1) by
taking $a=\langle\eta^{\otimes k},\eta^{\otimes k}\rangle$ and condition (3)
is (1) with $a=I$, we see that it suffices to prove (1).

To this end, set $X(\eta^{\ast})=C$ and apply the sufficiency assertion in
Theorem \ref{Theorem3.1} to conclude that if we take $k=2$, $B_{i}=I$,
$i=1,2$, $\eta_{1}=0$, $C_{1}=0$, $\eta_{2}=\eta$ and $C_{2}=C$. Then the
resulting matrix (\ref{positequat}) is $\left(
\begin{array}
[c]{cc}%
j & j\\
j & \Psi
\end{array}
\right)  $ where $j$ is the inclusion of $\sigma(M)^{\prime}$ in $B(H)$ and
$\Psi=(id-Ad(C,C))\circ(id-\theta_{\eta,\eta})^{-1}$. It follows from
Corollary \ref{Corollary3.4} that the matrix is completely positive if and
only if $\Psi-j$ is completely positive. We conclude that if $a\in
\sigma(M)^{\prime}$ is nonnegative and satisfies $\theta_{\eta,\eta}(a)\leq
a$, then $(\Psi-j)(a-\theta_{\eta,\eta}(a))\geq0$. If we write $a-\theta
_{\eta,\eta}(a)=a-\langle\eta,a\cdot\eta\rangle$, and expand $(\Psi
-j)(a-\langle\eta,a\cdot\eta\rangle)$, we obtain the inequality,%
\begin{align*}
0\leq(\Psi-j)(a-\langle\eta,a\cdot\eta\rangle)\\
=(id-Ad(C,C))(a)-(a-\langle\eta,a\cdot\eta\rangle)\\
=\langle\eta,a\cdot\eta\rangle-CaC^{\ast}\\
=\langle\eta,a\cdot\eta\rangle-X(\eta^{\ast})aX(\eta^{\ast})^{\ast}\text{,}%
\end{align*}
which completes the proof.
\end{proof}

\section{Examples of the Nevanlinna-Pick Theorem}

As we observed in Remark \ref{Remark3.2}, Theorem \ref{Theorem3.1} contains
the classical Nevanlinna-Pick theorem. In this section, we show how a number
of other generalizations of the Nevanlinna-Pick theorem that appear in the
literature are also consequences of Theorem \ref{Theorem3.1}.

We begin with the setting of Example \ref{Example2.1}, in which
$M=E=\mathbb{C}$ and we let $H$ be an arbitrary Hilbert space. Then
$\sigma:M\rightarrow B(H)$ is scalar multiplication, $\sigma(M)^{\prime}%
=B(H)$, and $E^{\sigma}=\{\eta:H\rightarrow\mathbb{C}\otimes H\mid\eta
a=(a\otimes I_{H})\eta$, $a\in M=\mathbb{C\}}$. Since $\mathbb{C}\otimes H$ is
naturally identified with $H$, we may identify $E^{\sigma}$ with $B(H)$ viewed
as a $W^{\ast}$-correspondences over $B(H)$ in the usual way; i.e.,
$E^{\sigma}$ is the identity correspondence $_{B(H)}B(H)_{B(H)}$. The algebra
$H^{\infty}(E)$ is just $H^{\infty}(\mathbb{T})$ and, given $\eta\in
E^{\sigma}$ $\left\|  \eta\right\|  <1$, the representation of $H^{\infty
}(E)=H^{\infty}(\mathbb{T})$ $\eta$ induces is the $H^{\infty}$-functional
calculus: $f\rightarrow f(\eta^{\ast})$, $f\in H^{\infty}(\mathbb{T})$.
Theorem \ref{Theorem3.1} then gives the following result.

\begin{theorem}
\label{Theorem4.2.1}Given operators $T_{1}$, $T_{2},\ldots,T_{k}\in B(H)$,
with $\left\|  T_{i}\right\|  <1$, $1\leq i\leq k$, and operators $B_{1}$,
$B_{2},\ldots,$ $B_{k}$, $C_{1}$, $C_{2},\ldots,C_{k}$ in $B(H)$, then there
is a function $f\in H^{\infty}(\mathbb{T})$, with $\left\|  f\right\|  \leq1$
such that $B_{i}f(T_{i})=C_{i}$, $1\leq i\leq k$, if and only if the map on
$M_{k}(B(H))$ defined by the matrix of maps
\[
\left(  \left(  Ad(B_{i},B_{j})-Ad(C_{i},C_{j})\right)  \circ\left(
id-\theta_{T_{i}^{\ast},T_{j}^{\ast}}\right)  ^{-1}\right)
\]
is completely positive.
\end{theorem}

Note that the map $\left(  id-\theta_{T_{i}^{\ast},T_{j}^{\ast}}\right)
^{-1}$ is given by the formula
\[
\left(  id-\theta_{T_{i}^{\ast},T_{j}^{\ast}}\right)  ^{-1}(A)=\sum
_{k=0}^{\infty}T_{i}^{k}AT_{j}^{\ast k}.
\]

Now we interchange the roles of $E$ and $E^{\sigma}$, so that $M=B(H)$,
$E=\,_{B(H)}B(H)_{B(H)}$ is the identity correspondence over $B(H)$,
$M^{\prime}=\mathbb{C}I_{H}$, and $E^{\sigma}=\mathbb{C}$. Then $H^{\infty
}(E)$ can be identified with $H^{\infty}(\mathbb{T})\otimes B(H)$, which in
turn may be viewed as the space of all bounded analytic $B(H)$-valued
functions on $\mathbb{D}$. The associated Fock space $\mathcal{F}(E)$ is just
the direct sum of copies of $_{B(H)}B(H)_{B(H)}$ and $\mathcal{F}(E)\otimes H$
can be identified in a natural way with $H^{2}(\mathbb{T})\otimes H$ on which
$H^{\infty}(\mathbb{T})\otimes B(H)$ acts as usual. Note that $\varphi
_{\infty}(M)$, which is a subalgebra of $H^{\infty}(E)$ is identified with
$\mathbb{C}\otimes B(H)$ - the constant $B(H)$-valued functions in $H^{\infty
}(\mathbb{T})\otimes B(H)$. Since $E^{\sigma}=\mathbb{C}$, $\mathbb{D}%
((E^{\sigma})^{\ast})$ is simply the open unit disc $\mathbb{D}$ in the
complex plane. Further, the representation of $H^{\infty}(\mathbb{T})\otimes
B(H)$ determined by an $\eta\in\mathbb{D}$ is given by the formula
$F\rightarrow F(\bar{\eta})$, where $F\in H^{\infty}(\mathbb{T})\otimes B(H)$
is viewed as a bounded analytic $B(H)$-valued function on $\mathbb{D}$. With
these observations, we see that the following theorem is an immediate
consequence of Theorem \ref{Theorem3.1}.

\begin{theorem}
\label{Theorem4.3.1}Given $z_{1},z_{2},\ldots,z_{k}$ in $\mathbb{D}$ and
operators $B_{1}$, $B_{2}$, $\ldots$, $B_{k}$, $C_{1}$, $C_{2}$, $\ldots$,
$C_{k}$ in $B(H)$, there is an element $G\in H^{\infty}(\mathbb{T})\otimes
B(H)$, with $\left\|  G\right\|  _{\infty}\leq1$, such that $B_{i}%
G(z_{i})=C_{i}$ if and only if the operator matrix
\[
\left(  \frac{B_{i}B_{j}^{\ast}-C_{i}C_{j}^{\ast}}{1-z_{i}\bar{z}_{j}%
}\right)
\]
is a positive element of $M_{k}(B(H))$.

Moreover, if all the $B_{i}$ and $C_{i}$ lie in a given von Neumann algebra
$N\subseteq B(H)$, and if the matrix $\left(  \frac{B_{i}B_{j}^{\ast}%
-C_{i}C_{j}^{\ast}}{1-z_{i}\bar{z}_{j}}\right)  $ is a positive element in
$M_{k}(N)$, then $G$ can be chosen in $H^{\infty}(\mathbb{T})\otimes N$.
\end{theorem}

\begin{proof}
The first part of the theorem follows immediately from Theorem
\ref{Theorem3.1} and Remark \ref{Remark3.2}. For the second, simply observe
that if the $B_{i}$ and $C_{i}$ all lie in $N\subseteq B(H)$, then by Theorem
\ref{Theorem3.1}, $G$ commutes with $\varphi_{\infty}(N^{\prime})$. As we
noted above, $\varphi_{\infty}(B(H))$ is identified with $\mathbb{C}\otimes
B(H)\subseteq H^{\infty}(\mathbb{T})\otimes B(H)$ and so $\varphi_{\infty
}(N^{\prime})$ is identified with $\mathbb{C}\otimes N^{\prime}$.
Consequently, $G\in H^{\infty}(\mathbb{T})\otimes N^{\prime\prime}=H^{\infty
}(\mathbb{T})\otimes N$.
\end{proof}

Of course this result contains the matricial Nevanlinna-Pick interpolation
theorem first proved by Sz.-Nagy and Koranyi in \cite{szNK56}. For an
elementary presentation of their result, which includes helpful references to
the literature, we recommend the note by Power \cite{sP89}.

Next consider the setting of Example \ref{Example2.2} in which $M=\mathbb{C}$
and $E=\mathbb{C}^{n}$. We fix a Hilbert space $H$ and the (one and only)
representation $\sigma$ of $M$ on $H$. Then $\sigma(M)^{\prime}=B(H)$ and
$E^{\sigma}=\{\eta:H\rightarrow\mathbb{C}^{n}\otimes H\mid\eta a=(a\otimes
I_{H})\eta$, $a\in M\}=B(H,\mathbb{C}^{n}\otimes H)$. The right action of
$\sigma(M)^{\prime}=B(H)$ is by right multiplication and the inner product is
given by the formula $\langle\eta_{1},\eta_{2}\rangle=\eta_{1}^{\ast}\eta_{2}%
$. Thus, $E^{\sigma}$ is simply column space over $B(H)$, $C_{n}(B(H))$ - that
is, the $n\times1$ matrices over $B(H)$. The algebra $H^{\infty}(E)$ is the
ultra weak closure of the noncommutative disc algebra of Popescu \cite{gP96a}
and was denoted by $\mathcal{L}_{n}$ by Davidson and Pitts \cite{DP98a}. Every
column matrix of operators $T:=(T_{i})_{i=1}^{n}\in C_{n}(B(H))$ with
$\left\|  T\right\|  <1$ (i.e., with $\sum_{i=1}^{n}T_{i}^{\ast}T_{i}%
\leq\varepsilon I_{H}$, $\varepsilon<1)$ defines a representation of
$E=\mathbb{C}^{n}$ via the formula $e_{i}\rightarrow T_{i}^{\ast}$, where the
$e_{i}$ denote the unit basis vectors in $\mathbb{C}^{n}$. This extends to a
representation $\pi$ of $\mathcal{L}_{n}$ which, on basis vectors
$e_{\text{\b{i}}}$ for words \thinspace$\underline{i}=(i_{1},i_{2},\ldots
i_{m})$ in the free semigroup on $n$ generators, is given by the formula
$\pi(e_{\underline{i}})=T_{i_{1}}^{\ast}T_{i_{2}}^{\ast}\cdots T_{i_{m}}%
^{\ast}$. We conclude from Theorem \ref{Theorem3.1} that the following theorem
is valid.

\begin{theorem}
\label{Theorem4.4.1}Let $T^{(1)}$, $T^{(2)}$, $\ldots$, $T^{(k)}\in E^{\sigma
}$ with $\left\|  T^{(i)}\right\|  <1$ be given and let the corresponding
representations of $\mathcal{L}_{n}$ be denoted by $\pi_{i}$, $i=1,2,\ldots
k$. Also let operators $B_{1}$, $B_{2}$, $\ldots$, $B_{k}$, $C_{1}$, $C_{2}$,
$\ldots$, $C_{k}$ in $B(H)$ be given. Then there is an $X\in\mathcal{L}_{n}$
with $\left\|  X\right\|  \leq1$ such that $B_{i}\pi_{i}(X)=C_{i}$, $1\leq
i\leq k$, if and only if for every $l<\infty$ and every choice of operators
$a_{(i,j)}\in B(H)$, $1\leq i\leq k$, $1\leq j\leq l$, the matrix%
\[
\left(  B_{i}t_{i,m}(a_{(i,j)}^{\ast}a_{(m,p)})B_{j}^{\ast}-C_{i}%
t_{i,m}(a_{(i,j)}^{\ast}a_{(m,p)})C_{m}^{\ast}\right)
\]
is positive, where%
\[
t_{i,m}(a)=\sum T_{\underline{p}}^{(i)\ast}aT_{\underline{p}}^{(m)}\text{,}%
\]
\underline{$p$}$=(p_{1},p_{2},\ldots,p_{j})$, $T_{\underline{p}}%
^{(i)}=T_{p_{1}}^{(i)}T_{p_{2}}^{(i)}\cdots T_{p_{j}}^{(i)}$, and the sum runs
over all $j$ and all \underline{$p$}$.$
\end{theorem}

Consider next what happens when we interchange the roles of $E$ and
$E^{\sigma}$ in the previous setting. For simplicity of notation, we will take
$M=B(H)$, i.e., we will let $\sigma$ be the identity representation of $B(H)$
on $H$, $E$ will be column space over $B(H)$, $C_{n}(B(H))$, viewed as a
correspondence over $B(H)$ in the usual fashion. Then $M^{\prime}%
=\mathbb{C}I_{H}$, and $E^{\sigma}=\mathbb{C}^{n}$. The algebra $H^{\infty
}(E)$ is $B(H)\otimes\mathcal{L}_{n}$ and, given $\eta\in\mathbb{C}^{n}$, with
$\left\|  \eta\right\|  <1$, one obtains a representation of $E$ on $H$ by
mapping $(T_{i}):=T$ to $\sum\bar{\eta}_{i}T_{i}$, where $\eta=(\eta_{1}%
,\eta_{2},\ldots,\eta_{n})$. The integrated form $\pi_{\eta}$ of this
representation of $E$ is the representation of $B(H)\otimes\mathcal{L}_{n}$
that maps $S\otimes X\in B(H)\otimes\mathcal{L}_{n}$ to $X(\bar{\eta})S\in
B(H)$, i.e., $\pi_{\eta}(S\otimes X)=X(\bar{\eta})S$. Here $X\mapsto
X(\bar{\eta})$ is the one-dimensional representation of $\mathcal{L}_{n}$
determined by $\eta$. Using Theorem \ref{Theorem3.1} and Remark
\ref{Remark3.2}, we obtain the following theorem that extends Theorem 4.1 of
\cite{gP98}, Corollary 2.8 of \cite{gP01}, Theorem 2.4 of \cite{AP00} and
Theorem 3.2 of \cite{DP98b}.

\begin{theorem}
\label{Theorem4.5.1}Given $\eta_{1}$, $\eta_{2}$, $\ldots$, $\eta_{k}$ in the
open unit ball of $\mathbb{C}^{n}$, and given operators $C_{1}$, $C_{2}$,
$\ldots$, $C_{k}\in B(H)$, there is a $Y\in B(H)\otimes\mathcal{L}_{n}$ with
$\left\|  Y\right\|  \leq1$ such that $\pi_{\eta_{i}}(Y)=C_{i}$, $1\leq i\leq
k$, if and only if the $k\times k$ operator matrix
\[
\left(  \frac{I-C_{i}C_{j}^{\ast}}{1-\langle\eta_{i},\eta_{j}\rangle}\right)
\]
is positive in $M_{k}(B(H))$.

Moreover, if each $C_{i}$ lies in the von Neumann algebra $N\subseteq B(H)$,
then we can choose $Y\in N\otimes\mathcal{L}_{n}$.
\end{theorem}

Consider the setting of quiver algebras discussed in Example \ref{Example2.3}
and let $E(C)$ be the correspondence associated to a matrix $C$ with
non-negative integer entries. We will use the notation of that example. Our
von Neumann algebra $M$ is the diagonal matrices $D_{n}$ and the
representation $\sigma$ will be determined by the multiplicity vector
$\mathbf{m}=(m_{1},m_{2},\ldots,m_{n})$, where all the $m_{i}$ equal $1$.
Then, of course, the Hilbert space $H$ of $\sigma$ is just $\mathbb{C}^{n}$
and $\sigma(M)^{\prime}=D_{n}$, also. The condition of Theorem
\ref{Theorem3.1} (and Remark \ref{Remark3.2}) refers to the complete
positivity of a map from $M_{k}(D_{n})$ into $M_{k}(B(H))$. To apply this
condition, we need an easy consequence of Choi's theorem \cite{mC75}.

\begin{lemma}
\label{Lemma4.6.1}Suppose $\varphi_{i,j}$, $1\leq i,j\leq k$, are maps from
$D_{n}$ to $B(H)$ and write $\varphi=\left(  \varphi_{ij}\right)  $ for the
map from $M_{k}(D_{n})$ to $M_{k}(B(H))$ that this family defines. For each
$m$, $1\leq m\leq n$, let $A_{m}=(\varphi_{ij}(e_{mm}))\in M_{k}(B(H))$. Then
$\varphi$ is completely positive if and only if $A_{m}\geq0$ for all $m$.
\end{lemma}

\begin{proof}
Suppose $\varphi$ is completely positive and let $B_{m}=(b_{ij})_{ij=1}^{k}\in
M_{k}(D_{n})$ be the matrix with $b_{ij}=e_{mm}$. Then clearly $B_{m}$ is a
positive element in $M_{k}(D_{n})$ and, evidently, $\varphi(B_{m})=A_{m}$.
Since $\varphi$ is assumed to be completely positive (in fact `positive'
suffices), we see that each $A_{m}$ is nonnegative.

For the converse, suppose the $A_{m}$ are all nonnegative. For each $i$ and
$j$, $1\leq i,j\leq k$, and for each $m$, $1\leq m\leq n$, define
$\varphi_{ij}^{(m)}:\mathbb{C}\rightarrow B(H)$ by the formula $\varphi
_{ij}^{(m)}(\lambda)=\varphi_{ij}(\lambda e_{mm})$. Let $\varphi^{(m)}=\left(
\varphi_{ij}^{(m)}\right)  $ be the induced map from $M_{k}(\mathbb{C})$ to
$M_{k}(B(H))$. Since $M_{k}(D_{n})$ is isomorphic to the direct sum
$M_{k}(\mathbb{C})\oplus M_{k}(\mathbb{C})\oplus\cdots\oplus M_{k}%
(\mathbb{C})$, we see easily that $\varphi=\varphi^{(1)}\oplus\varphi
^{(2)}\oplus\cdots\oplus\varphi^{(n)}$. Further, it is evident that $\varphi$
is completely positive if and only if each $\varphi^{(m)}$ is completely
positive. However, since each $A_{m}\geq0$, Choi's Theorem \cite[Theorem
2]{mC75} implies that $\varphi^{(m)}$ is completely positive.
\end{proof}

Recall that we are assuming that the multiplicity vector of our representation
$\sigma$ of $D_{n}$ is $(1,1,\ldots,1)$ and so $E^{\sigma}=E(C^{t})$, by
Example \ref{Example2.3}. If $\eta\in E(C^{t})$ with $\left\|  \eta\right\|
<1$, we can write $\eta=\left(  \eta_{ij}\right)  $, $\eta_{ij}\in H_{ji}$ and
$\sum\left\|  \eta_{ij}\right\|  ^{2}<1$.

The following theorem now follows from the discussion above and from
Theorem~\ref{Theorem3.1}.

\begin{theorem}
\label{NPquiver} Let $C$ and $E(C)$ be as above and let $\sigma$ be the
identity representation of $D_{n}$ on $\mathbb{C}^{n}$. Given $\eta_{1}^{\ast
},\eta_{2}^{\ast},\ldots,\eta_{k}^{\ast}\in\mathbb{D}(E(C^{t})^{\ast})$ and
$n\times n$ matrices $C_{1},C_{2},\ldots,C_{k}$, there is an element $X\in
H^{\infty}(E(C))$ such that $\Vert{X}\Vert\leq1$ and such that $X(\eta
_{i}^{\ast})=C_{i}$ for all $1\leq i\leq k$ if and only if the $n$ $k\times k$
matrices $A_{1},A_{2},\ldots A_{n}$, defined by
\[
A_{m}=(t_{\eta_{i},\eta_{j}}(e_{mm})-C_{i}t_{\eta_{i},\eta_{j}}(e_{mm}%
)C_{j}^{\ast})_{i,j}%
\]
are nonnegative matrices.
\end{theorem}

\begin{example}
\label{Example4.7}Suppose $M$ is a von Neumann subalgebra of $B(H)$ and that
$\alpha$ is a normal, unital endomorphism of $M$. Following Example
\ref{Example2.4}, we let $E$ be $_{\alpha}M$. Then, after taking our
representation $\sigma$ of $M$ to be the identity representation, we see that
$E^{\sigma}=\{\eta\in B(H)\mid\eta a=\alpha(a)\eta$, $a\in M\}$. Given
$\eta\in E^{\sigma}$ with $\left\|  \eta\right\|  <1$, the map $S(\eta)$ in
the covariant representation $(S(\eta),\sigma)$ of $E=\,_{\alpha}M$ that it
defines is given by the formula $S(\eta)(\xi)=\eta^{\ast}\xi$, $\xi
\in\,_{\alpha}M$. Every $X\in H^{\infty}(E)$, has a ``Fourier series'',
$X=\sum X_{n}$, which is Cesaro summable to $X$, and each coefficient $X_{n}$
belongs to $E^{\otimes n}$. Since $E^{\otimes n}\simeq\,_{\alpha^{n}}M$ and
since $_{\alpha^{n}}M=M$ as sets, we may think of the $X_{n}$'s as elements of
$M$. Then the integrated form, $\sigma\times S(\eta)$, of $(S(\eta),\sigma)$
applied to $X$ yields the equation%
\[
X(\eta^{\ast})=\sigma\times S(\eta)=\sum(\eta^{\ast})^{n}X_{n}\text{.}%
\]
So, given operators $C_{1}$, $C_{2}$, $\ldots$, $C_{k}$ in $B(H)$ and
$\eta_{1},\eta_{2}$, $\ldots$, $\eta_{k}\in E^{\sigma}$, with $\left\|
\eta_{i}\right\|  <1$, Theorem \ref{Theorem3.1} gives a necessary and
sufficient condition for the existence of an $X\in H^{\infty}(E)$ satisfying
$\left\|  X\right\|  \leq1$ such that%
\[
\sum_{n=0}^{\infty}(\eta_{i}^{\ast})^{n}X_{n}=C_{i}\text{,}%
\]
$1\leq i\leq k$.
\end{example}

Our final example concerns Nevanlinna-Pick interpolation in nest algebras.
Recall that a \emph{nest }of subspaces in a Hilbert space $H$ is a family
$\mathcal{N}$ of subspaces of $H$ that is totally ordered by inclusion. We
identify a subspace with the orthogonal projection onto it and so think of
$\mathcal{N}$ as a family of projections on $H$. Also, we assume that
$\mathcal{N}$ contains $0$ and $I$ and is closed in the ultraweak topology on
$B(H)$. We write $Alg\,\mathcal{N}$ for $\{T\in B(H)\mid(I-N)TN=0$,
$N\in\mathcal{N}\}$. Then $Alg\,\mathcal{N}$ is an ultraweakly closed algebra
called the \emph{nest algebra }determined by $\mathcal{N}$. We follow
\cite{kD88} for matters relating to nest algebras.

In general, it does not seem possible to write $Alg\,\mathcal{N}$ as
$H^{\infty}(E)$ for some correspondence. However, if $\mathcal{N}$ is finite
then it is, as we shall show. This allows us to apply our theory to finite
nest algebras. We can then bring our analysis to bear on a general nest
algebra by observing that $Alg\,\mathcal{N}=\cap\{Alg\,\mathcal{M}%
\mid\mathcal{M}\subseteq\mathcal{N}$, $\mathcal{M}$ finite$\}$.

Observe that if $\mathcal{N}=\{0=P_{0}<P_{1}<P_{2}<\cdots<P_{n}=I\}$ is a
finite nest of orthogonal projections in $B(H)$, then $\mathcal{D}%
:=\mathcal{N}^{\prime}$ is a von Neumann subalgebra of $B(H)$. Write
$Q_{k}=P_{k}-P_{k-1}$, $k=1,2,\ldots,n$ and let $E:=\{T\in B(H)\mid
Q_{k}T=TQ_{k+1}=Q_{k}TQ_{k+1}$, $1\leq k\leq n-1\}$. Then if we think of
$B(H)$ as being written as block matrices associated with the decomposition of
$H$ as $H=\sum^{\oplus}Q_{k}H$, $E$ may be viewed as the collection of those
matrices that are supported on diagonal immediately above the main diagonal.
It is clear that $E$ is an ultraweakly closed bimodule over $\mathcal{D}$
under operator multiplication and that $E$ has the $\mathcal{D}$-valued inner
product making $E$ into a $W^{\ast}$-correspondence over $\mathcal{D}$ given
by the formula: $\langle T,S\rangle=\Phi(T^{\ast}S)$, where $\Phi$ is the
conditional expectation of $B(H)$ onto $\mathcal{D}$, $\Phi(T)=\sum_{k}%
Q_{k}TQ_{k}$. It is also evident that we may identify $E\otimes_{\mathcal{D}%
}E$ with the space of products $E^{2}\subseteq B(H)$, which should be viewed
as the second superdiagonal in the block decomposition of $B(H)$. And, of
course, more generally $E^{\otimes k}$ may be identified with the
``$k$-superdiagonal'', $E^{k}\subseteq B(H)$. Of course, $E^{\otimes n}%
=E^{n}=0$. Also, note that the subspaces $E^{k}=E^{\otimes k}$ and
$E^{l}=E^{\otimes l}$ are orthogonal, when $k\neq l$, i.e., if $T\in E^{k}$
and $S\in E^{l}$, then $\Phi(T^{\ast}S)=0$. Hence the Fock space
$\mathcal{F}(E)$ is isomorphic, as a $W^{\ast}$-correspondence, to
$Alg\,\mathcal{N\,}=\mathcal{D}+E+E^{2}+\cdots+E^{n-1}$ viewed as a subspace
of $B(H)$ with the inner product $\langle T,S\rangle=\Phi(T^{\ast}S)$. It is
now evident from the definitions of $\mathcal{T}_{+}(E)$ and $H^{\infty}(E)$
that both of these algebras are completely isometrically isomorphic to
$Alg\,\mathcal{N}$. Indeed, $\mathcal{T}_{+}(E)=H^{\infty}(E)$ viewed as left
multiplication by elements in $Alg\,\mathcal{N}$ acting on $Alg\,\mathcal{N}$.
Our primary inspiration for the following theorem and corollary (Corollary
\ref{Corollary4.8.2}) comes from the paper \cite{KMT93} by Katsoulis, Moore
and Trent. See their Theorems 2 and 3 in particular. The theorem, itself, is
related to Theorem 8 of \cite{BG85}.

\begin{theorem}
\label{Theorem4.8.1}Let $\mathcal{N}$ be a nest of subspaces in a Hilbert
space $H$ and let $B$ and $C$ be operators in $B(H)$. Then there is an
operator $X\in Alg\,\mathcal{N}$ such that $\left\|  X\right\|  \leq1$ and
such that $BX=C$ if and only if
\begin{equation}
CNC^{\ast}\leq BNB^{\ast} \label{InequalCB}%
\end{equation}
for every $N\in\mathcal{N}$.
\end{theorem}

\begin{proof}
First observe that the inequality (\ref{InequalCB}) is clearly necessary. For
if $X\in Alg\,\mathcal{N}$, with $\left\|  X\right\|  \leq1$ and $BX=C$. Then
for $N\in\mathcal{N}$, we have $XNX^{\ast}=NXNX^{\ast}N\leq N$. Consequently,
$CNC^{\ast}\leq BNB^{\ast}$ for all $N\in\mathcal{N}$.

For sufficiency, assume first that the theorem has been proved for all finite
nests. Given a general nest $\mathcal{N}$ and operators $B$ and $C$ in $B(H)$
satisfying $CNC^{\ast}\leq BNB^{\ast}$ for all $N\in\mathcal{N}$, then for
every finite subnest $\mathcal{M}$ contained in $\mathcal{N}$, we may find an
operator $X_{\mathcal{M}}\in Alg\,\mathcal{M}$ such that $\left\|
X_{\mathcal{M}}\right\|  \leq1$ an $BX_{\mathcal{M}}=C$. Let $X$ be a
weak-$\ast$ limit point of the set $\{X_{\mathcal{M}}\mid\mathcal{M}%
\subseteq\mathcal{N}$, $\mathcal{M}$ finite$\}$. It is then straightforward to
check that $BX=C$, $X\in\cap\{Alg\,\mathcal{M}$ $\mid\mathcal{M}%
\subseteq\mathcal{N}$, $\mathcal{M}$ finite$\}=Alg\,\mathcal{N}$ and, of
course, $\left\|  X\right\|  \leq1$.

It therefore suffices to assume that $\mathcal{N}$ is a finite nest. We want
to apply Theorem \ref{Theorem3.1}. For this purpose, we let our representation
$\sigma$ be the identity representation. Also, we have $k$, the number of
operators $B$ and $C$ involved, equal to $1$. We want to find $X\in H^{\infty
}(E)$ ($\simeq Alg\,\mathcal{N}$) and an element $\eta\in E^{\sigma}$ so that

\begin{enumerate}
\item $\left\|  X\right\|  \leq1$,

\item $BX(\eta^{\ast})=C$,

\item $X(\eta^{\ast})\in Alg\,\mathcal{N}$, and

\item $\left\|  X(\eta^{\ast})\right\|  \leq1$.
\end{enumerate}

\noindent Of course condition (4) would be immediate from (1) if $\eta$ could
be chosen with $\left\|  \eta\right\|  <1$, which is the only situation
considered in Theorem \ref{Theorem3.1}. However, we need to rework the
arguments in the proof of Theorem \ref{Theorem3.1} slightly to fit the current
setting to allow an $\eta$ that has norm equal to one. This is possible
because $E^{\otimes n}=0$, and so all putatively infinite series are, in fact,
finite. The $\eta$ we want is defined by the formula%
\[
\eta^{\ast}(T\otimes h)=Th\text{,}%
\]
$T\in E$. We need to check that $\eta\in E^{\sigma}$. But this is easy. Recall
that $\sigma$ is the identity and that the right and left actions of
$\mathcal{D}$ on $E$ are given by multiplication. That is, for $S\in
\mathcal{D}$ and $T\in E$, $T\cdot S=TS$, while $\varphi(S)T=ST$.
Consequently, for $S\in\mathcal{D}$, and $T\otimes h\in E\otimes_{\mathcal{D}%
}H$, $\sigma^{E}\circ\varphi(S)(T\otimes h)=ST\otimes h$. So for
$S\in\mathcal{D}$ and $T\otimes h\in E\otimes_{\mathcal{D}}H$, we have
\begin{multline*}
\eta^{\ast}\sigma^{E}\circ\varphi(S)(T\otimes h)=\eta^{\ast}(ST\otimes
h)=\eta^{\ast}STh\\
=S\eta^{\ast}(T\otimes h)=\sigma(S)\eta^{\ast}(T\otimes h)\text{.}%
\end{multline*}
After taking adjoints, this shows that $\eta\in E^{\sigma}$. Note that as an
operator, $\eta^{\ast}$ is an isometry and so as an element of $E^{\sigma}$,
the norm of $\eta$ is $1$. Nevertheless, observe that for all $X\in H^{\infty
}(E)$, ``$X(\eta^{\ast})=X$''. More precisely, if we identify $\mathcal{F}(E)$
with $Alg\,\mathcal{N}$ as above, so that the action of $H^{\infty}(E)$ on
$\mathcal{F}(E)$ is left multiplication by $Alg\,\mathcal{N}$ on itself, then
for $X\in H^{\infty}(E)$, $X(\eta^{\ast})$ is operator on $H$ obtained by
thinking of $Alg\,\mathcal{N}$ as acting on $H$. To see this, recall how
$X(\eta^{\ast})$ is defined (see the discussion just before the statement of
Theorem \ref{Theorem3.1}): $X(\eta^{\ast})$ is defined to be $\sigma\times
S(\eta)(X)$, where, in our situation $\sigma$ is the identity representation
of $\mathcal{D}$, and $S(\eta)$ is the representation of $E$ on $H$ such that
$\widetilde{S(\eta)}=\eta^{\ast}$. That is, for all $X\in E$ and $h\in H$,
$X(\eta^{\ast})h=S(\eta)(X)h=\widetilde{S(\eta)}(X\otimes h)=\eta^{\ast
}(X\otimes h)=Xh$. Since $\sigma$ is the identity representation of
$\mathcal{D}$ on $H$ and $\sigma\times S(\eta)$ is an algebra homomorphism,
this calculation shows that $X(\eta^{\ast})h=\sigma\times S(\eta)(X)h=Xh$ for
all $X\in H^{\infty}(E)=Alg\,\mathcal{N}$. Thus, we see that inequality (4)
will be satisfied if we can satisfy inequality (1) and, of course, the
equation (3) will be satisfied if we can find the desired $X\in H^{\infty}%
(E)$. That is, the proof will be complete if we can find an $X\in H^{\infty
}(E)$, $\left\|  X\right\|  \leq1$, so that $CX(\eta^{\ast})=B$. To use
Theorem \ref{Theorem3.1} for this purpose, it suffices to see if
$(Ad(B)-Ad(C))\circ(id-\theta_{\eta,\eta})^{-1}$ is completely positive on
$\sigma(\mathcal{D})^{\prime}=\mathcal{D}^{\prime}$. And to do this we need to
show that $\eta$ is given by the formula:%
\begin{equation}
\eta h=\sum_{k=1}^{n-1}\frac{1}{\left\|  Q_{k+1}h\right\|  ^{2}}\left(
Q_{k}h\otimes(Q_{k+1}h)^{\ast}\right)  \otimes h\text{,} \label{alteta}%
\end{equation}
$h\in H$, where $\left(  Q_{k}h\otimes(Q_{k+1}h)^{\ast}\right)  $ denotes the
rank one operator determined by the vectors $Q_{k}h$ and $Q_{k+1}h$. Indeed,
for $h\in H$ and $S\otimes k\in E\otimes_{\mathcal{D}}H$, we have
\begin{multline*}
\langle\sum_{j=1}^{n-1}\frac{1}{\left\|  Q_{j+1}h\right\|  ^{2}}\left(
Q_{j}h\otimes(Q_{j+1}h)^{\ast}\right)  \otimes h,S\otimes k\rangle\\
=\sum_{j=1}^{n-1}\frac{1}{\left\|  Q_{j+1}h\right\|  ^{2}}\langle
h,[(Q_{j+1}h)\otimes(Q_{j}h)^{\ast}]Sk\rangle\\
=\sum_{j=1}^{n-1}\frac{1}{\left\|  Q_{j+1}h\right\|  ^{2}}\langle h,\langle
Sk,Q_{j}h\rangle Q_{j+1}h\rangle\\
=\sum_{j=1}^{n-1}\frac{1}{\left\|  Q_{j+1}h\right\|  ^{2}}\langle
h,Q_{j+1}h\rangle\langle Q_{j}h,Sk\rangle\\
=\sum_{j=1}^{n-1}\langle Q_{j}h,Sk\rangle\text{.}%
\end{multline*}
However, since $Q_{n}E=0$, this last sum really is $\sum_{j=1}^{n}\langle
Q_{j}h,Sk\rangle=\langle h,Sk\rangle=\langle h,\eta^{\ast}(S\otimes
k)\rangle=\langle\eta h,S\otimes k\rangle$. This verifies equation
(\ref{alteta}).

Now $\mathcal{D}^{\prime}$ is $\mathcal{N}^{\prime\prime}$, which is the span
of $\mathcal{N}$. This, in turn, is the span of the $Q_{i}$'s, $i=1,2,\ldots
,n$. By definition $\theta_{\eta,\eta}(T)=\eta^{\ast}(T\cdot\eta)\in B(H)$ for
all $T\in\mathcal{D}^{\prime}$. However, from equation (\ref{alteta})
$(Q_{j}\cdot\eta)h=\eta(Q_{j}h)=\frac{1}{\left\|  Q_{j}h\right\|  ^{2}}\left(
Q_{j-1}h\otimes(Q_{j}h)^{\ast}\right)  \otimes Q_{j}h$ (where this is to be
taken as $0$ if $Q_{j}h=0$ or if $j=1$). Consequently, $\theta_{\eta,\eta
}(Q_{j})h=\eta^{\ast}(Q_{j}\cdot\eta)h=Q_{j-1}h$; i.e., $\theta_{\eta,\eta
}(Q_{j})=Q_{j-1}$. Thus, $(id-\theta_{\eta,\eta})^{-1}(Q_{j})=\sum_{k=0}%
^{n-1}\theta_{\eta,\eta}^{k}(Q_{j})=P_{j}$ and so%
\[
(Ad(B)-Ad(C))\circ(id-\theta_{\eta,\eta})^{-1}(Q_{j})=BPjB^{\ast}%
-CP_{j}C^{\ast}\text{.}%
\]
Since the map $(Ad(B)-Ad(C))\circ(id-\theta_{\eta,\eta})^{-1}$ is defined on
the commutative algebra $\mathcal{D}^{\prime}$, to show that it is completely
positive, we need only show that it is positive. And for this, it clearly
suffices to show that $(Ad(B)-Ad(C))\circ(id-\theta_{\eta,\eta})^{-1}$ maps
each $Q_{j}$ to a positive operator. However, since our hypothesis is that
$BPjB^{\ast}-CP_{j}C^{\ast}\geq0$ for all $j$, this is the case, and the proof
is complete.
\end{proof}

As noted earlier, the following corollary was inspired in particular by
\cite[Theorem 3]{KMT93}. However, there are numerous variations of it in the
literature, dating back at least to \cite{cL69}. For a selection of more
recent analyses see \cite{mA92, AKMT92, JKMT02, eK95}.

\begin{corollary}
\label{Corollary4.8.2}Let $\mathcal{N}$ be a nest of subspaces of the Hilbert
space $H$ and let $u_{1},u_{2},\ldots,u_{m}$ and $v_{1},v_{2},\ldots,v_{m}$ be
$2m$ vectors in $H$. Then there is an operator $X\in Alg\,\mathcal{N}$ with
norm at most one such that $Xu_{i}=v_{i}$, $i=1,2,\cdots,m$, if and only if
for every $N\in\mathcal{N}$, the matrix inequality%
\[
\left(  \langle N^{\perp}v_{i},N^{\perp}v_{j}\rangle\right)  \leq\left(
\langle N^{\perp}u_{i},N^{\perp}u_{j}\rangle\right)
\]
is satisfied in $M_{m}(\mathbb{C)}$, where $N^{\perp}=I-N$.
\end{corollary}

\begin{proof}
Suppose first that we can find the desired $X$. Then $\left(  \langle
N^{\perp}v_{i},N^{\perp}v_{j}\rangle\right)  \allowbreak=\left(  \langle
X^{\ast}N^{\perp}Xu_{i},N^{\perp}u_{j}\rangle\right)  $. However, $X^{\ast
}N^{\perp}X=N^{\perp}X^{\ast}N^{\perp}XN^{\perp}\leq N^{\perp}$. Hence
$\left(  \langle N^{\perp}v_{i},N^{\perp}v_{j}\rangle\right)  =\left(  \langle
X^{\ast}N^{\perp}Xu_{i},N^{\perp}u_{j}\rangle\right)  \leq\left(  \langle
N^{\perp}u_{i},N^{\perp}u_{j}\rangle\right)  $.

Conversely, suppose that $\left(  \langle N^{\perp}v_{i},N^{\perp}v_{j}%
\rangle\right)  \leq\left(  \langle N^{\perp}u_{i},N^{\perp}u_{j}%
\rangle\right)  $ is satisfied for all $N\in\mathcal{N}$. Then in particular,
for $N=0$, we get $\left(  \langle v_{i},v_{j}\rangle\right)  \leq\left(
\langle u_{i},u_{j}\rangle\right)  $. This implies in particular that if $\sum
c_{i}u_{i}=0$, then $\sum c_{i}v_{i}=0$. Therefore, renumbering if necessary,
we may assume that $\{u_{1},u_{2},\ldots,u_{k}\}$ is a maximal linearly
independent set of the $u_{i}$'s and conclude that if we can find an $X\in
Alg\,\mathcal{N}$ with norm at most one such that $Xu_{i}=v_{i}$,
$i=1,2,\ldots k$, then $Xu_{i+1}=v_{i+1}$, $\ldots$, $Xu_{m}=v_{m}$,
automatically. Therefore, we may assume at the outset that all the $u_{i}$'s
are linearly independent and thus that the dimension of $H$ is at least $m$.
Choose and orthonormal family $\{e_{1},e_{2},\ldots,e_{m}\}$ in $H$ and set
$B=\sum e_{i}\otimes u_{i}^{\ast}$ and $C=\sum e_{i}\otimes v_{i}^{\ast}$.
Then for every $N\in\mathcal{N}$, we see that
\begin{align*}
CN^{\perp}C^{\ast}  &  =\sum_{i,j}(e_{i}\otimes v_{i}^{\ast})N^{\perp}%
(v_{j}\otimes e_{j}^{\ast})\\
&  =\sum_{i,j}(e_{i}\otimes v_{i}^{\ast})(N^{\perp}v_{j}\otimes e_{j}^{\ast
})=\sum_{i,j}\langle v_{i},N^{\perp}v_{j}\rangle e_{i}\otimes e_{j}^{\ast}%
\end{align*}
and similarly $BN^{\perp}B^{\ast}=\sum_{i,j}\langle u_{i},N^{\perp}%
u_{j}\rangle e_{i}\otimes e_{j}^{\ast}$. Consequently, the inequality $\left(
\langle N^{\perp}v_{i},N^{\perp}v_{j}\rangle\right)  \leq\left(  \langle
N^{\perp}u_{i},N^{\perp}u_{j}\rangle\right)  $ implies that $CN^{\perp}%
C^{\ast}\leq BN^{\perp}B^{\ast}$ for all $N\in\mathcal{N}$. By Theorem
\ref{Theorem4.8.1}, but with $\mathcal{N}$ replaced by the nest $\mathcal{N}%
^{\perp}:=\{N^{\perp}\mid N\in\mathcal{N}\}$, we can find a $Y\in
Alg\,\mathcal{N}^{\perp}$ ($=(Alg\,\mathcal{N})^{\ast}$) of norm at most one
such that $BY=C$. If we set $X=Y^{\ast}$ then $XB^{\ast}=C^{\ast}$. Since
$B^{\ast}e_{j}=\sum(u_{i}\otimes e_{i}^{\ast})e_{j}=u_{j}$, we conclude that
$Xu_{j}=XB^{\ast}e_{j}=C^{\ast}e_{j}=v_{j}$.
\end{proof}

\begin{remark}
\label{Remark4.8.2}One can use the Nevanlinna-Pick theorem to compute distance
formulas. For example, one can prove Arveson's distance formula \cite{wA75}
for nest algebras. (It is known that it suffices to prove the distance formula
for finite nests.) More generally, one can compute the distance of an element
in a (finite) nest algebra from an ideal. Since the results are well known, we
will not spell out all the details, but just point out how one can use Theorem
\ref{Theorem3.1}. So, let $\mathcal{N}=\{0=P_{0}<P_{1}<P_{2}<\cdots<P_{n}=I\}$
and let $J$ be a (2-sided) ideal in $Alg\,\mathcal{N}$. Then $J$ can be
written as $\cap\{\ker(\Psi_{G})\}$, where $G$ runs over a certain finite set
of intervals of the nest $\mathcal{N}$ and where $\Psi_{G}$ is the map given
by the formula $\Psi_{G}(T)=GTG$, $T\in Alg\,\mathcal{N}$. Let $\Psi
=\sum^{\oplus}\Psi_{G}$, a representation of $Alg\,\mathcal{N}$ on
$\sum^{\oplus}GH$. We want to show that for $T$ in $Alg\,\mathcal{N}$,
$dist(T,J)=\left\|  \Psi(T)\right\|  $. For this, it will suffice to prove
that given $T\in Alg\,\mathcal{N}$ with $\left\|  \Psi(T)\right\|  \leq1$ then
one can find an $X\in Alg\,\mathcal{N}$ with norm at most one such that
$\Psi(X)=\Psi(T)$. By Theorem \ref{Theorem3.1}, a necessary and sufficient
condition for this is the complete positivity of the map%
\[
(id-Ad(\Psi(T)))\circ(id-\theta_{\eta,\eta})^{-1}%
\]
where $\eta$ is the vector associated with the representation $\Psi$ of
$Alg\,\mathcal{N}$ $=H^{\infty}(E)$, as in Theorem \ref{Theorem4.8.1}.
(Note:$(id-Ad(\Psi(T)))\circ(id-\theta_{\eta,\eta})^{-1}$ is defined on
$\Psi(\mathcal{D})^{\prime}$.) Straightforward but tedious calculations reveal
that it is; we omit the details.
\end{remark}

\section{CNC representations\label{CNCReps}}

Throughout this section we will continue with our standing assumption that $M$
is a von Neumann algebra and we will fix a $W^{\ast}$-correspondence $E$ over
$M$. As we pointed out in the Introduction and in Remark \ref{keyproblem}, not
all completely contractive representations of $\mathcal{T}_{+}(E)$ extend to
$H^{\infty}(E)$. We showed in Theorem \ref{lem1.10Cor1.11} that if
$\sigma\times T$ is a completely contractive representation of $\mathcal{T}%
_{+}(E)$ on a Hilbert space $H$, where $\sigma$ is assumed to be a normal
representation of $M$ on $H$, then $\sigma\times T$ extends to $H^{\infty}(E)$
whenever $\left\|  \tilde{T}\right\|  <1$. In this section, we want to extend
this result to certain representations $\sigma\times T$, where $\left\|
\tilde{T}\right\|  =1$ but are otherwise constrained by a condition that we
call ``completely non-coisometric''. We believe this is the greatest level of
generality that can be achieved with current technology. We also want to
extend our generalization of the Nevanlinna-Pick theorem to accommodate such representations.

To understand the idea behind the notion, recall that if $T$ is a contraction
operator on a Hilbert space $H$, then the set%
\[
H_{1}:=\{h\in H\mid\left\|  (T^{\ast})^{k}h\right\|  =\left\|  h\right\|
\text{, for all }k\}
\]
is a subspace of $H$ that is invariant under $T^{\ast}$ and is the largest
subspace of $H$ on which $T^{\ast}$ acts like an isometry. Alternatively,
$H_{1}$ is the intersection of the kernels of the operators $\Delta_{\ast
}T^{\ast k}$, $k=1,2,\ldots$, where $\Delta_{\ast}=(I-TT^{\ast})^{1/2}$. The
operator $T$ is called completely non-coisometric in case $H_{1}=0$. Observe
that if the powers of $T^{\ast}$ tend strongly to zero, then $T$ is completely
non-coisometric. So, in particular, all strict contractions are completely
non-coisometric. The importance of this notion for us is that it generalizes
to our situation and it does so in such a way that can be detected in the
geometry of the dilation spaces we consider. This was observed first, for
single contractions, in the work of Sz.-Nagy and Foia\c{s} (see \cite{szNF70}%
). Subsequently, in the study of row contractions, Popescu identified the
relevant geometric facts in \cite[Proposition 2.9]{gP89}. Thus, we are led to
propose the following definition.

\begin{definition}
\label{complnoncoisorep}Suppose $(T,\sigma)$ is a completely contractive
covariant representation of $E$ on the Hilbert space $H$. We call $(T,\sigma)$
\emph{completely non-coisometric} in case the subspace
\[
H_{1}:=\cap_{n=0}^{\infty}\ker((I_{E^{\otimes n}}\otimes\Delta_{\ast}%
)\tilde{T}_{n}^{\ast})
\]
is the zero subspace, where $\Delta_{\ast}:=(I-\tilde{T}\tilde{T}^{\ast
})^{1/2}$ and where, recall, $\tilde{T}_{n}$ denotes the $n^{th}$ generalized
power of $\tilde{T}$ (formula (\ref{iterates})). We shall abreviate
``completely non-coisometric'' by ``cnc''. The space $H_{1}$, itself, will be
called the \emph{coisometric subspace} of $H$.
\end{definition}

\begin{remark}
\label{altcnc}The space $H_{1}$ has alternate descriptions:%
\begin{align*}
H_{1}  &  =\cap_{n=0}^{\infty}\ker((I_{E^{\otimes n}}\otimes\Delta_{\ast
})\tilde{T}_{n}^{\ast})\\
&  =\cap_{n=1}^{\infty}\{\ker\Delta_{\ast}T(\xi_{1})^{\ast}T(\xi_{2})^{\ast
}\cdots T(\xi_{n})^{\ast}\mid\xi_{1},\xi_{2},\ldots\xi_{n}\in E\}\\
&  =\{h\in H\mid\left\|  \tilde{T}_{n}^{\ast}h\right\|  =\left\|  h\right\|
\text{, for all }n\}\text{.}%
\end{align*}
Further, $H_{1}$ is the largest subspace of $H$ that is invariant under
$(T\times\sigma)(\mathcal{T}_{+}(E))^{\ast}$ and on which $\tilde{T}^{\ast}$
acts isometrically.
\end{remark}

Our primary objective is to prove the following theorem, which, among other
things, generalizes a result due to Popescu \cite[Theorem 4.3]{gP95a}.

\begin{theorem}
\label{cncext}Let $(T,\sigma)$ be a completely contractive covariant
representation of $E$ on the Hilbert space $H$ and assume that $(T,\sigma)$ is
completely non-coisometric. Then

\begin{enumerate}
\item $T\times\sigma$ can be extended to an ultraweakly continuous, completely
contractive representation of the Hardy algebra $H^{\infty}(E)$ on $H$.

\item  If $T_{r}=rT,$ $0\leq r<1$, then each $(T_{r},\sigma)$ is a completely
contractive covariant representation of $E$ on $H$, with $\left\|  \tilde
{T}_{r}\right\|  =r\left\|  \tilde{T}\right\|  <1$, and for all $X\in
H^{\infty}(E)$,%
\[
\lim_{r\rightarrow1}(T_{r}\times\sigma)(X)=(T\times\sigma)(X)
\]
in the strong operator topology on $B(H)$.

\item  If $\{X_{r}\}_{0\leq r<1}$ is a bounded net in $H^{\infty}(E)$
converging to $X$ in the ultraweak topology on $H^{\infty}(E)$, then%
\[
\lim_{r\rightarrow1}(T_{r}\times\sigma)(X_{r})=(T\times\sigma)(X)
\]
in the weak operator topology on $B(H)$.
\end{enumerate}
\end{theorem}

The proof of this result requires some analysis of the geometry of the space
for the minimal isometric dilation of $T\times\sigma$. This can be developed
in a series of straightforward lemmas. In order not to lose sight of the
connection between Theorem \ref{cncext} and our Nevanlinna-Pick analysis, we
postpone the proof until after reaping two corollaries: an enhancement of our
Nevanlinna-Pick theorem and the accompanying Schwarz lemma.

\begin{theorem}
\label{Theorem2.1}Let $E$ be a $W^{\ast}$-correspondence over the von Neumann
algebra $M$ and let $\sigma$ be a representation of $M$ on the Hilbert space
$H$. Given $\eta_{1}^{\ast},\eta_{2}^{\ast},\cdots,\eta_{k}^{\ast}$ in the
\emph{closure} of $\mathbb{D((}E^{\sigma})^{\ast})$ such that $(\sigma
,\eta_{i}^{\ast})$ are completely non-coisometric and given two $k$-tuples of
operators in $B(H)$, $B_{1},B_{2},\cdots,B_{k}$ and $C_{1},C_{2},\cdots,C_{k}%
$, then there is an element $X\in H^{\infty}(E)$ such that $\left\|
X\right\|  \leq1$ and such that%
\begin{equation}
B_{i}X(\eta_{i}^{\ast})=C_{i}\text{,} \label{Cond1}%
\end{equation}
$i=1,2,\cdots k$, if and only if for each $i$ there is a net $\{C_{i}%
(r)\mid0\leq r<1\}$ in $B(H)$ converging strongly to $C_{i}$ as $r\rightarrow
1$ such that for every $r$, $0\leq r<1$, the map from $M_{k}(\sigma
(M)^{\prime})$ to $M_{k}(B(H))$ defined by the $k\times k$ matrix in
$M_{k}(B(\sigma(M)^{\prime},B(H)))$,%
\begin{equation}
\left(  \left(  Ad(B_{i},B_{j})-Ad(C_{i}(r),C_{j}(r))\right)  \circ\left(
id-r^{2}\theta_{\eta_{i},\eta_{j}}\right)  ^{-1}\right)  \label{Cond2}%
\end{equation}
is completely positive.

Moreover, if $N$ is a von Neumann subalgebra of $M$, if for each $i$,
$i=1,2,\cdots,k$, $C_{i}$ and $B_{i}$ lie in the von Neumann algebra generated
by $\sigma(N)$ and $\sigma(M)^{\prime}$, $\sigma(N)\bigvee\sigma(M)^{\prime}$,
and if for each $i$ there is a net $\{C_{i}(r)\mid0\leq r<1\}$ in
$\sigma(N)\bigvee\sigma(M)^{\prime}$ converging strongly to $C_{i}$ such that
for each $r$, $0\leq r<1$, the matrix (\ref{Cond2}) represents a completely
positive operator, a solution $X$ to the interpolation equation (\ref{Cond1})
can be found in $H^{\infty}(E)$ that commutes with $\varphi_{\infty}(M\cap
N^{\prime})^{\prime}$.
\end{theorem}

\begin{proof}
Suppose nets $\{C_{i}(r)\mid0\leq r<1\}$, $i=1,2,\ldots,k$, can be found such
that the matrix (\ref{Cond2}) is completely positive for all $r$, $0\leq r<1$.
Then we may apply Theorem \ref{Theorem3.1} to the representations $(r\eta
_{i}^{\ast},\sigma)$ and the operators $B_{i}$ and $C_{i}(r)$ in $B(H)$ to
find operators $\{X_{r}\mid0\leq r<1\}$ in $H^{\infty}(E)$, all of norm at
most $1$, such that $B_{i}X_{r}(r\eta_{i}^{\ast})=C_{i}(r)$ for all
$i=1,2,\ldots,k$. Since the $X_{r}$ are uniformly bounded by $1$, we may pass
to a subnet, if necessary, and assume that they converge to an element $X\in
H^{\infty}(E)$ in the ultraweak topology. Then, applying the third assertion
of Theorem \ref{cncext}, we conclude that for all $i=1,2,\ldots,k$,
\[
B_{i}X(\eta_{i}^{\ast})=C_{i}\text{.}%
\]
Note that if the operators $B_{i}$, $C_{i}$, and $C_{i}(r)$, $0\leq r<1$, all
lie in $\sigma(N)\bigvee\sigma(M)^{\prime}$, then Theorem \ref{Theorem3.1}
guarantees that the $X_{r}$ may all be chosen to commute with $\varphi
_{\infty}(M\cap N^{\prime})^{\prime}$ in $H^{\infty}(E)$. Any limit, $X$, then
also commutes with $\varphi_{\infty}(M\cap N^{\prime})^{\prime}$.

Conversely, suppose there is an $X\in H^{\infty}(E)$ that satisfies equation
(\ref{Cond1}) for all $i$, $i=1,2,\ldots,k$ and set $C_{i}(r)=X(r\eta
_{i}^{\ast})=(r\eta_{i}^{\ast}\times\sigma)(X)$. Then the strong convergence
of $C_{i}(r)$ to $C_{i}$ follows from the second condition in Theorem
\ref{cncext}, \ and the complete positivity of the matrix (\ref{Cond2}) is a
consequence of Theorem \ref{Theorem3.1}.
\end{proof}

The ``Schwartz lemma'', Theorem \ref{Theorem3.5}, can be extended to the
setting of completely non-coisometric representations, too.

\begin{theorem}
\label{Theorem2.2}Suppose an element $X$ of $H^{\infty}(E)$ has norm at most
one and satisfies the equation $X(0)=0$. Then for every $\eta^{\ast}%
\in(E^{\sigma})^{\ast}$ such that $(\eta^{\ast},\sigma)$ is completely
non-coisometric the following assertions are valid:

\begin{enumerate}
\item  If $a$ is a nonnegative element in $\sigma(M)^{\prime}$, and if
$\langle\eta,a\cdot\eta\rangle\leq a$, then
\[
X(\eta^{\ast})aX(\eta^{\ast})^{\ast}\leq\langle\eta,a\cdot\eta\rangle.
\]

\item  If $\eta^{\otimes k}$ denotes the element $\eta\otimes\eta\otimes
\cdots\otimes\eta\in E^{\otimes k}$, then
\[
X(\eta^{\ast})\langle\eta^{\otimes k},\eta^{\otimes k}\rangle X(\eta^{\ast
})^{\ast}\leq\langle\eta^{\otimes k+1},\eta^{\otimes k+1}\rangle\text{.}%
\]

\item $X(\eta^{\ast})X(\eta^{\ast})^{\ast}\leq\langle\eta,\eta\rangle$.
\end{enumerate}
\end{theorem}

\begin{proof}
As in the proof of Theorem \ref{Theorem3.5}, we need only attend to the first
assertion. For this purpose, consider $X(r\eta^{\ast})=(r\eta^{\ast}%
\times\sigma)(X)$, $0\leq r<1$, and suppose $a$ is a nonnegative element of
$\sigma(M)^{\prime}$ satisfying the inequality $\langle\eta,a\cdot\eta
\rangle\leq a$. Then certainly $r^{2}\langle\eta,a\cdot\eta\rangle\leq a$, and
we may apply Theorem \ref{Theorem3.5} to conclude that $X(r\eta^{\ast
})aX(r\eta^{\ast})^{\ast}\leq r^{2}\langle\eta,a\cdot\eta\rangle=r^{2}%
\eta^{\ast}(I_{E}\otimes a)\eta$ for all $r$. If we let $b$ be the positive
square root of $a$, then this inequality guarantees that there is a $Y(r)$ of
norm at most $1$ such that $X(r\eta^{\ast})b=r\eta^{\ast}(I_{E}\otimes
b)Y(r)$. Passing to a subnet of $\{Y(r)\}_{0\leq r<1}$, if necessary, we may
assume that there is a contraction $Y$ to which $\{Y(r)\}_{0\leq r<1}$
converges weakly. It follows that $X(\eta^{\ast})b=\eta^{\ast}(I_{E}\otimes
b)Y$, which in turn implies that $X(\eta^{\ast})aX(\eta^{\ast})^{\ast}\leq
\eta^{\ast}(I_{E}\otimes a)\eta=\langle\eta,a\cdot\eta\rangle$.
\end{proof}

We turn now to the details of the proof of Theorem \ref{cncext}. For this
purpose, we fix a completely contractive covariant representation $(T,\sigma)$
of $E$ on the Hilbert space $H$ and we let $(V,\rho)$ be its minimal isometric
dilation acting on $K=H\oplus\mathcal{F}(E)\otimes_{\sigma_{1}}\mathcal{D}$.
We will follow the notation developed immediately after Theorem \ref{Theorem
1.13}. As we have noted in Remark \ref{cpendo}, $V$ determines an endomorphism
$L$ of the von Neumann algebra $\rho(M)^{\prime}$ via the formula
\[
L(x)=\tilde{V}(I_{E}\otimes x)\tilde{V}^{\ast}\text{,}%
\]
$x\in\rho(M)^{\prime}$ and the powers of $L$, $\{L^{n}\}_{n\geq0}$ are given
by the formula%
\[
L^{n}(x)=\tilde{V}_{n}(I_{E^{\otimes n}}\otimes x)\tilde{V}_{n}^{\ast}\text{,}%
\]
$x\in\rho(M)^{\prime}$, where $\tilde{V}_{n}$ denotes the generalized $n^{th}$
power of $\tilde{V}$. The following terminology comes from \cite{MS99}. The
notation, however, is a bit different.

\begin{definition}
\label{wandering}A subspace $\frak{M}$ of $K$ is called a \emph{wandering
subspace }of $K$ for $V$ if it is invariant under $\rho(M)$, so the projection
onto it, $P_{\frak{M}}$, lies in $\rho(M)^{\prime}$, and if the projections
$\{L^{n}(P_{\frak{M}})\}_{n\geq0}$ are mutually orthogonal. In this event, we
write $L_{\infty}(\frak{M)}$ for the range of $\sum_{n=0}^{\infty}%
L^{n}(P_{\frak{M}})$.
\end{definition}

Following \cite{MS99}, we shall write $P_{n}$ for the image of $I_{K}$ under
$L^{n}$, i.e., $P_{n}=L^{n}(I_{K})=\tilde{V}_{n}\tilde{V}_{n}^{\ast}$. Then
$\{P_{n}\}_{n\geq0}$ is a decreasing sequence of projections in $\rho
(M)^{\prime}$ whose infimum we denote by $P_{\infty}$. Also, we set
$Q_{n}=P_{n}-P_{n+1}$, $n\geq0$. Then the range $\frak{Q}_{0}$ of $Q_{0}$ is a
wandering subspace and $L_{\infty}(\frak{Q}_{0})=(I-P_{\infty})(K)$. The
spaces $L_{\infty}(\frak{Q}_{0})$ and $(I-P_{\infty})(K)$ reduce $(V,\rho)$.
If $\rho_{0}$ is given by the formula $\rho_{0}(a)=\rho(a)|\frak{Q}_{0}$,
$a\in M$, then the restriction of $(V,\rho)$ to $L_{\infty}(\frak{Q}_{0})$ is
unitarily equivalent to $(V_{ind},\rho_{ind})$ on $\mathcal{F}(E)\otimes
_{\rho_{0}}\frak{Q}_{0}$ where, recall, $V_{ind}(\xi)=T_{\xi}\otimes
I_{\frak{Q}_{0}}$ and $\rho_{ind}(a)=\varphi_{\infty}(a)\otimes I_{\frak{Q}%
_{0}}$. The restriction, $(V_{\infty},\rho_{\infty})$, of $(V,\rho)$ to
$(I-P_{\infty})(K)$ is fully coisometric. This decomposition, $(V,\rho
)=(V_{ind},\rho_{ind})\oplus(V_{\infty},\rho_{\infty})$, is called the
\emph{Wold decomposition} of $(V,\rho)$ and its existence is proved and other
properties are developed in \cite[Theorem 2.9]{MS99}.

\begin{remark}
\label{EnhancedWold}Of course, if $\frak{M}\subseteq K$ is a wandering
subspace for $(V,\rho)$, then $L_{\infty}(\frak{M})$ is invariant under
$(V,\rho)$, as may be seen from Lemma 2.7 of \cite{MS99} and the discussion
just presented shows that restriction of $(V,\rho)$ to $L_{\infty}(\frak{M})$
is unitarily equivalent to the representation induced by the restriction of
$\rho$ to $\frak{M}$. See, in particular, Corollary 2.10 of \cite{MS99}.
\end{remark}

\begin{lemma}
\label{Lemma1.2}In regard to the decomposition of the space $K$ as
$H\oplus\mathcal{F}(E)\otimes_{\sigma_{1}}\mathcal{D}$, we have%
\[
\frak{Q}_{0}=\overline{Q_{0}(H)}=\overline{\{\Delta_{\ast}^{2}h\oplus
(-\Delta\tilde{T}^{\ast}h)\mid h\in H\}}\subseteq H\oplus\mathcal{D}\text{,}%
\]
where $\Delta_{\ast}:=(I_{H}-\tilde{T}\tilde{T}^{\ast})^{1/2}$.
\end{lemma}

\begin{proof}
From the minimality of $K$, $I=\bigvee_{n=0}^{\infty}L^{n}(P_{H})=P_{H}\bigvee
P_{1}$. Since $Q_{0}=I-P_{1}$, we see that $\frak{Q}_{0}=\overline{Q_{0}(H)}$.
On the other hand, the matricial representation of $\tilde{V}$ (\ref{Vtilde})
gives the equation%
\[
Q_{0}=I_{K}-\tilde{V}\tilde{V}^{\ast}=\left(
\begin{array}
[c]{cccc}%
I_{H}-\tilde{T}\tilde{T}^{\ast} & -\tilde{T}\Delta & 0 & \ldots\\
-\Delta\tilde{T}^{\ast} & I_{\mathcal{D}}-\Delta^{2} & 0 & \\
0 & 0 & 0 & \\
\vdots &  &  & \ddots
\end{array}
\right)
\]
from which the rest of the lemma follows.
\end{proof}

The following is the key geometric fact to which we alluded above. It
generalizes analyses in \cite{szNF70} and Proposition 2.9 of \cite{gP89}.

\begin{proposition}
\label{Lemma1.4}The subspaces $\mathcal{D}$ and $\frak{Q}_{0}$ are wandering
subspaces of $K$ and the coisometric subspace $H_{1}$ of $H$, viewed as
contained in $K$ is the orthogonal complement of $L_{\infty}(\mathcal{D}%
)\bigvee L_{\infty}(\frak{Q}_{0})$. Thus $(T,\sigma)$ is completely
non-coisometric if and only if $K=L_{\infty}(\mathcal{D})\bigvee L_{\infty
}(\frak{Q}_{0})$.
\end{proposition}

\begin{proof}
It is clear from the definition of $K$ that $\mathcal{D}$ is a wondering
subspace and that $L_{\infty}(\mathcal{D})=\mathcal{F}(E)\otimes_{\sigma_{1}%
}\mathcal{D}=H^{\perp}$. Of course, $\frak{Q}_{0}$ is wandering, as we have
noted before. The key point of the proposition is the relation between $H_{1}$
and $L_{\infty}(\mathcal{D})\bigvee L_{\infty}(\frak{Q}_{0})$. Since $K\ominus
L_{\infty}(\mathcal{D)}=H$, we may write%
\begin{multline*}
K\ominus(L_{\infty}(\mathcal{D})\bigvee L_{\infty}(\frak{Q}_{0}))=H\cap
(K\ominus L_{\infty}(\frak{Q}_{0}))=H\cap P_{\infty}(K)\\
=\{h\in H\mid\langle h,V(\xi_{1})V(\xi_{2})\cdots V(\xi_{n})k\rangle=0\text{,
}\forall\xi_{i}\in E,\text{ }n\geq0,\;\forall k\in\frak{Q}_{0}\}\text{.}%
\end{multline*}
Since $V(\xi)^{\ast}$ leaves $H$ invariant and $V(\xi)^{\ast}|H=T(\xi)^{\ast}%
$, by Theorem \ref{Theorem 1.13}, we conclude that $H\cap P_{\infty}(K)$
consists of all vectors $h\in H$ such that
\[
0=\langle h,V(\xi_{1})V(\xi_{2})\cdots V(\xi_{n})k\rangle=\langle T(\xi
_{n})^{\ast}T(\xi_{n-1})^{\ast}\cdots T(\xi_{1})^{\ast}h,P_{H}k\rangle
\]
for all choices of $\xi_{i}$'s in $E$ and for all $k\in\frak{Q}_{0}$. However,
by Lemma \ref{Lemma1.2}, $\overline{P_{H}\frak{Q}_{0}}=\overline{\Delta_{\ast
}^{2}H}$, and this completes the proof.
\end{proof}

\begin{lemma}
\label{Lemma1.5}A completely contractive covariant representation $(T,\sigma)$
is completely non-coisometric if and only if $P_{\infty}(K)=\overline
{P_{\infty}(L_{\infty}(\mathcal{D}))}$ if and only if $\overline{P_{\infty
}(H)}\subseteq\overline{P_{\infty}(L_{\infty}(\mathcal{D}))}$.
\end{lemma}

\begin{proof}
The second ``if and only if'' is obvious since $H=K\ominus L_{\infty
}(\mathcal{D})$. For the first, note that if $(T,\sigma)$ is completely
non-coisometric, so that $K=L_{\infty}(\mathcal{D})\bigvee L_{\infty}%
(\frak{Q}_{0})$, by Proposition \ref{Lemma1.4}, then since $P_{\infty}$
vanishes on $L_{\infty}(\frak{O}_{0})$, $P_{\infty}(K)=\overline{P_{\infty
}(L_{\infty}(\mathcal{D}))}$. For the converse, observe that if $P_{\infty
}(K)=\overline{P_{\infty}(L_{\infty}(\mathcal{D}))}$, then $P_{\infty}%
(H_{1})\subseteq\overline{P_{\infty}(L_{\infty}(\mathcal{D}))}$. So, if $h\in
H_{1}$, we may write $P_{\infty}h=\lim P_{\infty}(g_{n})$ for a sequence
$\{g_{n}\}\subseteq L_{\infty}(\mathcal{D})$. Since $h$ is orthogonal to
$L_{\infty}(\frak{Q}_{0})$, by Proposition \ref{Lemma1.4}, we have
$h=P_{\infty}h$ and $g_{n}-(I-P_{\infty})g_{n}=P_{\infty}g_{n}\rightarrow h$.
However, each $g_{n}-(I-P_{\infty})g_{n}$ lies in $L_{\infty}(\mathcal{D}%
)+L_{\infty}(\frak{Q}_{0})\subseteq K\ominus H_{1}$, again by Proposition
\ref{Lemma1.4}, and we conclude that $h=0$. Thus $H_{1}=\{0\}$ and
$(T,\sigma)$ is completely non-coisometric.
\end{proof}

Next, we show that if $\frak{M}\subseteq K$ is a wandering subspace for
$(V,\rho)$, then the restriction of $(V,\rho)$ to $L_{\infty}(\frak{M})$
satisfies the conclusion of Theorem \ref{cncext}.

\begin{lemma}
\label{Lemma1.7a}Suppose $\frak{M}\subseteq K$ is a wandering subspace for
$(V,\rho)$ and let $(V_{1},\rho_{1})$ denote the restriction of $(V,\rho)$ to
$L_{\infty}(\frak{M})$. Then:

\begin{enumerate}
\item $V_{1}\times\rho_{1}$ can be extended to an ultraweakly continuous,
completely contractive representation of the Hardy algebra $H^{\infty}(E)$ on
$L_{\infty}(\frak{M})$.

\item  If $V_{1r}=rV_{1},$ $0\leq r<1$, then each $(V_{1r},\rho_{1})$ is a
completely contractive covariant representation of $E$ on $H$, with $\left\|
\tilde{V}_{1r}\right\|  =r<1$ (so $(V_{1r}\times\rho)$ can be extended to all
of $H^{\infty}(E)$), and for all $X\in H^{\infty}(E)$,%
\[
\lim_{r\rightarrow1}(V_{1r}\times\rho)(X)=(V_{1}\times\rho_{1})(X)
\]
in the strong operator topology on $B(L_{\infty}(\frak{M}))$.

\item  If $\{X_{r}\}_{0\leq r<1}$ is a bounded net in $H^{\infty}(E)$
converging to $X$ in the ultraweak topology on $H^{\infty}(E)$, then%
\[
\lim_{r\rightarrow1}(V_{1r}\times\rho_{1})(X_{r})=(V_{1}\times\rho_{1})(X)
\]
in the weak operator topology on $B(L_{\infty}(\frak{M}))$.
\end{enumerate}
\end{lemma}

\begin{proof}
As noted in Remark \ref{EnhancedWold}, $(V_{1},\rho_{1})$ is (unitarily
equivalent to) the representation induced by the restriction $\rho_{0}$ of
$\rho$ to $\frak{M}$, $\rho_{0}^{\mathcal{F}(E)}$. That is, up to unitary
equivalence, $\rho=\rho_{0}^{\mathcal{F}(E)}\circ\varphi_{\infty}$ and
$V_{1}(\xi)=\rho_{0}^{\mathcal{F}(E)}(T_{\xi})$. So we will assume that
$(V_{1},\rho_{1})$ has this form. Since $\rho_{0}^{\mathcal{F}(E)}$ is a
normal representation of $\mathcal{L}(\mathcal{F}(E))$ (see \cite{mR74b}), its
restriction to $H^{\infty}(E)$ is an ultraweakly continuous, completely
contractive representation of $H^{\infty}(E)$ and, evidently, it extends
$V_{1}\times\rho_{1}$. This proves (1).

For (2) recall that each $X\in H^{\infty}(E)$ has a Fourier development
$X=\sum_{n=0}^{\infty}X_{n}$, where $X_{n}=\Phi_{n}(X)$ and $\Phi_{n}$ is
given by formula (\ref{shift}). As we noted, this series is Cesaro summable to
$X$ in the ultraweak topology on $\mathcal{L}(\mathcal{F}(E))$. In fact, owing
to the implementation of $\gamma_{t}$ by the unitary group (\ref{unitarygroup}%
), a straightforward calculation reveals that
\[
V_{1r}\times\rho_{1}(X)=\sum_{n=0}^{\infty}r^{n}(V_{1}\times\rho_{1}%
)(X_{n})\text{,}%
\]
$0\leq r<1$. The strong continuity of the unitary group (\ref{unitarygroup})
in this induced representation guarantees that the series representing $X$ is
Abel summable in the strong operator topology in the space of this
representation. That is, $V_{1r}\times\rho_{1}(X)\rightarrow V_{1}\times
\rho_{1}(X)$ in the strong operator topology. Alternatively, for $0\leq r<1$,
define $R(r)$ on $B(\mathcal{F}(E)\otimes_{\rho_{0}}\frak{M})$ by the formula%
\[
R(r)(\xi_{1}\otimes\xi_{2}\otimes\cdots\xi_{n}\otimes h)=r^{n}\xi_{1}%
\otimes\xi_{2}\otimes\cdots\xi_{n}\otimes h\text{,}%
\]
$\xi_{1}\otimes\xi_{2}\otimes\cdots\xi_{n}\otimes h\in E^{\otimes n}%
\otimes_{\rho_{0}}\frak{M}$. Then $R(r)\rightarrow I$ strongly on
$\mathcal{F}(E)\otimes_{\rho_{0}}\frak{M}$. A straightforward calculation
shows that for all $g\in E^{\otimes n}\otimes_{\rho_{0}}\frak{M}$,%
\[
(V_{1r}\times\rho_{1})(X)g=r^{-n}R(r)(V_{1}\times\rho_{1})(X)g\rightarrow
(V_{1}\times\rho_{1})(X)g\text{.}%
\]
Since the set of operators $V_{1r}\times\rho_{1}(X)$ is uniformly bounded by
the norm of $X$, we conclude that $V_{1r}\times\rho_{1}(X)g\rightarrow
(V_{1}\times\rho_{1})(X)g$ for all $g\in\mathcal{F}(E)\otimes_{\rho_{0}%
}\frak{M}$.

Finally, to prove (3), fix a bounded net $\{X_{r}\}_{0\leq r<1}$ in
$H^{\infty}(E)$ that converges to $X$ ultraweakly. Then for all $g\in
E^{\otimes n}\otimes_{\rho_{0}}\frak{M}$ and all $h\in E^{\otimes k}%
\otimes_{\rho_{0}}\frak{M}$ we have
\begin{multline*}
\langle(V_{1r}\times\rho_{1})(X_{r})g,h\rangle=\langle r^{-n}R(r)(V_{1}%
\times\rho_{1})(X_{r})g,h\rangle\\
=r^{k-n}\langle\rho_{0}^{\mathcal{F}(E)}(X_{r})g,h\rangle\rightarrow
\langle\rho_{0}^{\mathcal{F}(E)}(X)g,h\rangle=\langle(V_{1}\times\rho
_{1})(X)g,h\rangle\text{.}%
\end{multline*}
Since the net $\{(V_{1r}\times\rho_{1})(X_{r})\}_{0\leq r<1}$ is uniformly
bounded, it must converge weakly to $(V_{1}\times\rho_{1})(X)$.
\end{proof}

\begin{lemma}
\label{Lemma1.8a}In the notation of Lemma \ref{Lemma1.7a}, suppose $S$ is a
bounded operator from $L_{\infty}(\frak{M})$ to $K$ that intertwines
$(V_{1},\rho_{1})$ and $(V,\rho)$, i.e., suppose $S(V_{1}\times\rho
_{1})(X)=(V\times\rho)(X)S$ for all $X\in\mathcal{T}_{+}(E)$, and let
$\frak{N}$ denote the space $\overline{SL_{\infty}(\frak{M})}$. Then
$\frak{N}$ is invariant under $(V,\rho)$ and the restriction of $(V,\rho)$ to
$\frak{N}$ is an isometric representation. If $(W,\tau)$ denotes this
restriction, then $W\times\tau$ admits a unique extension to a representation
of $H^{\infty}(E)$ on $\frak{N}$ that satisfies the conclusions (1), (2) and
(3) of Lemma \ref{Lemma1.7a}.
\end{lemma}

\begin{remark}
It is tempting to assert that $(W,\tau)$ is an induced representation, so that
the result follows immediately from Lemma \ref{Lemma1.7a}. However, $(W,\tau)$
need not be induced. One may find examples in the classical setting of a
single contraction. There, the matter reduces to the question: Does there
exist a unitary operator $W$ and an operator $S$ with zero kernel and dense
range such that $WS=SU_{+}$, where $U_{+}$ is the unilateral shift? The answer
is easily seen to be ``yes''. Just let $X$ be a subset of $\mathbb{T}$ of
measure different from $0$ and $1$, let $W$ be multiplication by $z$ on
$L^{2}(X)$ and let $S$ be the operator of multiplication by the indicator
function, $1_{X}$, viewed as an operator from $H^{2}(\mathbb{T})$ to
$L^{2}(X)$. Then since nonzero $H^{2}(\mathbb{T})$-functions cannot vanish on
sets of positive measure, it is easy to see that $S$ and $W$ have the desired properties.
\end{remark}

\begin{proof}
(of Lemma \ref{Lemma1.8a}) It is clear that $\frak{N}$ is invariant under
$(V,\rho)$ and that the restriction, $(W,\tau)$, is an isometric
representation. Further, since $S(V_{1}\times\rho_{1})(X)=(V\times\rho)(X)S$
for all $X\in\mathcal{T}_{+}(E)$, we may write
\begin{equation}
S(V_{1}\times\rho_{1})(X)g=(W\times\tau)(X)Sg\label{intertwining}%
\end{equation}
for all $X\in\mathcal{T}_{+}(E)$ and all $g\in L_{\infty}(\frak{M})$. We find,
then, that $S(V_{1r}\times\rho_{1})(X)g=(W_{r}\times\tau)(X)Sg$ for all
$X\in\mathcal{T}_{+}(E)$, all $g\in L_{\infty}(\frak{M})$, and all $r$, $0\leq
r<1$. By Lemma \ref{Lemma1.7a}, $(W_{r}\times\tau)(X)Sg\rightarrow(W\times
\tau)(X)Sg$ for all $X\in\mathcal{T}_{+}(E)$, all $g\in L_{\infty}(\frak{M})$.
However, $(W_{r}\times\tau)(X)$ makes sense for all $X\in H^{\infty}(E)$ by
Corollary \ref{Corollary1.12} and the equation $S(V_{1r}\times\rho
_{1})(X)g=(W_{r}\times\tau)(X)Sg$ holds for all $g\in L_{\infty}(\frak{M)}$
and all $X\in H^{\infty}(E)$. Thus we may define $(W\times\tau)(X)Sg$ to be
$S(V_{1}\times\rho_{1})(X)g$. Since%
\begin{multline*}
\left\|  S(V_{1}\times\rho_{1})(X)g\right\|  =\lim_{r\rightarrow1}\left\|
(W_{r}\times\tau)(X)Sg\right\|  \\
\leq\lim\sup\left\|  W_{r}\times\tau\right\|  \left\|  X\right\|  \left\|
Sg\right\|  \leq\left\|  X\right\|  \left\|  Sg\right\|  \text{,}%
\end{multline*}
we see that $(W\times\tau)(X)$ extends to a bounded operator on all of
$\frak{N}$ of norm dominated by $\left\|  X\right\|  $. This shows that
$W\times\tau$ extends to $H^{\infty}(E)$ \emph{and} that equation
(\ref{intertwining}) is satisfied for all $X$ in $H^{\infty}(E)$. This
equation, in turn, shows that $(W\times\tau)$ is an algebra homomorphism,
which must be contractive. Further, since $(W_{r}\times\tau)(X)Sg\rightarrow
(W\times\tau)(X)Sg$ for all $X\in H^{\infty}(E)$ and all $g\in L_{\infty
}(\frak{M})$, the fact that $\left\|  (W_{r}\times\tau)\right\|  \leq1$ for
all $r$ implies that $(W_{r}\times\tau)(X)\rightarrow(W\times\tau)(X)$
strongly on $\frak{N}$ for all $X\in H^{\infty}(E)$. That is, condition (2)
from Lemma \ref{Lemma1.7a} is satisfied. Since $W\times\tau$ is the strong
limit of the \emph{completely contractive} maps $W_{r}\times\tau$,
$W\times\tau$ must itself be a completely contractive representation of
$H^{\infty}(E)$.

To prove condition (3) from Lemma \ref{Lemma1.7a}, take $g,h\in L_{\infty
}(\frak{M})$ and let $\{X_{r}\}_{0\leq r<1}$ be a bounded net in $H^{\infty
}(E)$ converging to $X$ in the ultraweak topology on $H^{\infty}(E)$. Then
\begin{align*}
\langle(W_{r}\times\tau)(X_{r})Sg,Sh\rangle &  =\langle S(V_{1r}\times
\rho)(X_{r})g,Sh\rangle=\langle(V_{1r}\times\rho)(X_{r})g,S^{\ast}Sh\rangle\\
&  \rightarrow\langle(V_{1}\times\rho)(X)g,S^{\ast}Sh\rangle=\langle
(W\times\tau)(X)Sg,Sh\rangle
\end{align*}
by Lemma \ref{Lemma1.7a}, the definition of $W\times\tau$ and equation
(\ref{intertwining}). Since the set $\{(W_{r}\times\tau)(X_{r})\}_{0\leq r<1}$
is uniformly bounded, the weak convergence of $(W_{r}\times\tau)(X_{r})$ to
$(W\times\tau)(X)$ follows.

A similar argument shows, finally, that $W\times\tau$ is ultraweakly
continuous. Indeed, for this it suffices to consider a \emph{bounded }net
$\{X_{\alpha}\}$ in $H^{\infty}(E)$ converging ultraweakly to the element $X$,
say. Then equation (\ref{intertwining}) shows that for all vectors $g$ and $h$
in $L_{\infty}(\frak{M})$ we have%
\begin{multline*}
\langle(W\times\tau)(X_{\alpha})Sg,Sh\rangle=\langle S(V_{1}\times\rho
_{1})(X_{\alpha})g,Sh\rangle\\
\rightarrow\langle S(V_{1}\times\rho_{1})(X)g,Sh\rangle=\langle(W\times
\tau)(X)Sg,Sh\rangle\text{.}%
\end{multline*}
Since the net $\{(W\times\tau)(X_{\alpha})\}$ and since vectors of the form
$Sg$ and $Sh$ are dense in $\frak{N}$ we conclude that $\{(W\times
\tau)(X_{\alpha})\}$ converges ultraweakly to $W\times\tau$.
\end{proof}

With these preliminaries out of the way, we may prove Theorem \ref{cncext}.

\begin{proof}
(of Theorem \ref{cncext}) Since $(T,\sigma)$ is assumed to be completely
non-coisometric, we may apply Lemma \ref{Lemma1.5} to conclude that
\[
K=L_{\infty}(\frak{Q}_{0})\oplus P_{\infty}(K)=L_{\infty}(\frak{Q}_{0}%
)\oplus\overline{P_{\infty}(L_{\infty}(\mathcal{D)})}\text{.}%
\]
The two summands in this decomposition are reducing for $V\times\rho$ and we
have the Wold decomposition of $(V,\rho)$ written as $(V,\rho)=(V_{ind}%
,\rho_{ind})\oplus(V_{\infty},\rho_{\infty})$. Note, in particular, that
$P_{\infty}$ commutes with $(V,\rho)$ and intertwines $(V,\rho)$ and the
restriction of $(V,\rho)$ to $L_{\infty}(\mathcal{D})$. Since this restriction
is an induced representation we may apply Lemma \ref{Lemma1.8a} to it. But
also, since the restriction of $(V,\rho)$ to $L_{\infty}(\frak{Q}_{0})$ is the
induced representation $(V_{ind},\rho_{ind})$, we may apply Lemma
\ref{Lemma1.7a} to it. Hence, using these two lemmas, we may conclude that
statements (1)-(3) of Theorem \ref{cncext} are valid for $(V,\rho)$. Finally,
since $(T,\sigma)$ is the compression of $(V,\rho)$ to $H$ and $H$ is a
semi-invariant subspace for $(V\times\rho)(H^{\infty}(E))$, we may define
$T\times\sigma(X)$, $X\in H^{\infty}(E)$, to be $P_{H}(V\times\rho)(X)|H$,
$X\in H^{\infty}(E)$, to complete the proof.
\end{proof}

In the Sz.-Nagy--Foia\c{s} theory \cite{szNF70}, a contraction $T$ is said to
belong to the class $C_{\cdot0}$ in case the powers of $T^{\ast}$ tend to zero
in the strong operator topology. Evidently, such a contraction is completely
non-coisometric. In \cite{gP89} Popescu generalized the notion of $C_{\cdot0}$
contractions to the setting of row contractions. He called the corresponding
class $C_{0}$. This idea makes sense in our setting, too.

\begin{definition}
\label{Cpointzero}A completely contractive covariant representation
$(T,\sigma)$ of $E$ on a Hilbert space $H$ is said to belong to the class
$C_{\cdot0}$ in case $\left\|  \tilde{T}_{k}^{\ast}h\right\|  \rightarrow0$
for all $h\in H$.
\end{definition}

Of course, such a representation is completely non-coisometric. Moreover, just
as in the case of single contractions and row contractions, one may identify
when a covariant representation belongs to $C_{\cdot0}$ in terms of its
minimal isometric dilation.

\begin{proposition}
Let $(T,\sigma)$ be a completely contractive covariant representation of $E$
on a Hilbert space $H$ and let $(V,\rho)$ be its minimal isometric dilation
acting on the Hilbert space $K$. Then $(T,\sigma)$ belongs to the class
$C_{\cdot0}$ if and only if $(V,\rho)$ belongs to $C_{\cdot0}$ and this
happens if and only if $(V,\rho)$ is an induced representation.
\end{proposition}

\begin{proof}
The operators $\tilde{T}_{k}\tilde{T}_{k}^{\ast}$ on $H$ form a decreasing
sequence, as do the operators $\tilde{V}_{k}\tilde{V}_{k}^{\ast}$ on $K$. The
representation $(T,\sigma)$ (resp. $(V,\rho)$) belongs to class $C_{\cdot0}$
if and only if $\tilde{T}_{k}\tilde{T}_{k}^{\ast}\rightarrow0$ strongly (resp.
$\tilde{V}_{k}\tilde{V}_{k}^{\ast}\rightarrow0$ strongly). In the case of
$(V,\rho)$, the limit of the $\tilde{V}_{k}\tilde{V}_{k}^{\ast}$ is just the
projection $P_{\infty}$ and so $(V,\rho)$ belongs to class $C_{\cdot0}$ if and
only if $P_{\infty}=0$. That is, $(V,\rho)$ belongs to class $C_{\cdot0}$ if
and only if $(V,\rho)$ is induced \cite[Corollary 2.10]{MS99}. Now if
$(V,\rho)$ belongs to class $C_{\cdot0}$, then since $P_{H}\tilde{V}_{k}%
\tilde{V}_{k}^{\ast}P_{H}=\tilde{T}_{k}\tilde{T}_{k}^{\ast}$ we see that
$(T,\sigma)$ belongs to class $C_{\cdot0}$. On the other hand, if $(T,\sigma)$
belongs to class $C_{\cdot0}$, this equation shows that $\left\|  \tilde
{V}_{k}^{\ast}h\right\|  \rightarrow0$ for all $h\in H$, which in turn means
that $H$ is orthogonal to the range of $P_{\infty}$. Since the range of
$P_{\infty}$ reduces $(V,\rho)$, the minimality of $(V,\rho)$ implies that
$P_{\infty}=0$.
\end{proof}

\begin{remark}
\label{ConcludingRemark}For those familiar with the Sz.-Nagy--Foia\c{s}
theory, the question, ``Is there an analogue of \emph{completely non-unitary}
contractions in the setting of completely contractive covariant
representations $(T,\sigma)$?'' is natural and compelling. This class of
contractions is larger than the class of completely non-coisometric
contractions. More important, however, is the fact that every contraction
decomposes into the direct sum of a unitary operator and a completely
non-unitary contraction. Further, the representation of the disc algebra
associated to a completely non-unitary contraction extends to a representation
of $H^{\infty}(\mathbb{T})$. This extension, of course, is the
Sz.-Nagy--Foia\c{s} functional calculus and the reason for its existence is
the fact that the spectral measure of the unitary part of the minimal
isometric dilation of a completely non-unitary contraction is mutually
absolutely continuous with respect to Lebesgue measure on the circle. Thus,
fancifully perhaps, we wonder if every covariant representation $(T,\sigma)$
of $E$ has a similar decomposition into the direct sum $(T,\sigma
)=(T_{cnu},\sigma_{cnu})\oplus(T_{c},\sigma_{c})$, where $T_{cnu}\times
\sigma_{cnu}$ extends to $H^{\infty}(E)$ in such a fashion that the
conclusions of Theorem \ref{cncext} hold, and where $T_{c}\times\sigma_{c}$ is
the restriction to $\mathcal{T}_{+}(E)$ of a $C^{\ast}$-representation of the
Cuntz-Pimsner algebra $\mathcal{O}(E)$, with $E$ treated simply as a $C^{\ast
}$-correspondence. This speculation is meaningful and of considerable
intersest even in the context of row contractions. We note that in the
profound study \cite{BVp03}, Ball and Vinnikov have formulated and analyzed a
notion of completely non-unitary row contractions. However, we do not know if
these have the extension property we desire. On the other hand, the beautiful
paper by Davidson, Katsoulis and Pitts \cite{DKP01} gives insights into
aspects of this problem that need to be addressed even in the context of
isometric row contractions, i.e., in the context of Cuntz-Toeplitz families of
isometries. Although they do not tackle our problem head on, they give an
analysis of the situations when a Cuntz-Toeplitz family of isometries
$(V_{1},V_{2},\cdots,V_{n})$ generates a weakly closed algebra that is
completely isometrically isomorphic, and ultraweakly homeomorphic to
$H^{\infty}(\mathbb{C}^{n})$, where $\mathbb{C}^{n}$ is viewed as a $W^{\ast}%
$-correspondence over $\mathbb{C}$. In particular they raise interesting
issues about what ought to be a concept of ``absolute continuity'' for
representations of the Cuntz algebra $\mathcal{O}_{n}$.
\end{remark}

\end{document}